\documentclass[11pt,oneside,reqno]{article}
\usepackage[width=15.7cm,height=24cm,centering]{geometry}
\usepackage{amsmath,amsthm,amssymb,color,bm}
\usepackage[authoryear]{natbib}
\usepackage{booktabs}

\usepackage{tikz}
\RequirePackage{amsthm,amsmath,amsfonts,amssymb,mathrsfs,dsfont,mathtools,thmtools}
\RequirePackage[colorlinks,citecolor=blue,urlcolor=blue,linkcolor=blue]{hyperref}
\usepackage{bookmark}
\usepackage[capitalize]{cleveref}
\RequirePackage{graphicx}

\crefname{equation}{}{}
\crefformat{equation}{\textup{(#2#1#3)}}
\crefrangeformat{equation}{\textup{(#3#1#4)--(#5#2#6)}}
\crefdefaultlabelformat{#2\textup{#1}#3}
\crefname{lemma}{Lemma}{Lemmas}
\crefname{page}{p.}{pp.}
\usepackage[normalem]{ulem}
\usepackage{bbm}
\usepackage{mathrsfs}
\usepackage{graphicx}
\usepackage{tikz}
\usepackage{enumitem}
\usetikzlibrary{arrows, automata}
\usetikzlibrary{calc}
\usetikzlibrary{positioning}

\numberwithin{equation}{section}
\allowdisplaybreaks[4]

\theoremstyle{plain}
\newtheorem{theorem}{Theorem}[section]
\newtheorem{proposition}{Proposition}[section]
\newtheorem{lemma}{Lemma}[section]
\newtheorem{corollary}{Corollary}[section]

\newtheorem{remark}{Remark}[section]

\theoremstyle{definition}

\makeatletter

\newcount\minute
\newcount\hour
\newcount\hourMins
\def\now{%
\minute=\time%
\hour=\time \divide \hour by 60%
\hourMins=\hour \multiply\hourMins by 60%
\advance\minute by -\hourMins%
\zeroPadTwo{\the\hour}:\zeroPadTwo{\the\minute}%
}
\def\zeroPadTwo#1{\ifnum #1<10 0\fi#1}

\renewcommand{\cite}{\citet}

\def\^#1{\ifmmode {\mathaccent"705E #1} \else {\accent94 #1} \fi}
\def\~#1{\ifmmode {\mathaccent"707E #1} \else {\accent"7E #1} \fi}

\def\*#1{#1^\ast}
\edef\-#1{\noexpand\ifmmode {\noexpand\bar{#1}} \noexpand\else \-#1\noexpand\fi}
\def\>#1{\vec{#1}}
\def\.#1{\dot{#1}}

\def\atop{\@@atop}
\def\*#1{\mathscr{#1}}

\renewcommand{\leq}{\leqslant}
\renewcommand{\le}{\leqslant}
\renewcommand{\geq}{\geqslant}
\renewcommand{\ge}{\geqslant}

\newcommand{\IE}{\mathbb{E}}

\newcommand{\Var}{\mathop{\mathrm{Var}}\nolimits}
\newcommand{\Cov}{\mathop{\mathrm{Cov}}}

\let\SNCLTStandardMid\mid
\def\mid{\vert}

\def\beqn#1\eeqn{\begin{align}#1\end{align}}
\def\beq#1\eeq{\begin{align*}#1\end{align*}}

\usepackage{graphicx}
\usepackage{latexsym}
\usepackage{amsmath,amsthm,amssymb,amscd}
\usepackage{epsf,amsmath}
\DeclareMathOperator*{\argmax}{argmax}

\def\E{{\IE}}

\renewcommand\section{\@startsection {section}{1}{\z@}%
{-2.8ex \@plus -.8ex \@minus -.2ex}%
{1.0ex \@plus.2ex}%
{\center\small\sc\mathversion{bold}}}
\def\subsection#1{\@startsection {subsection}{2}{0pt}%
{-2.8ex \@plus -.8ex \@minus -.2ex}%
{.8ex \@plus.2ex}%
{\bf\mathversion{bold}}{#1}}

\def\subsubsection#1{\@startsection{subsubsection}{3}{0pt}%
{\medskipamount}%
{-10pt}%
{\normalsize\itshape}{\texorpdfstring{\kern-2.2ex. #1.}{#1}}}

\def\blfootnote{\xdef\@thefnmark{}\@footnotetext}

\makeatother

\let\mid\SNCLTStandardMid

\providecommand{\Pp}{\mathbb{P}}
\providecommand{\R}{\mathbb{R}}

\providecommand{\1}{\mathbf{1}}

\providecommand{\distrect}{d_{\mathrm{rect}}}
\providecommand{\barPhi}{\overline\Phi}
\providecommand{\calR}{\mathcal R}

\providecommand{\Del}{\Delta_{n,q,K}}

\usepackage{array}
\defcitealias{LIULOPESROBUST}{Liu and Lopes (forthcoming)}

\hypersetup{pdfpagemode=UseOutlines,
  pdftitle={Sharp High-Dimensional Self-Normalized Central Limit Theorems},
  pdfauthor={Zi-Yao Su and Yi-Kun Zhao},
  pdfsubject={Self-normalized Gaussian approximation under marginal moments},
  pdfkeywords={self-normalized sums, Gaussian approximation, correlation gaps}}
\bookmarksetup{depth=3,open=true,openlevel=2}
\begin{document}

\title{Sharp High-Dimensional\\Self-Normalized Central Limit Theorems}
\author{Zi-Yao Su, Yi-Kun Zhao}
\date{\it The Chinese University of Hong Kong}
\maketitle

\noindent{\bf Abstract:}
We study Gaussian approximation over all coordinate rectangles for
coordinatewise self-normalized sums of iid centered random vectors. The Gaussian
reference has the population correlation matrix, which may be singular.
Under a fixed bound $K>1$ on the standardized marginal third absolute
moments, uniform approximation over the corresponding moment class holds
exactly when $d=o(\sqrt n)$ under arbitrary coordinate dependence. If
$R_{jk}\le1-a$ for $j\ne k$ and a fixed $a\in(0,1]$, the necessary and
sufficient condition becomes $\log d=o(n^{1/3})$.

For marginal moment order $q\ge3$, we also allow the correlation gap
$a_n$ to depend on $n$. For every fixed $0<\beta<1/3$, uniform
approximation over the corresponding class in dimension
$d_n=\lfloor e^{n^\beta}\rfloor$ holds if and only if
$a_n\log n\to\infty$. Under arbitrary dependence, we determine the
sharp error order $1\wedge\{d/n^{(q-2)/2}\}^{1/q}$ for every
$q\ge3$ and $d\ge2$. At a fixed positive gap, we identify, for every $q>2$,
the supremum of the decay exponents in uniform error bounds that allow
fixed logarithmic factors. The proofs replace observations one at a time,
keeping their contributions to both numerator and denominator in the
same expansion. The resulting weights depend only on the retained
sample, allowing the remainders to be controlled by marginal moments.



\section{Introduction}\label{sec:introduction}

Self-normalization gives accurate normal tail approximations under weak
moment assumptions. We ask how far this accuracy extends jointly across
many coordinates. Even when the Gaussian reference matches the population
correlation matrix, dependence can prevent approximation in dimensions
far below the exponential range suggested by the marginal tails.
We determine how restrictions on that dependence change the dimension
range for the coordinatewise self-normalized statistic.

Let $X_1,\ldots,X_n$ be independent copies of a centered vector in
$\R^d$, $d\ge2$, with unit marginal variances, and set
\[
 T_{n,j}=\frac{\sum_{i=1}^nX_{ij}}
                    {(\sum_{i=1}^nX_{ij}^2)^{1/2}},\qquad 1\le j\le d,
 \qquad 0/0=0.
\]
We compare $T_n$ with $G_P\sim N_d(0,R_P)$, where $R_P$ is the
population correlation matrix. The error is measured over all coordinate
rectangles, uniformly over laws with $\max_j\E|X_j|^q\le K$ for
fixed $q>2$ and $K>1$. The law may depend on $n$, and $R_P$ may be
singular. We keep $T_n$ fixed when studying the effect of dependence;
approximations based on modified observations are discussed in
Section~\ref{sec:related-work}.

The main contrast is already present under a bounded third absolute
moment. Without a restriction on coordinate dependence, uniform
approximation requires polynomial dimension. A fixed positive one-sided
correlation gap, $R_{jk}\le1-a$ for $j\ne k$, changes the condition
to exponential dimension:
\[
 \begin{array}{ll}
 d=o(\sqrt n),&\text{under arbitrary coordinate dependence},\\[2pt]
 \log d=o(n^{1/3}),&\text{under a fixed positive one-sided gap}.
 \end{array}
\]
Each condition is necessary and sufficient over its respective class
when $q=3$; see Theorems~\ref{thm:main-polynomial} and
\ref{thm:main-shrinking}. Thus the change in dimension comes from
the dependence restriction, with the same statistic and marginal
moment bound. The gap bounds correlations from above; it permits
correlations close to $-1$ and imposes no dimension-independent
spectral bounds.

The next question is how much the gap can decrease while preserving
exponential dimension. For every $q\ge3$, fixed $0<\beta<1/3$ and
deterministic sequence $a_n\in[0,1]$, uniform approximation over the
corresponding moment and correlation class in dimension
$d_n=\lfloor e^{n^\beta}\rfloor$ holds if and only if
$a_n\log n\to\infty$. This includes fixed positive gaps and
identifies the transition when the correlations are allowed to
approach one.

The occurrence of $\log n$ can be understood from the lower-bound
construction in Proposition~\ref{prop:rare-mark-unified}. A marginal
$q$th-moment bound permits positive observations of size
$A\asymp\sqrt n\log n$ with probability of order $A^{-q}$.
Taking $m\asymp n^{(q-2)/2}(\log n)^{q+1}$ makes the expected
number of such observations among the $nm$ entries of the block
of order $\log n$. In this construction, with probability tending
to one, at least one occurs and its coordinate has a self-normalized
ratio above a threshold tending to one from below.
The rectangle need constrain only this polynomial
block, even when the ambient dimension is exponential.

The relevant Gaussian scale is visible in an equicorrelated block:
\[
 \begin{aligned}
 G_j&=\sqrt{1-a}\,Z_0+\sqrt a\,Z_j,\qquad 1\le j\le m,\\
 \max_{j\le m}G_j&=\sqrt{1-a}\,Z_0+\sqrt a\max_{j\le m}Z_j,
 \end{aligned}
\]
where $Z_0,Z_1,\ldots,Z_m$ are independent standard normal variables.
The variance correction caused by the rare observations is
$o(1/\log n)$, so its effect on the Gaussian maximum tends to zero
in probability.
In the displayed model, the independent part
shifts the maximum by order $\sqrt{a\log m}\asymp\sqrt{a\log n}$,
whereas the common component retains fluctuations of order one as
$a\to0$. If $a\log n$ stays bounded, the Gaussian block maximum
has a nonvanishing probability of lying below the threshold, while
the sample block maximum exceeds it with probability tending to one.
This explains the obstruction in the lower bounds. The sufficient
condition follows from comparison estimates uniform over the full
correlation class.

We also obtain quantitative refinements of these dimension results.
For every $q\ge3$, the sharp order of the uniform error under
arbitrary dependence is
\[
 1\wedge\left(\frac d{n^{(q-2)/2}}\right)^{1/q}.
\]
For $q>3$ and gaps satisfying $a_n\to0$ and
$a_n\log n\to\infty$, we determine the critical power of the gap
in dimension bounds near $\log d\asymp n^{1/3}$; the sufficient
and necessary bounds still differ by logarithmic factors, as stated
in Section~\ref{sec:main}. For each $q>2$ and fixed positive gap, we
also determine the supremum of the sample-size decay exponents in uniform
error bounds allowing fixed logarithmic factors. Section~\ref{sec:main-sample-rate}
gives the precise exponent and explains when it is attained.

\subsection{Related work}\label{sec:related-work}

\paragraph{The original self-normalized statistic.}
A closely related result under unrestricted coordinate dependence is
due to \citet[Lemma~A.7 and Theorem~A.8]{KP2026FREE}. For the ordinary
Studentized mean, they prove uniform approximation over all rectangles
under iid sampling, a marginal $q>4$ moment bound and a positive
marginal variance lower bound, allowing arbitrary, possibly singular,
population correlation matrices. Their sufficient range includes
$d=o(n^{(q-2)/2-\xi})$ for fixed $\xi\in(0,1)$.
The exact transformation between their centered sample-variance
denominator and our self-normalized statistic is given in
Section~\ref{sec:inference-scope}. For $T_n$, we obtain the necessary
and sufficient dimension condition for every $q\ge3$, together with
the sharp error order and no fixed loss in the power of $n$. We then
identify how a one-sided correlation gap enlarges this range.

Results under more specific dependence assumptions help locate the
source of the difficulty. When the coordinates within each observation
are iid, \citet{DAS2024} proves approximation over all rectangles;
under third moments, the rate is of order
$n^{-1/2}\log^{3/2}(dn)$ in its stated range, yielding
$\log d=o(n^{1/3})$. \citet{QCS2025} study relative tail errors for
maxima under uncorrelated coordinates or additional correlation
conditions. Their setting permits dependent, uncorrelated coordinates,
and their relative tail criterion differs from the uniform rectangle
criterion used here. \citet{CTUK2025} retain the random denominators
and treat the maximum absolute self-normalized sum under general
dependence, using either an optimized Gaussian reference or the
population correlation matrix. Their simple moment bounds involve
$\E\max_j|X_j|^q$, which can be as large as $Kd$ when only
$\max_j\E|X_j|^q\le K$ is assumed.

For unrestricted dependence, \citet{WU2026} gives counterexamples to
the conjecture in \citet[Remark~2.4]{QCS2025}, using the original
self-normalized statistic and its matching population covariance.
One example has dimension of order $n^q$ under a fixed marginal
$q$th moment bound, $q>2$; another has all fixed marginal moments
bounded and error tending to one when $\log d$ is of order
$(\log n)^3$. The paper also constructs a fixed infinite-dimensional
distribution. Our lower bounds locate the obstruction at the smaller
scale $n^{(q-2)/2}$, match the uniform error order for $q\ge3$, and
quantify its dependence on the correlation gap.
Proposition~\ref{prop:qcs-counterexample} also gives a fixed-threshold
example in dimension $\lfloor n^{\tau/2}\rfloor$ under a marginal
$(2+\tau)$th moment bound, $0<\tau\le1$.

\paragraph{Changing the statistic.}
Winsorization and trimming offer a different way to accommodate heavy
tails. \citet[Theorems~3.1 and~4.1]{KOCK2026BARRIER} establishes a
dimension barrier for this broader problem. Over the mean-zero
class with marginal variances bounded below and marginal $q$th
absolute moments bounded above, $q>2$, he proves that
\[
 \log d=o\!\left(n^{(q-2)/(3q-2)}\right)
\]
is necessary for uniform Gaussian approximation, even when the
statistic can be chosen or randomized. The Gaussian reference has
the population covariance or correlation matrix, and the procedure
cannot depend on the unknown law.
For the population covariance matrix, he establishes approximation
over all rectangles for quantile-winsorized means under this same
condition when $2<q\le4$. Thus the sufficient and necessary dimension
conditions match in this moment range. For $q>4$, the sufficient
condition has the additional logarithmic factor
\[
 \log d=o\!\left(n^{(q-2)/(3q-2)}/(\log n)^{q/(3q-2)}\right).
\]
Our results address the complementary question of retaining $T_n$:
we determine its dimension limit over the unrestricted moment class,
and the larger range obtained by restricting the correlations.

Other robust results provide several forms of calibration.
\citet{LOUWU2017} use clipping at a moment-dependent level and obtain
Gaussian approximation for the absolute maximum of clipped means
when $\log d=o(n^{\theta/(4+3\theta)})$ under a bounded marginal
$(2+\theta)$th moment, $0<\theta\le1$, without restricting coordinate
dependence. \citet{KP2025ROBUST} use empirical quantiles to winsorize
the observations without knowledge of the moment order or bound.
Under iid sampling and a fixed marginal $q$th moment bound, $q>2$,
their sufficient condition is
$\log d=o(n^{(q-2)/(5q-2)})$; their results include all rectangles,
full studentization and Gaussian multiplier bootstrap.
\citet{RESENDE2024} obtains the same dimension exponent for a
supremum of trimmed means, with empirical bootstrap and simultaneous
confidence intervals formed by dividing trimmed means by estimated
marginal scales \citep[Section~4.1]{RESENDE2024}.
Table~\ref{tab:third-moment-comparison} compares the dimension
conditions when the marginal third absolute moments are bounded.

Faster bootstrap rates can be obtained with additional structure.
\citetalias{LIULOPESROBUST} use moment equivalence in all directions,
variance decay, density and correlation assumptions, together with
partial standardization. For tests, \citet[Theorems~3.1 and~4.1]{KP2026FREE}
establish equal asymptotic power for robust and ordinary max-tests
under common conditions where both are valid.
Section~\ref{sec:inference-scope} also records the conservative
calibration of $T_n$ available from marginal tail bounds.

\begin{table}[htbp]
\centering
\small
\caption{Dimension conditions under iid sampling and bounded marginal
third absolute moments, with the positive variance bounds of the cited
results. All observations are uncontaminated. The conditions for $T_n$
are necessary and sufficient uniformly over the corresponding classes.
The robust-estimator rows give sufficient conditions; for Kock's
approximation with population covariance, the sufficient condition
matches the necessary dimension condition.}
\label{tab:third-moment-comparison}
\setlength{\tabcolsep}{5pt}
\begin{tabular}{@{}>{\raggedright\arraybackslash}p{3.9cm}>{\raggedright\arraybackslash}p{6.5cm}>{\raggedright\arraybackslash}p{3.9cm}@{}}
\toprule
Statistic and reference & Approximation and assumptions & Dimension condition\\
\midrule
Original $T_n$, this paper
& Arbitrary coordinate dependence; all rectangles; population correlation
& $d=o(\sqrt n)$\\[4pt]
Original $T_n$, this paper
& Fixed positive one-sided correlation gap; all rectangles; population correlation
& $\log d=o(n^{1/3})$\\[4pt]
\midrule
Quantile-winsorized means \citep{KOCK2026BARRIER}
& Arbitrary coordinate dependence; all rectangles; population covariance
& $\log d=o(n^{1/7})$\\[4pt]
Clipped means \citep{LOUWU2017}
& Arbitrary coordinate dependence; absolute maximum; moment-dependent cutoff
& $\log d=o(n^{1/7})$\\[4pt]
Winsorized means \citep{KP2025ROBUST}
& Arbitrary coordinate dependence; all rectangles; full studentization and bootstrap
& $\log d=o(n^{1/13})$\\[4pt]
Trimmed means \citep{RESENDE2024}
& Arbitrary coordinate dependence; supremum; empirical bootstrap and studentized trimmed-mean intervals
& $\log d=o(n^{1/13})$\\
\bottomrule
\end{tabular}
\end{table}

\paragraph{Scalar estimates and comparison methods.}
The scalar theory describes the accuracy gained from
self-normalization. \citet{BG1996} give a Berry--Esseen bound for
Student's statistic, and \citet{JSW2003} establish Cram\'er-type
deviations under finite $(2+\delta)$th moments. Refined expansions
identify the skewness correction \citep{WANG2005,SG2017}.
For self-normalized processes and general numerator--denominator
pairs, see \citet{SZ2016,GSS2022} and the broader treatment in
\citet{SW2013}.

For ordinary means, marginal moment bounds already impose a
polynomial restriction on dimension \citep{ZW2017,KP2024}.
In the common range $q\ge4$, \citet{KP2024} identify the power
$q/2-1$, with fixed power margins in the sufficient and failure
conditions. Their obstruction can occur with independent heavy-tailed
coordinates. Self-normalization permits exponential dimension under
coordinate independence, so the lower bounds here concern the
interaction between the random denominators and coordinate dependence.

The comparison theory for ordinary sums provides Gaussian and bootstrap
approximations allowing singular covariance matrices
\citep{CCK2013,CCK2017,CCKK2022,FK2021}. Refined bounds describe how
accuracy depends on covariance degeneracy and distributional structure
\citep{LOPES2022,CCK2023,AL2026,FKLZ2023,FK2024}.
For general covariance, \citet{LM2026,CCKK2026} obtain near-$n^{-1/2}$
rates under their respective light-tail assumptions.

To apply these methods with only marginal polynomial moments, we
must control large observations while retaining their contributions
to the random denominators. We use the bounded-sum comparison of
\citet{CCKK2026} and adapt the two-stage, third-moment comparison of
\citet{GZ2026} to numerator--denominator pairs.
The pairwise gap in \citet[Corollary~4]{CCKK2026} is
$1-\max_{j\ne k}|R_{jk}|$, which also controls nearly opposite
coordinates. Our one-sided restriction allows such coordinates and
requires separate estimates for them. Section~\ref{sec:proof-outline}
explains how the comparison is combined with the denominator and
Gaussian boundary estimates.

\subsection{Proof outline}\label{sec:proof-outline}

\paragraph{Retaining the random denominators.}
We compare smooth rectangle tests by replacing one independent
observation at a time, using the Lindeberg method of \citet{CHAT2006};
its application to self-normalized sums is also proposed in
\citet[Section~3]{CTUK2025}. Two simpler reductions incur dimension
costs. Excluding every sample entry larger than $\tau$ gives the union
bound $ndK\tau^{-q}$, which at $\tau\asymp\sqrt n$ tends to zero
only in polynomial dimension. Replacing the denominators by their
population limits is also costly: even when the maximal normalized
numerator is $O(\sqrt{\log(ed)})$, a uniform error $\eta\le1/2$ in
the empirical second moments costs $O(\eta\log(ed))$ after Gaussian
anti-concentration. Obtaining a sufficiently small $\eta$ is restrictive
when the squared observations have only $q/2$ moments.

We keep the contribution of the observation to both numerator and
denominator. After deleting one observation, let $S$ and $Q>0$ be the
sum and sum of squares in one coordinate, and put $t=S/\sqrt Q$.
Inserting $y$ gives
\[
 \frac{S+y}{\sqrt{Q+y^2}}
 =\frac{t+y/\sqrt Q}{\sqrt{1+y^2/Q}}.
\]
The quadratic term includes the change in the denominator.
Conditional on the retained sample, the Taylor coefficients are fixed.
The means and covariances of the original and Gaussian incoming
observations agree, so the linear and quadratic terms have equal
conditional expectations. The inserted ratio has squared positive part
at most $t_+^2+1$ for every $y$. We use a smooth increasing
approximation to the squared positive-part function whose increase
after insertion is also bounded by an absolute constant.

Differentiating each smooth maximum gives nonnegative coordinate
weights $w_j$ with $\sum_jw_j=1$. Evaluating them before insertion
makes them depend only on the retained sample and the fixed rectangle
thresholds. Combined with the insertion bound, the weighted derivative
estimates control the Taylor remainder through weighted $q$th powers
when $2<q\le3$; see Lemma~\ref{lem:core-weighted}. For an independent
incoming observation $Y$ from the original sample,
\[
 \E\!\left[\sum_j w_j|Y_j|^q\,\middle|\,
                 \text{retained sample}\right]\le K,
 \qquad \sum_jw_j=1.
\]
Thus the conditional expectation uses only marginal moments.
Gaussian observations satisfy the same bound with a constant depending
on $q$. Independence and the moment bounds persist through the
replacement, so this construction applies at every step and yields
Proposition~\ref{prop:core-comparison}.

For the sharp error bound under arbitrary dependence and $q\ge3$,
we balance the smoothing error with the contribution of large
observations. This yields the same order as the lower bound.
Section~\ref{sec:unrestricted} gives the quantitative estimates and
the choice of smoothing width.

\paragraph{Controlling the Gaussian boundary error.}
Passing from a smooth test back to a rectangle requires two Gaussian
estimates. The first bounds the probability of a thin boundary strip.
For inequalities $G_j\le b_j$, we measure width using the increasing
transformation $g_N(t)=-\log\barPhi(t)$.
Lemma~\ref{lem:score-strip} bounds the probability that
$\max_j\{g_N(G_j)-g_N(b_j)\}$ lies in an interval by the interval's
width. This estimate uses the conditional monotonicity argument of
\citet{CCK2015} and holds for arbitrary covariance, including singular
covariance; it does not require a correlation gap.

The second estimate controls how many coordinates can be close to
their own thresholds at the same time. Here the correlation gap is
essential. We condition the Gaussian law on the full rectangle and
use the entropy representation of \citet{FOLLMER1985,LEHEC2013}.
This gives a nonnegative measurable function on the underlying
Gaussian space, determined by the conditioning rectangle.
Lemma~\ref{lem:core-energy} bounds these counts pointwise in terms
of this same function, for every strip width and every subset satisfying
the gap condition. Its integral bound controls, by Markov's
inequality, the probability that the function is large. Together,
these estimates bound the smoothing error without separate union
bounds over the subsets and widths. The drift representation of
the function is given with the lemma.

\paragraph{Using moments above order three.}
For $q>3$, the cubic term must be treated before taking absolute
values. We use the two-stage comparison of \citet[Section~4]{GZ2026}
and extend the matching to the moments involving the denominator.
In the first stage, each observation is replaced by an auxiliary
numerator--denominator pair. The numerator matches the joint third
moments of the original observation; the denominator contribution is
chosen so that the mixed moments in the expansion also agree. The
expected cubic terms therefore cancel. In the second stage, the auxiliary sample is compared
with the Gaussian reference. The cubic contribution of the original
skewness remains and must be bounded.

An independent Gaussian component in the auxiliary observations
permits integration by parts in the second stage without inverting
the covariance matrix. A weighted estimate adapted from
\citet[Section~H.2]{CCKKSUP2026} then controls the contribution of
the truncated third moments. The auxiliary construction, the mixed
moment identities and the remaining tail estimates are developed in
Section~\ref{sec:shrinking-proofs}.
The lower bound shows why the skewness constraint remains. In the
bounded, skewed independent-coordinate example, the scalar correction
has order $x^3/\sqrt n$ \citep{WANG2005,GSS2022}. At the threshold
$x\asymp\sqrt{\log d}$, it explains the boundary
$\log d\asymp n^{1/3}$ even when higher moments exist.

Section~\ref{sec:main} states the dimension and rate results.
Sections~\ref{sec:proofs-main} and \ref{sec:shrinking-proofs} prove
them, using the comparison and Gaussian boundary estimates of
Section~\ref{sec:inputs}. Section~\ref{sec:extensions} first treats
partitions with a correlation gap within each group, then the
additional structure of bounded-variation factor models.
Appendix~\ref{sec:inference-details} supplies the short derivations
behind the inference remarks in Section~\ref{sec:inference-scope}.

\section{Main results}\label{sec:main}

\subsection{The model}

For $q>2$, $K>1$, $d\ge2$ and $a\in[0,1]$, let
$\mathcal P_{d,q}(K,a)$ be the class of laws satisfying
\begin{equation}\label{eq:model}
 \E X_j=0,\qquad \E X_j^2=1,\qquad
 \max_{j\le d}\E|X_j|^q\le K,\qquad
 R_{jk}:=\E X_jX_k\le1-a\quad(j\ne k).
\end{equation}
For independent copies $X_1,\ldots,X_n$ with common law $P$, set
\begin{equation}\label{eq:statistic}
 T_{n,j}=\frac{\sum_{i=1}^nX_{ij}}
                 {(\sum_{i=1}^nX_{ij}^2)^{1/2}},\qquad
 G_P\sim N_d(0,R_P),\qquad 0/0=0.
\end{equation}
Let $\mathcal R_d$ consist of all closed coordinate rectangles, with
infinite endpoints and zero widths allowed. Define
\begin{equation}\label{eq:rectangle-error}
 \begin{split}
 \distrect(U,V)&=\sup_{A\in\mathcal R_d}
       |\Pp(U\in A)-\Pp(V\in A)|,\\
 \Delta_{n,q,K}(d,a)&=\sup_{P\in\mathcal P_{d,q}(K,a)}
                         \distrect(T_n,G_P).
 \end{split}
\end{equation}
The law in the supremum may depend on $n$. All uniform statements below
use this convention. The restriction on $R_{jk}$ is one-sided; negative
correlations and singular covariance matrices are allowed. It imposes
no dimension-independent spectral bounds. For example, a matrix with
all off-diagonal entries equal to $1-a$, $0<a<1$, has largest eigenvalue
$1+(d-1)(1-a)$. In two dimensions, the correlation $-1+\varepsilon$
is compatible with a fixed positive one-sided gap, while the smallest
eigenvalue is $\varepsilon$.

We write $\Phi$ and $\phi$ for the standard normal distribution function
and density, and $\barPhi=1-\Phi$. The symbols $c,C$ denote positive
constants whose values may change. Their dependence is restricted to the
fixed parameters specified in the relevant statement. We write
$x_+=\max(x,0)$ and $\log_+x=\max(\log x,0)$.

The unit-variance assumption only fixes scale. For marginal variances
$0<\sigma_j^2<\infty$, replace the moment bound by
$\max_j\E|X_j/\sigma_j|^q\le K$ and let $R$ be the correlation matrix.
Dividing a column by $\sigma_j$ leaves its ratio unchanged. Zero-variance
coordinates may be removed. The restriction $K>1$ is needed for the
general lower bounds: at $K=1$, unit variance forces $|X_j|=1$ almost
surely, and the denominator is deterministic.

\subsection{From polynomial to exponential dimension}

\begin{theorem}[Unrestricted coordinate dependence]\label{thm:main-polynomial}
Fix $q\ge3$ and $K>1$. There are constants $c,C>0$ such that, for all
sufficiently large $n$ and every $d\ge2$,
\begin{equation}\label{eq:main-polynomial}
 c\left[1\wedge\left(\frac d{n^{(q-2)/2}}\right)^{1/q}\right]
 \le\Delta_{n,q,K}(d,0)
 \le C\left[1\wedge\left(\frac d{n^{(q-2)/2}}\right)^{1/q}\right].
\end{equation}
The upper bound holds for every $a\in[0,1]$. More generally, if
$a_n\log n=O(1)$, then
\begin{equation}\label{eq:main-polynomial-criterion}
 \Delta_{n,q,K}(d_n,a_n)\to0
 \quad\Longleftrightarrow\quad d_n=o(n^{(q-2)/2}).
\end{equation}
\end{theorem}

The two-sided bound in \eqref{eq:main-polynomial} has no logarithmic
loss. It concerns the class with $a=0$. For positive gaps satisfying
$a_n\log n=O(1)$, \eqref{eq:main-polynomial-criterion} identifies the
dimension condition for convergence; it does not assert the same
two-sided error bound. Theorem~\ref{thm:polynomial-envelope} gives the
proof and the lower bound used in this distinction.

\begin{theorem}[Correlation gaps and dimension]\label{thm:main-shrinking}
\label{cor:dimensions}
Fix $q\ge3$, $K>1$ and $0<\beta<1/3$. For every deterministic sequence
$a_n\in[0,1]$,
\begin{equation}\label{eq:main-exponential-criterion}
 \Delta_{n,q,K}(\lfloor e^{n^\beta}\rfloor,a_n)\to0
 \quad\Longleftrightarrow\quad a_n\log n\to\infty.
\end{equation}
If $\liminf_n a_n>0$, then
\begin{equation}\label{eq:main-cubic-criterion}
 \Delta_{n,q,K}(d_n,a_n)\to0
 \quad\Longleftrightarrow\quad\log d_n=o(n^{1/3}).
\end{equation}
If $a_n\to0$ and $a_n\log n\to\infty$, then, for every fixed $c_0,C_0>0$,
\begin{equation}\label{eq:main-leading-gap}
 \log\Delta_{n,q,K}(d,a_n)
 =-\left(\frac{q-2}{2}+o(1)\right)a_n\log n
\end{equation}
uniformly over $c_0n^{(q-2)/2}\le d\le\exp(C_0n^\beta)$.
\end{theorem}

Equation~\eqref{eq:main-exponential-criterion} identifies the scale
$1/\log n$ at which the dimension condition changes. The last assertion
also gives the leading term of the logarithmic error when the gap
vanishes. Equivalently, in its stated dimension range,
$\Delta_{n,q,K}(d,a_n)=n^{-\{(q-2)/2+o(1)\}a_n}$. The $o(1)$ is
inside the exponent; it does not imply that the remaining factor is
bounded by a fixed power of $\log n$.

\subsection{Dimensions near \texorpdfstring{$\log d\asymp n^{1/3}$}{log d of order n\string^(1/3)}}

The preceding theorem treats $\log d=n^\beta$ for fixed $\beta<1/3$.
For $q>3$, the next result describes how a vanishing gap changes the
dimension condition when $\log d$ is closer to $n^{1/3}$.

\begin{theorem}[The critical power of the correlation gap]\label{thm:main-critical}
Fix $q>3$ and $K>1$. Suppose $a_n\to0$ and $a_n\log n\to\infty$.
For a fixed real number $b$, let
$\log d_n\asymp n^{1/3}a_n^b$. Then
\begin{equation}\label{eq:main-critical-power}
 \Delta_{n,q,K}(d_n,a_n)\to0
 \quad\Longleftrightarrow\quad b\ge\frac4{3(q-2)}.
\end{equation}
If $b<4/[3(q-2)]$, the error tends to one.
\end{theorem}

The theorem determines the power of $a_n$, including its critical
value. The full dimension threshold is bounded as follows.
Theorem~\ref{thm:shr-critical-upper} proves convergence if
\begin{equation}\label{eq:main-log-upper}
 \log(2d_n)=o\!\left[
 n^{1/3}a_n^{4/[3(q-2)]}
 \left\{\frac{\log(1/a_n)}{\log\log(1/a_n)}\right\}^{
           2\min\{1,q-3\}/[3(q-2)]}\right].
\end{equation}
Theorem~\ref{thm:newlower-sphere} shows that convergence requires
\begin{equation}\label{eq:main-log-lower}
 \log(2d_n)=o\!\left[
 n^{1/3}a_n^{4/[3(q-2)]}
 \{\log(1/a_n)\}^{\min\{1,\,2(q-1)/[3(q-2)]\}}\right].
\end{equation}
These two bounds imply Theorem~\ref{thm:main-critical}, including
the boundary case $b=4/[3(q-2)]$. The logarithmic interval between
\eqref{eq:main-log-upper} and \eqref{eq:main-log-lower} remains open.
Thus the sharp assertion here concerns the power of the gap.

\subsection{The optimal sample-size exponent}\label{sec:main-sample-rate}

We now fix the correlation gap and ask for the largest possible power
of $n$ in a uniform error bound. We permit a fixed power of
$\log(end)$ in that bound. This criterion concerns the sample-size
exponent; the resulting dimension dependence need not yield the full
range in \eqref{eq:main-cubic-criterion}.

\begin{theorem}\label{thm:main-fixed}
Fix $q>2$, $K>1$ and $a\in(0,1]$. The supremum of all exponents $\gamma$
for which
\begin{equation}\label{eq:main-sample-class}
 \Delta_{n,q,K}(d,a)\le Cn^{-\gamma}\log^B(end),\qquad n,d\ge2,
\end{equation}
holds with finite constants $C,B$, independent of $n,d$, is
\begin{equation}\label{eq:main-sample-exponent}
 \gamma_*(q,a)=\min\left\{\frac{q-2}{2q},
                         \frac{(q-2)a}{2(1-a)}\right\},
\end{equation}
where the second term is $+\infty$ at $a=1$. If $a>1/(q+1)$, the
supremum is attained with $B=2$. If $0<a\le1/(q+1)$, every smaller
positive exponent is attained with $B=2$, but the endpoint is not
attained for any finite $B$.
\end{theorem}

The two terms in \eqref{eq:main-sample-exponent} have different sources.
A lower bound based on nearly opposite coordinates limits the first
exponent. Rare perturbations of a common component limit the second.
Section~\ref{sec:proofs-main} gives both constructions and the
corresponding upper bounds.

The endpoint has an explicit correction. For $0<a\le1/(q+1)$ there are
integers $d_n\asymp n^{(q-2)/2}(\log n)^{q+1}$ such that
\begin{equation}\label{eq:main-endpoint}
 \log\Delta_{n,q,K}(d_n,a)
 =-\frac{(q-2)a}{2(1-a)}\log n
   +\frac{\sqrt{a(q-2)}}{1-a}\sqrt{\log n}
   +O(\log\log n).
\end{equation}
The factor $\exp\{\sqrt{a(q-2)\log n}/(1-a)\}$ exceeds every fixed
power of $\log n$ and prevents attainment of the endpoint. The full bounds and
proof are given in Theorem~\ref{thm:fixed-sample-rate}. This result allows
every $q>2$; the cubic dimension threshold in
Theorem~\ref{thm:main-shrinking} assumes $q\ge3$.

\subsection{Failure of a dimension-uniform tail approximation}

In \citet[Remark~2.4, (2.16)--(2.17)]{QCS2025}, the conjectured
relative tail approximation is uniform in dimension under a bounded
marginal $(2+\tau)$th absolute moment, with no restriction on the
correlation matrix. The Gaussian reference has that same matrix.
\citet{WU2026} gives counterexamples for self-normalized and
Studentized maxima, including a fixed infinite-dimensional population
with bounded marginal moments of every fixed order and positive
definite finite-dimensional covariance matrices.
Under a third-moment bound, \citet[Remark~4.1]{KOCK2026BARRIER} also
deduces failure of the uniform conjecture from his lower bound for
arbitrary statistics. The next proposition
gives failure at fixed thresholds in dimension
$d_n=\lfloor n^{\tau/2}\rfloor$, under a marginal $(2+\tau)$th
moment bound. The counterexamples address this conjecture; the
proved results of \citet{QCS2025} use additional assumptions.

\begin{proposition}\label{prop:qcs-counterexample}
For every $0<\tau\le1$ and $K>1$, there are centered laws with unit marginal variances
in dimension $d_n=\lfloor n^{\tau/2}\rfloor$, with marginal
$(2+\tau)$th moments at most $K$, for which
\begin{align}
 \lim_{n\to\infty}
 \frac{\Pp(\max_{j\le d_n}T_{n,j}>0)}
      {\Pp(\max_{j\le d_n}G_{P,j}>0)}&>1,
       \label{eq:qcs-positive}\\
 \lim_{n\to\infty}
 \frac{\Pp(\max_{j\le d_n}|T_{n,j}|>1/2)}
      {\Pp(\max_{j\le d_n}|G_{P,j}|>1/2)}&>1.
       \label{eq:qcs-absolute}
\end{align}
Every finite-dimensional population covariance in this construction is
positive definite.
\end{proposition}

The construction uses rare positive observations added to a common
symmetric sign. The Gaussian reference has positive definite covariance,
so singularity is not responsible for the failure. The proof is given
in Section~\ref{sec:qcs-proof}.

\subsection{Interpretation for mean inference}\label{sec:inference-scope}

The following two remarks explain the inferential scope of the
approximation results. Their derivations are given in
Appendix~\ref{sec:inference-details}.

\paragraph{Studentization.}
Let $\bar X_j=n^{-1}\sum_iX_{ij}$,
$s_{n,j}^2=n^{-1}\sum_i(X_{ij}-\bar X_j)^2$, and
$U_{n,j}=\sqrt n\,\bar X_j/s_{n,j}$, with $U_{n,j}=0$ when
$s_{n,j}=0$. On $E_n=\{\min_j s_{n,j}>0\}$,
\begin{equation}\label{eq:studentization-transform}
 U_{n,j}=\frac{T_{n,j}}{\sqrt{1-T_{n,j}^2/n}},\qquad
 T_{n,j}=k_n(U_{n,j}),\qquad
 k_n(x)=\frac{x}{\sqrt{1+x^2/n}}.
\end{equation}
For $G\sim N_d(0,R_P)$, including singular $R_P$,
\begin{equation}\label{eq:studentization-transfer}
 \distrect(T_n,G)\le\distrect(U_n,G)+\Pp(E_n^c)
                   +C\frac{\log^2(edn)}n.
\end{equation}
The uniform variance bounds of \citet{KP2026FREE} make
$\Pp(E_n^c)=o(1)$ in their polynomial dimension range, so their
sufficient condition transfers to $T_n$ through this inequality.

\paragraph{Conservative calibration from marginal tails.}
For $q>2$, put $\delta=\min(q-2,1)$ and, for $0<\alpha_n<1$, set
$c_n=\Phi^{-1}(1-\alpha_n/(2d_n))$.
The scalar moderate-deviation bound of
\citet[Theorem~2.3]{JSW2003} and a union bound imply
\begin{equation}\label{eq:marginal-calibration}
 \begin{split}
 \log(d_n/\alpha_n)&=o\!\left(n^{\delta/(2+\delta)}\right)
 \\
 &\Longrightarrow\quad
 \sup_{P\in\mathcal P_{d_n,q}(K,0)}
 \Pp\!\left(\max_{j\le d_n}|T_{n,j}|>c_n\right)
 \le\alpha_n\{1+o(1)\}.
 \end{split}
\end{equation}
For $q\ge3$ and fixed $\alpha_n=\alpha$, this allows
$\log d_n=o(n^{1/3})$ under arbitrary coordinate dependence.
This tail calibration uses no joint Gaussian approximation and hence
does not conflict with the dimension barrier of
\citet[Theorem~3.1]{KOCK2026BARRIER}.

The necessary conditions in our main theorems concern the joint
rectangle approximation of the original $T_n$ by its matching
$G_P$. Robust estimators and the conservative marginal calibration
above provide alternative routes to mean inference. Joint Gaussian
calibration can use the dependence structure to select critical
values; when $R_P$ is unknown, its implementation additionally
requires a covariance estimator and a Gaussian comparison bound,
or an appropriate bootstrap theorem.

\section{Proofs of the main results}\label{sec:proofs-main}

Section~\ref{sec:unrestricted} proves
Theorem~\ref{thm:main-polynomial}. For fixed positive gaps, we first
prove the sample-size bounds of Theorem~\ref{thm:main-fixed} in
Section~\ref{sec:fixed-rates}, then extend convergence to the full
range $\log d=o(n^{1/3})$ in Section~\ref{sec:fixed-cubic}. This order
allows the latter proof to use the sample-size bounds on groups of
polynomially many low-threshold constraints. Section~\ref{sec:shrinking-rate}
treats vanishing gaps and completes the proof of
Theorem~\ref{thm:main-shrinking}; it also reuses the rare-observation
construction from Section~\ref{sec:fixed-rates}. The finer dimension
bounds near the cubic boundary are proved in
Section~\ref{sec:shrinking-proofs}. Section~\ref{sec:qcs-proof} gives the
counterexample. The comparison and smoothing lemmas used throughout
are collected in Section~\ref{sec:inputs}.

\subsection{Arbitrary dependence: the error bound and dimension threshold}
\label{sec:unrestricted}
Fix $q\geq3$ and $K>1$. The relevant error scale is
\[
 D=\left(\frac d{n^{(q-2)/2}}\right)^{1/q}.
\]
Constants may depend on $q,K$; any dependence on $M$ is stated
explicitly. We measure distance from a rectangle boundary using the
score $g(x)=-\log\overline\Phi(x)$ from Lemma~\ref{lem:score-strip}.

\begin{theorem}[A quantitative lower bound for small correlation gaps]
\label{thm:polynomial-envelope}
The bounds and convergence criterion in
Theorem~\ref{thm:main-polynomial} hold. In addition, for each finite
$M$ there are $c_M>0,n_M$, depending only on $q,K,M$, such that,
whenever $n\geq n_M$, $2\leq d\leq n^{(q-2)/2}$, and
$0\leq a\leq M/\log n$,
\begin{equation}\label{eq:small-gap-lower-new}
 \Delta_{n,q,K}(d,a)\geq c_M\left\{n^{-(q-2)/(2q)}
 +\frac{x}{(x^{1/q}+\sqrt{2a\log(2d)})^{q-1}}\right\},
 \qquad x=\frac d{n^{(q-2)/2}}.
\end{equation}
\end{theorem}

\subsubsection{Gaussian comparison under a coupling}

\begin{lemma}[Gaussian comparison under a coordinatewise coupling]
\label{lem:coupled-Gaussian-new}
Let $(G,H)$ be jointly centered Gaussian $d$-vectors with $\mathbb E G_j^2=1$ and $\max_j\|H_j-G_j\|_2\leq r$. Then
\[
 d_{\mathrm{rect}}(G,H)\leq C\min\{1,r\log(ed)\}.
\]
No condition is imposed on covariance rank.
\end{lemma}
\begin{proof}
It suffices to consider $r<1/4$. A rectangle event can be written as a
lower-orthant event in $m\leq2d$ signed coordinates. Use $G,H$ also for
the corresponding signed Gaussian vectors, and let $Z\sim N_m(0,I_m)$
be independent of both. Set $Y_t=(1-t)G+tH+\varepsilon Z$ and $B=H-G$.
For $0<\varepsilon\leq1/4$, the marginal standard deviations of $Y_t$
lie in $[1/2,2]$, and its covariance matrix is positive definite. Write
$f_j$ for the density of $Y_{t,j}$ and
\[
 u_j=\mathbb P(Y_{t,-j}\leq b_{-j}\mid Y_{t,j}=b_j).
\]
Apply Gaussian anti-concentration for rectangles \citep{CCK2017} and
differentiate the rectangle probability with respect to a common
threshold shift. This gives
\begin{equation}\label{eq:face-sum-new}
 \sum_{j=1}^m f_j(b_j)u_j\leq C\sqrt{\log(em)}.
\end{equation}
The density transport identity
\[
 \partial_t f_t=-\nabla\!\cdot\{f_t\mathbb E(B\mid Y_t=\cdot)\}
\]
follows by differentiating smooth compactly supported test functions. Integration over the orthant therefore bounds the absolute derivative of its probability by
\[
 \sum_j f_j(b_j)\,
 \mathbb E\bigl[|B_j|\mathbf1_{\{Y_{t,-j}\leq b_{-j}\}}\mid Y_{t,j}=b_j\bigr].
\]
Conditionally on $Y_{t,j}=b_j$, $B_j$ is normal with mean of magnitude at most $2r|b_j|$ and standard deviation at most $r$. Gaussian tail integration bounds the displayed conditional expectation by
\[
 Cr u_j\{1+|b_j|+\sqrt{\log(1/u_j)}\},
\]
with the product interpreted as zero when $u_j=0$. The sum $\sum_j f_j(b_j)u_j|b_j|$ is at most $C\log(em)$. Indeed, use
\eqref{eq:face-sum-new} when $|b_j|\leq C_0\sqrt{\log(em)}$ and the
density bound $f_j(b)\leq C e^{-b^2/8}$ otherwise, with $C_0$ sufficiently
large. The sum with $|b_j|$ replaced by
$\sqrt{\log(1/u_j)}$ has the same bound. Split it at $u_j=e^{-2\log(em)}$: above this value use
\eqref{eq:face-sum-new}, and below it use
$u\sqrt{\log(1/u)}\leq C\sqrt u$.

Integrating in $t$ proves the result for the regularized vectors. Let
$\varepsilon\downarrow0$. The marginal variances of both endpoints are
positive when $r<1/4$, so the boundary of every fixed rectangle has
probability zero. This also covers rectangles of zero width.
\end{proof}

\subsubsection{Gaussian approximation for the truncated statistic}

We use the following Gaussian approximation bound for sums of bounded
vectors. If the vectors $\xi_i$ are independent and centered,
$\operatorname{Var}(\sum_i\xi_{ij})=1$, and
$\max_i\|\xi_i\|_\infty\leq b$, then
\begin{equation}\label{eq:CCKK-input-new}
 d_{\mathrm{rect}}\!\left(\sum_i\xi_i,G\right)
 \leq C\min\{1,b\log^2(ed)\}.
\end{equation}
Here $G$ has the covariance of $\sum_i\xi_i$. The bound follows from
\citet[Theorem~1, equation~(4)]{CCKK2026} and
\citet[Proposition~G.1]{CCKKSUP2026}; the formulation used here is stated
in \eqref{eq:external-cckk}. It allows singular covariance matrices.

We may assume $d\leq n^{(q-2)/2}$, since the discrepancy is always at
most one. Truncate each observation at the level
\[
 T=\frac{(nd)^{1/q}}{\log^2(ed)},\quad h=\frac T{\sqrt n}=\frac D{\log^2(ed)},
\]
\[
 Z_{ij}=X_{ij}\mathbf1_{\{|X_{ij}|\leq T\}},\qquad
 U_j=\frac{\sum_i Z_{ij}}{(\sum_i Z_{ij}^2)^{1/2}},\qquad
 S=\frac1{\sqrt n}\sum_i(Z_i-\mathbb EZ_i).
\]
For fixed $q$, $T\geq c_q n^{1/q}$. Hence the truncated marginal variances lie in $[1/2,2]$ for large $n$. Applying \eqref{eq:CCKK-input-new} after coordinatewise standardization gives
\[
 d_{\mathrm{rect}}(S,G_{\operatorname{Cov}(Z_1)})\leq Ch\log^2(ed)\leq CD.
\]
Construct this Gaussian and $G_R$ jointly with the covariance of
$(Z_1-\mathbb EZ_1,X_1)$. The standard deviation of each coordinate
difference is at most
\[
 Cr_0,\qquad r_0=T^{1-q/2}=\frac{D\log^{q-2}(ed)}{\sqrt d}.
\]
Lemma~\ref{lem:coupled-Gaussian-new} bounds the resulting Gaussian
comparison error by $Cr_0\log(ed)\leq C_qD$. The mean shift due to truncation satisfies
\[
 b_0:=\sqrt n\max_j|\mathbb EZ_{1j}|
 \leq C\sqrt nT^{1-q}=\frac{CD\log^{2q-2}(ed)}d,
 \qquad b_0\sqrt{\log(ed)}\leq C_qD.
\]
Bernstein's inequality and a union bound, using $\log(nd)\leq C_q\log(en)$, show that outside an event of probability $O(n^{-2})$,
\begin{align*}
 \|S\|_\infty&\leq C(\sqrt{\log(en)}+h\log(en)),\\
 \max_j\left|n^{-1}\sum_i Z_{ij}^2-1\right|
 &\leq e_n:=C\{r_0^2+v_n\sqrt{\log(en)}+h^2\log(en)\},
\end{align*}
where $v_n=n^{-1/2}\{1+T^{(4-q)_+/2}\}$. These quantities satisfy
\[
 e_n=o(1),\qquad h\sqrt{\log(en)}=O(1),\qquad
 e_n\sqrt{\log(en)\log(ed)}\leq C_qD.
\]
For completeness,
\[
 \frac{v_n}{D}\leq C_q n^{-\min\{1/q,(q-2)/(2q)\}}.
\]
If $d<n^{(q-2)/4}$ then $D\leq n^{-(q-2)/(4q)}$ and $h,r_0\leq C_qD$, which controls the remaining terms after multiplication by $\sqrt{\log(en) \log(ed)}$. If $d\geq n^{(q-2)/4}$, then $\log(ed)\asymp_q\log(en)$, $h=D/\log^2(ed)$, and $r_0^2=D^2\log^{2q-4}(ed)/d$, giving the same conclusions. On the preceding event from
Bernstein's inequality,
\[
 \|U-S\|_\infty\leq C(b_0+e_n\sqrt{\log(en)}).
\]
Gaussian rectangle anti-concentration, together with $n^{-2}\leq D$, proves
\begin{equation}\label{eq:clipped-core-new}
 d_{\mathrm{rect}}(U,G_R)\leq CD.
\end{equation}
Write a fixed rectangle in signed coordinates as $y_j\leq b_j$,
and set $M(y)=\max_j\{g(y_j)-g(b_j)\}$. Lemma~\ref{lem:score-strip}
and the preceding comparison give
\begin{equation}\label{eq:clipped-score-strip-new}
 \mathbb P\{|M(U)|\leq\delta\}\leq2\delta+CD.
\end{equation}
This application of \eqref{eq:clipped-core-new} is valid because a
sublevel set of the maximum score is a rectangle or the empty set.

\subsubsection{The error from restoring truncated observations}

Fix a coordinate $j$. Condition on
$I_j=\{i:|X_{ij}|>T\}$ and on the values at those indices. The remaining
$N=n-|I_j|$ observations, denoted by $Y_1,\ldots,Y_N$, are independent
and have the conditional law of $X_{1j}$ given $|X_{1j}|\leq T$.
Their common mean satisfies $|\mu_j|\leq CT^{1-q}$, their variance is
$1+o(1)$ uniformly in $j$, and their third absolute moment is bounded
by a constant depending only on $q,K$. Since
\[
 n\mathbb P(|X_{1j}|>T)\leq KnT^{-q}=\frac{K\log^{2q}(ed)}d\leq C_{q,K},
\]
the probability that any column has $N<n/2$ is at most $d2^n(C_{q,K}/n)^{n/2}=O(n^{-2})$.

For $N\geq n/2$, the iid self-normalized Berry--Esseen bound of \citet{BG1996}, followed by the centering correction below, gives
\begin{equation}\label{eq:conditional-scalar-new}
 \sup_t\left|\mathbb P\left\{
 \frac{\sum_{i=1}^N Y_i}{(\sum_{i=1}^N Y_i^2)^{1/2}}\leq t
 \,\middle|\, I_j,(X_{ij})_{i\in I_j}\right\}-\Phi(t)\right|
 \leq C(n^{-1/2}+\sqrt nT^{1-q}).
\end{equation}
To control the centering correction, first note that both
$\sum_iY_i^2$ and $\sum_i(Y_i-\mu_j)^2$ exceed $n/7$ except on an event
of probability at most $Ce^{-cn}$. To see this, truncate each square at
a sufficiently large fixed constant. The third absolute moment bound
ensures a uniform positive lower bound on its truncated mean, after
which Hoeffding's inequality applies. Subtracting $\mu_j$ changes the numerator by at most $N|\mu_j|$ and the norm by at most $\sqrt N|\mu_j|$. Since a centered self-normalized ratio is bounded in absolute value by $\sqrt N$, the ratio changes by at most $C\sqrt n|\mu_j|$ on this event. A bounded Gaussian density then yields \eqref{eq:conditional-scalar-new}. Student's statistic is a monotone transform of the centered
self-normalized ratio. That transform changes the Gaussian comparison
by $O(N^{-1})$, so either formulation of the scalar theorem yields the
stated bound.

Set
\[
 W_j=n^{-1/2}\sum_{i\in I_j}|X_{ij}|,
 \qquad Q=\sum_{i=1}^NY_i^2.
\]
On $Q\geq n/7$, reinsertion transforms the retained ratio $t$ into
\[
 f(t)=\frac{t+v}{\sqrt{1+r}},\qquad
 |v|\leq w:=\sqrt7W_j,\qquad 0\leq r\leq w^2.
\]
The same argument applies to a lower constraint after reversing the
signs of the observations and threshold. On $Q\geq n/7$, these deterministic bounds on $v,r$ hold for every
retained sample. The argument below does not condition on $Q$.

We next bound the probability that restoring the omitted observations
changes whether a threshold is exceeded. Fix $0<\delta\leq1/4$ and a
threshold $b$. Consider the values of $t$ for which some
$|v|\leq w$, $0\leq r\leq w^2$ make $t$ and $f(t)$ give different outcomes for the constraint
$t\leq b$ and satisfy $|g(f(t))-g(t)|>\delta$. These values lie in at most
two intervals with total standard normal probability at most
\begin{equation}\label{eq:crossing-mass-new}
 C_q\delta^{1-q}w^q.
\end{equation}
Indeed, for $w\leq1/8$, crossing implies
\[
 |t-b|\leq w+w^2|b|.
\]
This interval has normal mass at most $Cw(1+w|b|)e^{-cb^2}\leq Cw$. Since $|f(t)-t|\leq w+w^2|t|$ and $g'(t)\leq C(1+|t|)$, a score displacement exceeding $\delta$ requires $w(1+|t|)>c\delta$. Intersecting the interval with this set gives at most two intervals. If $w<c'\delta$, the intersection is empty unless $|b|\geq c''\delta/w$, and its mass is at most
\[
 Cw\exp(-c'''\delta^2/w^2)\leq C_q\delta^{1-q}w^q.
\]
For $c'\delta\leq w\leq1/8$, use $Cw\leq C_q\delta^{1-q}w^q$; for $w>1/8$, the unit bound suffices. This proves \eqref{eq:crossing-mass-new} using sets determined only by $b,w,\delta$.

Outside $\{|M(U)|\leq\delta\}$, a change in rectangle membership
requires a coordinate to cross a threshold with score displacement
greater than $\delta$. If $U$ is inside the rectangle, select a
constraint violated by the full statistic. If $U$ is outside, select a
constraint attaining $M(U)$. Apply \eqref{eq:conditional-scalar-new} to
the intervals in \eqref{eq:crossing-mass-new}. Only coordinates with
$I_j\ne\varnothing$ contribute. Summing these conditional errors and
the exceptional probabilities for the denominators gives at most
\[
 CndT^{-q}(n^{-1/2}+\sqrt nT^{1-q}+e^{-cn})\leq CD.
\]
To see the last inequality, its three terms are respectively bounded by $C\log^{2q}(ed)n^{-1/2}$, $CD\log^{4q-2}(ed)/d$, and $C\log^{2q}(ed)e^{-cn}$. Divide the first by $D$ and use $\log^{2q}(ed)/(nd)^{1/q}=O_q(1)$; every fixed power of $\log(ed)$ divided by $d$ is bounded.

Apply Rosenthal's inequality \citep{Rosenthal1970} to
$A_i=|X_{ij}|\mathbf1_{\{|X_{ij}|>T\}}$ after centering, and add the
mean term. This gives
\[
 n^{q/2}\mathbb E W_j^q
 \leq C_q\{n\mathbb EA_1^q+(n\mathbb EA_1^2)^{q/2}+(n\mathbb EA_1)^q\}
 \leq C_{q,K}n.
\]
The last inequality follows from $\mathbb EA_1^r\leq KT^{r-q}$ for
$r=1,2,q$ and $nT^{-q}\leq C_q$. Thus
$\sum_j\mathbb EW_j^q\leq CD^q$. Combining \eqref{eq:clipped-core-new}, \eqref{eq:clipped-score-strip-new}, and \eqref{eq:crossing-mass-new} yields
\[
 \Delta_{n,q,K}(d,0)\leq C\{D+\delta+D^q\delta^{1-q}\}.
\]
The terms $\delta$ and $D^q\delta^{1-q}$ bound, respectively, the
boundary strip and the crossings caused by large observations. They
are of the same order at $\delta=D$. Use this choice when $D\leq1/4$
and the bound by one otherwise. This proves the upper bound in
Theorem~\ref{thm:polynomial-envelope}.

\subsubsection{Lower bounds with the prescribed Gaussian covariance}

\begin{lemma}[A nearly antipodal pair]
\label{lem:antipodal-new}
For all sufficiently large $n$, every $d\geq2$, and every $a\in[0,1]$,
\[
 \Delta_{n,q,K}(d,a)\geq c_{q,K}n^{-(q-2)/(2q)}.
\]
\end{lemma}
\begin{proof}
Let $U$ be a symmetric sign and let $V$, independent of $U$, equal each of $\pm n^{1/q}$ with probability $\lambda/(2n)$ and zero otherwise. Choose a fixed $\lambda>0$ sufficiently small in terms of $q,K$, and standardize $(U+V,-U+V)$ to unit variances. Minkowski's inequality verifies the moment bound, since $\|U\pm V\|_q\leq1+\lambda^{1/q}$. The correlation is negative for large $n$, so every $a\leq1$ is allowed.

For independent rows put $h=n^{-(q-2)/(2q)}$, $J=n^{-1/q}\sum_iV_i$, and $W=n^{-1/2}\sum_iU_i$. For the lower orthant with both thresholds equal to zero, the
denominators are positive and the event reduces to $hJ\leq-|W|$. The variables $J,W$ are independent, $J$ is symmetric, and $\mathbb E|J|\leq\lambda$. The scalar Berry--Esseen bound, applied with closed endpoints, gives
\[
 \mathbb P(hJ\leq-|W|)
 \leq\mathbb E\{\Phi(h|J|)-1/2\}+Cn^{-1/2}
 \leq\frac{\lambda h}{\sqrt{2\pi}}+Cn^{-1/2}.
\]
The matched Gaussian can be represented, up to a common positive scale, as $(Z+h\sqrt\lambda Z',-Z+h\sqrt\lambda Z')$. By rotational symmetry its zero-orthant probability equals $\pi^{-1}\arctan(h\sqrt\lambda)$, which is at least $h\sqrt\lambda/(2\pi)$ for large $n$. Choose $\lambda$ so that $\sqrt\lambda/(2\pi)>2\lambda/\sqrt{2\pi}$. Since $n^{-1/2}=o(h)$, the discrepancy is at least $c h$. Append independent symmetric sign coordinates and leave them unconstrained in the test rectangle to obtain every $d\geq2$.
\end{proof}

We next prove \eqref{eq:small-gap-lower-new}. All thresholds in this
construction are zero, so the test depends only on the signs of the
numerators. It therefore also gives a lower bound for the ordinary
normalized sum; see \citet{KP2024} for related dimension lower bounds.
Assume $2\leq d\leq n^{(q-2)/2}$ and $a\log n\leq M$. Put
\[
 x=d/n^{(q-2)/2},\quad z=x^{1/q},\quad
 t=\sqrt{2a\log(2d)},\quad \rho=1-a,\quad b=8(z+t),
\]
\[
 A=b\sqrt n,\qquad \mu=\eta x/b^q,\qquad
 \pi_0=\mu/(nd),\qquad v=A^2\pi_0(1-\pi_0).
\]
Choose $\eta>0$ sufficiently small in terms of $q,K$. Then $\mu\leq1/8$ and $A^q\pi_0=\eta$. Let $C_i$ be independent symmetric signs, each shared by all
coordinates in row $i$, and let $V_{ij}$ be independent signs of mean
$\sqrt\rho$, independent of the $C_i$. Put $Y_{ij}=C_iV_{ij}$. Thus each $Y_{ij}$ is a symmetric sign, and distinct coordinates in
the same row have correlation $\rho$. With independent Bernoulli$(\pi_0)$ marks, define
\[
 X_{ij}=\frac{Y_{ij}+A(I_{ij}-\pi_0)}{\sqrt{1+v}}.
\]
The variables are centered and have unit variance; their off-diagonal correlations are $\rho/(1+v)\leq1-a$. Since $A\|I-\pi_0\|_q\leq2\eta^{1/q}$, Minkowski's inequality verifies
the moment bound. All denominators are positive for large $n$.

Write $S_j=n^{-1/2}\sum_iY_{ij}$, $N_j=\sum_iI_{ij}$, and $e_0=b\mu/d$. Because the denominators are positive, the zero-threshold sample
event is $\max_j(S_j+bN_j)\leq e_0$. Nonnegative marks can only decrease its probability relative to $\max_jS_j\leq e_0$. The probability of exactly one mark is at least $\mu/2$. Conditional on its column $j$, the crossing probability is
\[
 \mathbb P(S_k\leq e_0\ \forall k)
 -\mathbb P(S_k\leq e_0\ (k\ne j),\ S_j\leq e_0-b).
\]
By \eqref{eq:CCKK-input-new}, each probability differs from its Gaussian counterpart by at most
\[
 e_n=C_qn^{-1/2}\log^2(en).
\]
The Gaussian vector with the covariance of $Y_i$ is
$G_j^0=\sqrt\rho Z_0+\sqrt a Z_j$. The independent events
\[
 -t-b/2\leq\sqrt\rho Z_0\leq-t-b/4,
 \qquad 0\leq Z_j\leq1,
 \qquad \max_{k\ne j}Z_k\leq\sqrt{2\log(2d)}
\]
have joint probability at least $c_{q,M}b$. Here $t,b$ remain in a fixed bounded interval, $\rho\geq1/2$ eventually, and the final probability is at least $1/2$. On this event, every $G_k^0$ is negative, whereas adding $b$ to
coordinate $j$ raises it above $e_0$. Here we use $b\geq8t$ and
$e_0\leq b/8$. It follows that the marks decrease the sample rectangle
probability by at least $c_{q,M}\mu b-Ce_n$.

The Gaussian with the covariance of $X_i$ has coordinates $(\sqrt\rho Z_0+\sqrt{a+v}Z_j)/\sqrt{1+v}$. Conditioning on $\max_jZ_j$ and using the bounded density of $\sqrt\rho Z_0$ shows that its zero-orthant probability differs from $\mathbb P(\max_jG_j^0\leq e_0)$ by at most
\[
 C\{e_0+(\sqrt{a+v}-\sqrt a)\sqrt{\log(ed)}\}.
\]
As $d\to\infty$, this is $o(\mu b)$ uniformly over the stated range
of $n,a$: when $t\leq z$, $\mu$ is bounded below and the ratio of the second term to $\mu b$ is at most $C_q\sqrt{\log(ed)/d}$; when $t>z$, use $\sqrt{a+v}-\sqrt a\leq v/\sqrt a$, $b\leq16t$, and $v\leq b^2\mu/d$ to bound the ratio by $C_q\log(ed)/d$. Also $e_0/(\mu b)=1/d$. Consequently, for all $d\geq d_0(q,K,M)$,
\[
 \Delta_{n,q,K}(d,a)\geq c\mu b-Ce_n
 =c'\frac{x}{(z+t)^{q-1}}-Ce_n.
\]
Lemma~\ref{lem:antipodal-new} and $e_n=o(n^{-(q-2)/(2q)})$ absorb the last error and give the sum in \eqref{eq:small-gap-lower-new}. For $d<d_0$, its second term is at most $x^{1/q}=d^{1/q}n^{-(q-2)/(2q)}$, so the antipodal lemma alone suffices.

When $a=0$, \eqref{eq:small-gap-lower-new} yields $cD$ for $d\leq n^{(q-2)/2}$. For larger $d$, use $\lfloor n^{(q-2)/2}\rfloor$ coordinates and
append independent symmetric signs. Finally, if $a_n\log n$ is bounded and $d_n/n^{(q-2)/2}$ does not tend to zero, pass to a subsequence and use a fixed positive multiple of $n^{(q-2)/2}$ coordinates in \eqref{eq:small-gap-lower-new}. Then its second term stays bounded below. This proves the dimension criterion and completes the proof of
Theorem~\ref{thm:polynomial-envelope}.

\subsection{Fixed correlation gaps: rates and lower bounds}
\label{sec:fixed-rates}
Fix $q>2$, $K>1$, and $a\in(0,1]$. The sample-size bounds here also
cover $2<q<3$; the dimension criterion in
Section~\ref{sec:fixed-cubic} assumes $q\geq3$. Define
\[
 B_q(n,d)=n^{-(q-2)/(2q)}
 \{\log^{3/2}(end)+\log^{3/2}(ed)\sqrt{\log(en)}\},
 \qquad a_c=\frac1{q+1}.
\]
Write $r=a/(1-a)$, with $r=\infty$ when $a=1$.
Constants in this subsection may depend on $q,K,a$.
 
\begin{theorem}[Optimal power of the sample size and the endpoint correction]\label{thm:fixed-sample-rate}
For every $n,d\geq2$,
\begin{equation}\label{eq:fixed-sample-upper}
 \Delta_{n,q,K}(d,a)\leq1\wedge C\{B_q(n,d)+E_{q,a}(n)\},
\end{equation}
where, for some finite $b=b(q,K,a)$,
\[
 E_{q,a}(n)=
 \begin{cases}
 n^{-(q-2)r/2}\exp\!\left\{\dfrac{\sqrt{a(q-2)}}{1-a}\sqrt{\log n}\right\}\log^b(en),
       &0<a\leq a_c,\\[2mm]
 0,&a_c<a\leq1.
 \end{cases}
\]
There is $c_{q,K}>0$ such that, for all sufficiently large $n$ and all $d\geq2$,
\begin{equation}\label{eq:fixed-universal-lower}
 \Delta_{n,q,K}(d,a)\geq c_{q,K}\left[1\wedge
 \{n^{-(q-2)/(2q)}\sqrt{\log(ed)}+n^{-1/2}\log^{3/2}(ed)\}\right].
\end{equation}
For $0<a\leq a_c$, there are integers
$m_n\asymp n^{(q-2)/2}(\log n)^{q+1}$ and a finite $b'>0$ such that, for all large $n$ and $d\geq m_n$,
\begin{equation}\label{eq:fixed-rare-lower}
 \Delta_{n,q,K}(d,a)\geq
 n^{-(q-2)r/2}\exp\!\left\{\frac{\sqrt{a(q-2)}}{1-a}\sqrt{\log n}\right\}\log^{-b'}(en).
\end{equation}
These bounds imply all assertions of Theorem~\ref{thm:main-fixed}.
For $0<a\leq a_c$, at the same dimensions $m_n$ they also give
\begin{equation}\label{eq:fixed-endpoint-log}
 \log\Delta_{n,q,K}(m_n,a)
 =-\frac{q-2}{2} r\log n+\frac{\sqrt{a(q-2)}}{1-a}\sqrt{\log n}
   +O(\log\log n).
\end{equation}
\end{theorem}

\subsubsection{The truncated sum and coordinates near its maximum}
For the upper bound, it suffices to treat $0<a\leq1/2$. If $a>1/2$,
enlarge the class by replacing $a$ with $1/2$. Since
$1/2>1/(q+1)$ for every $q>2$, this gives the required upper bound.
The lower bounds below are constructed at the stated value of $a$.

Set
\begin{align*}
 \varepsilon&=n^{-(q-2)/(2q)},&
 \tau&=2(Kn)^{1/q},\\
 Y_{ij}&=X_{ij}\mathbf1_{\{|X_{ij}|\leq\tau\}},&
 \mu_j&=\mathbb E Y_{ij},\\
 U_{ij}&=(Y_{ij}-\mu_j)/\sqrt n,&
 S_j&=\sum_{i=1}^nU_{ij}.
\end{align*}
Let $\delta=2\tau/\sqrt n\asymp\varepsilon$, and let $k_j$ count entries in column $j$ exceeding $\tau$ in absolute value. The marginal moment condition gives
\begin{equation}\label{eq:clip-moments-fixed}
 n\mathbb P(|X_j|>\tau)\leq\lambda_0:=2^{-q},\quad
 \sqrt n|\mu_j|\leq C\delta,\quad
 \mathbb E(X_j-Y_j)^2\leq C\delta^2,\quad |U_{ij}|\leq\delta.
\end{equation}
Put $B=B_q(n,d)$. We may assume that $B$ is less than a sufficiently small fixed constant and $n$ is sufficiently large; otherwise increase the constant in the asserted bound. In this range,
\begin{equation}\label{eq:fixed-small-regime}
 \log(end)\leq C\varepsilon^{-2/3}=o(n),\quad
 \delta\sqrt{\log(end)}=o(1),\quad \delta \log^{3/2}(end)\leq C.
\end{equation}
Write the constraints of a fixed rectangle as $X_j\leq b_j$ and
$-X_j\leq-c_j$, omitting infinite endpoints. We call these
signed-coordinate constraints faces and group them according to the
sign used. The following estimates apply to the resulting signed
vectors, separately within each sign group.

For a deterministic set $I$ of $k\leq C\log(end)$ row indices, let
$W=\sum_{i\notin I}U_i$ and $V=\operatorname{Cov}(W)$. Uniformly over
these choices of $I$, the marginal variances of $W$ lie in $[0.85,2]$
for large $n$, and
\begin{equation}\label{eq:fixed-cov-margin}
 V_{jj}-V_{jl}\geq a/2
\end{equation}
for distinct coordinates in the same sign group. Indeed, $V=(1-k/n)\operatorname{Cov}(Y)$, and each entry of $\operatorname{Cov}(Y)-R$ is $O(\delta)$ by Cauchy--Schwarz and \eqref{eq:clip-moments-fixed}.

The bounded-sum theorem of \citet[Theorem~1, equation~(4)]{CCKK2026}, after coordinate standardization, yields
\begin{equation}\label{eq:fixed-bounded-clt}
 e_W:=d_{\rm rect}(W,N(0,V))\leq C\delta \log^{3/2}(ed)\sqrt{\log(en)}.
\end{equation}
The same estimate holds after restricting to any nonempty set $J$ of
coordinates. After standardization, each row is bounded by $O(\delta)$,
and the expected Gaussian maximum is $O(\sqrt{\log(e|J|)})$. The
additional factor in the cited theorem
$[1+\log_+(1/(C\delta\sqrt{\log(e|J|)}))]^{1/2}$ is $O(\sqrt{\log(en)})$. To compare $N(0,\operatorname{Cov}(Y))$ with $G_R$, construct them
jointly with the covariance of $(Y-\mu,X)$. Each coordinate difference
then has standard deviation at most $C\delta$. A union bound at
$C\delta\sqrt{\log(end)}$, followed by Gaussian anti-concentration for
rectangles, gives
\begin{equation}\label{eq:fixed-S-compare}
 d_{\rm rect}(S,G_R)\leq CB.
\end{equation}
The anti-concentration step uses only positive marginal variances.

Bernstein's inequality, $\mathbb EY_j^4\leq C(1+\tau^{(4-q)_+})$, and a binomial tail bound give, except on an event of probability $O(n^{-3})$,
\begin{equation}\label{eq:fixed-global-good}
 \begin{aligned}
 \max_j k_j&\leq C\log(end),\qquad
 \|S\|_\infty\leq C\sqrt{\log(end)},\\
 \max_j\left|n^{-1}\sum_iY_{ij}^2-1\right|
 &\leq C\{\sqrt{\log(end)/n}+\delta^2\log(end)\}=:e_Q.
 \end{aligned}
\end{equation}
For $2<q<4$, the additional variance term
$\tau^{(4-q)/2}\sqrt{\log(end)/n}$ is $O(\delta^2\log(end))$; for $q\geq4$ it is $O(\sqrt{\log(end)/n})$.
Cauchy--Schwarz gives
\[
 |T_j|^2\leq
 \frac{(S_j+\sqrt n\mu_j)^2}{n^{-1}\sum_iY_{ij}^2}+k_j.
\]
Consequently, intersecting a test rectangle with $[-M,M]^d$, where
$M=C\sqrt{\log(end)}$, changes its probability under either $T$ or $S$
by at most $O(n^{-3})$. Choose $D=C\delta\log(end)$ large enough to
bound $Me_Q$, $C\log(end)\delta$, and $\sqrt n\|\mu\|_\infty$. Also $MD+D^2=O(1)$ by \eqref{eq:fixed-small-regime}.

\begin{lemma}[Counting coordinates near the maximum]
\label{lem:retained-face-count}
Let $\xi_1,\ldots,\xi_N$ be independent centered vectors in $\mathbb R^p$, with $|\xi_{ij}|\leq\delta$. Put $W=\sum_i\xi_i$ and $V=\operatorname{Cov}(W)$, and assume
\[
 \frac34\leq V_{jj}\leq2,\qquad V_{jj}-V_{jk}\geq a/2
 \quad(j\ne k),\qquad 0<a\leq1.
\]
For real thresholds $t_j$, define
\[
 M_W=\max_{j\leq p}(W_j-t_j),\qquad
 d_j=M_W-(W_j-t_j),\qquad e_W=d_{\rm rect}(W,N_p(0,V)).
\]
Fix $Q_0<\infty$, $\Gamma\geq\log(ep)$ and $M\geq0$, and let
\[
 \mathcal B=\left\{\|W\|_\infty\leq M,\quad
                    \max_j\sum_i\xi_{ij}^2\leq Q_0\right\}.
\]
There are $c>0$ and $A,C<\infty$, depending only on $Q_0$, such that, if $w\geq0$ and
\begin{equation}\label{eq:fixed-face-parameters}
 \delta M\leq ca,\qquad
 \lambda\geq A\left(\frac wa+\frac{\delta \Gamma}{a^2}\right),
 \qquad 0<\lambda\delta\leq1,
\end{equation}
then
\begin{align}
 \sum_{j\leq p}\mathbb P(\mathcal B,d_j\leq w)
 &\leq C\exp\{C\lambda(M+\lambda)\},
 \label{eq:fixed-face-count}\\
 \sum_{j:t_j\leq v}
 \mathbb P(\mathcal B,d_j\leq w,\ |W_j-t_j|\leq w)
 &\leq C\exp\{C\lambda(v_++\lambda)\}
             \min\{1,\lambda\sqrt{\Gamma}+e_W\},
 \qquad v\in\mathbb R.
 \label{eq:fixed-face-shell}
\end{align}
The covariance matrix may be singular. The constants $c,A,C$ do not depend on $a$.
\end{lemma}

\begin{proof}
Fix $j$. Independently for each row, with probability $\rho=\lambda\delta$ replace that row by a fresh draw whose density relative to its original law is $1+\xi_{ij}/\delta$; otherwise leave it unchanged. This is a probability density because the row is centered and $|\xi_{ij}|\leq\delta$. Denote the new sum by $W^{(j)}$.

Conditionally on the original array,
\[
 \mathbb E(W_l^{(j)}-W_l\mid \xi_1,\ldots,\xi_N)
 =\lambda V_{jl}-\lambda\delta W_l .
\]
A fresh tilted row satisfies $\mathbb E[(\xi_{il}^{(j)})^2]\leq2\mathbb E\xi_{il}^2$. Thus on $\mathcal B$ the sum of the conditional variances of the changes in coordinate $l$
is at most
\[
 2\lambda\delta\left(2V_{ll}+\sum_i\xi_{il}^2\right)
 \leq C(Q_0)\lambda\delta .
\]
Every centered change is bounded by $4\delta$. Bernstein's inequality gives, for a sufficiently small fixed $b_0>0$,
\[
 \mathbb P\left(
 |W_l^{(j)}-W_l-\lambda V_{jl}+\lambda\delta W_l|
 >b_0a\lambda\ \middle|\ \xi_1,\ldots,\xi_N\right)
 \leq 2\exp\{-c(Q_0)a^2\lambda/\delta\}.
\]
Choose $A$ sufficiently large in the lower bound on $\lambda$. The
union of these exceptional events over $p$ coordinates then has
conditional probability at most $1/2$. On $\mathcal B$, the absolute value of $\lambda\delta W_l$ is at
most $ca\lambda$. Decreasing $c$ gives, with
conditional probability at least $1/2$,
\begin{equation}\label{eq:fixed-face-success}
 \max_l|W_l^{(j)}-W_l-\lambda V_{jl}|\leq a\lambda/16.
\end{equation}
The constants also give $w\leq a\lambda/16$.

On $\{d_j\leq w\}$ and this success event, for $k\ne j$,
\[
 (W_j^{(j)}-t_j)-(W_k^{(j)}-t_k)
 \geq-w+\lambda(V_{jj}-V_{jk})-a\lambda/8>0.
\]
Thus $j$ uniquely maximizes $W_l^{(j)}-t_l$ after replacement, and
$W_j^{(j)}\leq M+3\lambda$ on $\mathcal B$. If also $|W_j-t_j|\leq w$ and $t_j\leq v$, then
\[
 \lambda/2\leq\max_l(W_l^{(j)}-t_l)\leq3\lambda,
 \qquad W_j^{(j)}\leq v+3\lambda.
\]
Here the lower bound uses $V_{jj}\geq3/4$ and
$w+a\lambda/16\leq\lambda/8$.

The unconditional new row law has density $1+\lambda\xi_{ij}$ relative to its original law. Hence the new product law $Q_j$ satisfies
\[
 \frac{dQ_j}{dP}(\xi_1,\ldots,\xi_N)
 =\prod_i(1+\lambda\xi_{ij})\leq e^{\lambda W_j}.
\]
This includes $\lambda\delta=1$, when a factor may vanish. On $\mathcal B$, the conditional probability of
\eqref{eq:fixed-face-success} is at least $1/2$. Thus each event before replacement has probability at most
twice that of the corresponding event under $Q_j$. The image events specify different unique maximizers, so they are
pairwise disjoint when evaluated under $P$. On each of them, $W_j\leq M+3\lambda$. Applying the density
bound and summing in $j$ proves \eqref{eq:fixed-face-count}.

For the second assertion, the same argument restricts all image
events to $\lambda/2\leq M_W\leq3\lambda$, with
$W_j\leq v+3\lambda$. Compare the two endpoints of this interval with their Gaussian
counterparts and apply Gaussian anti-concentration \citep{CCK2017}.
The resulting bound is the smaller of one and
$C\lambda\sqrt{\log(ep)}+2e_W$. Taking a limit from below includes a possible atom at the lower
endpoint; the Gaussian maximum is continuous because each marginal variance is positive. This proves \eqref{eq:fixed-face-shell}.
\end{proof}

We apply the lemma to each sign group of the retained vector, with
$M=C\sqrt{\log(end)}$, $Q_0=C$, and
$\lambda=A(w/a+\delta \log(end)/a^2)$. In the range \eqref{eq:fixed-small-regime},
$\delta\sqrt{\log(end)}\to0$, so $\delta M\leq ca$ eventually. For moderate entries below $\sqrt{n/\log(end)}$, the widths used below satisfy
\[
 w=C(\delta \log(end)+kv/\sqrt n),\qquad
 \lambda\leq C_a(\delta \log(end)+k/\sqrt{\log(end)}),\qquad k\leq C\log(end).
\]
It follows that $\lambda\delta=o(1)$ and
$\lambda(M+\lambda)\leq C_a(1+k)$. For a crossing into the rectangle caused by an entry above
$\sqrt{n/\log(end)}$, we use $w=2k\delta$; then
$\lambda(M+\lambda)\leq C_a$. For a crossing out of the rectangle,
with $v_0\leq t_j<2v_0$, we use \eqref{eq:fixed-face-shell} with
threshold bound $2v_0$. The inequalities
$v_0w+w^2\leq C(1+k)$ and $v_0\leq C\sqrt{\log(end)}$ imply
$\lambda(2v_0+\lambda)\leq C_a(1+k)$ and
$\lambda\delta=o(1)$. These bounds verify the hypotheses of
Lemma~\ref{lem:retained-face-count} in each case used below.

\subsubsection{Threshold crossings after conditioning on the remaining rows}
There are two ways in which $S$ and $T$ can disagree on a rectangle
$A$. If $S\notin A$ and $T\in A$, call this an inward crossing and
select a maximally violated face in a sign group violated by $S$. If
$S\in A$ and $T\notin A$, call this an outward crossing and select a
face violated by $T$. In the selected coordinate $j$, let $I$ be the
$k$ indices with $|X_{ij}|>\tau$, and delete these rows from every
coordinate of $S$ to obtain $W$.

We first derive necessary conditions on $W$ for the crossing. For each
fixed $j,I$, we then drop the requirement that no other entry in column
$j$ exceeds $\tau$. The remaining restriction on the large entries
depends only on the deleted rows and is independent of $W$. Its
probability can therefore be bounded before summing the estimates for
the candidate coordinates. This argument fixes $j,I$ before using
independence; it does not condition on a randomly selected maximizer.

Write $e_S$ for $e_W$ with $I=\varnothing$. If $k=0$, a crossing places the maximum of $S_j-t_j$ over the
relevant sign group in $[-D,D]$, with probability at most
$O(D\sqrt{\log(end)}+e_S)=O(B)$. Put $u=\sqrt{n/\log(end)}$. If the largest omitted entry lies in a dyadic interval $(v,2v]$ below $u$, $v\geq\tau$, then
\[
 |T_j-S_j|\leq C(D+kv/\sqrt n).
\]
Indeed the normalized numerator changes by at most $C\delta+2kv/\sqrt n$ and the normalized squared denominator by at most $4kv^2/n$. Use $n^{-1}\sum_iY_{ij}^2\geq1/2$ and $v\leq u$. Since $\|W-S\|_\infty\leq k\delta$, both crossing directions satisfy the event conditions in
\eqref{eq:fixed-face-shell}, with $w=C(D+kv/\sqrt n)$ and threshold
bound $M$. Now $\lambda(M+\lambda)\leq C(1+k)$ and $\lambda\delta=o(1)$. Summing over all choices of deleted positions gives probability at
most
\[
 nKv^{-q}\frac{\lambda_0^{k-1}}{(k-1)!}.
\]
The dyadic sums of $nKv^{-q}$ and $(v/\sqrt n)nKv^{-q}$ are $O(1)$ and $O(\delta)$. Multiply \eqref{eq:fixed-face-shell} by these probability bounds and
sum over $k$. The factorial denominator makes the sum finite despite the
exponential factors in $k$, and the total is $O(B)$.

It remains to treat deleted entries whose magnitude exceeds $u$. Set
\[
 z=nKu^{-q}=Kn^{-(q-2)/2}\log^{q/2}(end)=K(\varepsilon\sqrt{\log(end)})^q=o(1).
\]
For an inward crossing involving at least one entry above $u$, the
selected face of $W$ is within $2k\delta$ of the maximum in its sign
group. At this width, \eqref{eq:fixed-face-count} is $O(1)$ because
$k\leq C\log(end)$ and $\delta\log^{3/2}(end)=O(1)$. The probability
bound for the deleted entries is $z\lambda_0^{k-1}/(k-1)!$, so their
total contribution is $O(B)$.

For an outward crossing with $b\geq1$ entries above $u$ and $h=k-b$
other deleted entries, the probability bound, including all choices of
positions, is
\begin{equation}\label{eq:fixed-mark-weight}
 \frac{z^b\lambda_0^h}{b!h!}.
\end{equation}
If $t_j>\sqrt b$, deleting only the entries above $u$ from the ratio
leaves a value greater than $\sqrt{t_j^2-b}$. In fact the removed numerator is at most $\sqrt b$ times its Euclidean norm, and
$t\sqrt{x^2+y^2}-\sqrt b\,y\geq x\sqrt{t^2-b}$.
The remaining ratio differs from $S_j$ by at most $C(D+h/\sqrt{\log(end)})$. For a dyadic threshold bin $v\leq t_j<2v$, $v\geq\sqrt b$, apply the face-shell bound at
$w=C(D+k\delta+h/\sqrt{\log(end)}+b/v)$. Here $vw+w^2\leq C(1+k)$, because $h^2/\log(end)\leq Ch$ and $b^2/v^2\leq b$. Summing over the $O(\log(end))$ threshold intervals and using
\eqref{eq:fixed-mark-weight} gives $O(z\log(end))=O(B)$. The remaining
outward crossings involve a selected threshold $t_j\leq\sqrt b$.

\subsubsection{Rectangles with many low thresholds}
Select $m$ signed coordinates from one sign group: either $G_{R,j}$
for upper constraints or $-G_{R,j}$ for lower constraints. Let $F_m(t)$
be the probability that all selected coordinates are at most $t$.
Their pairwise correlations are at most $1-a$, so Slepian's inequality
\citep{Slepian1962} gives, for $0<a<1$,
\[
 F_m(t)\leq\Phi\!\left(\frac{t-\sqrt a\,v}{\sqrt{1-a}}\right)+\Phi(v)^m.
\]
Take $v=(2\log m-6\log\log m)^{1/2}$ for large $m$. Mills' inequalities and $(x-y)_+^2\geq x^2-2xy$ yield
\begin{equation}\label{eq:fixed-screen}
 F_m(t)\leq C m^{-r}\exp\{c_a t_+\sqrt{\log m}+C\log\log(eem)\},
 \qquad c_a=\frac{\sqrt{2a}}{1-a}.
\end{equation}
The term $\Phi(v)^m$ is smaller than $m^{-r-1}$ for large $m$;
increasing $C$ covers the remaining values of $m$. If $r>1/q$, the strict exponent inequality implies
\begin{equation}\label{eq:fixed-screen-fast}
 F_m(t)\leq Ce^{Ct_+^2}m^{-1/q}.
\end{equation}
When $a=1$, apply this bound to the larger class with $a=1/2$, for
which $r=1>1/q$.

Let $m_b$ count the constraints with $t_j\leq\sqrt b$ in one sign group. Since we have reduced to $a\leq1/2$, we have $0<r\leq1$. Put
\[
 y=\log(1/z),\qquad
 Q=z^r\exp\{c_a\sqrt y+C\log(e+y)\}.
\]
In the range \eqref{eq:fixed-small-regime},
$z\leq Cn^{-(q-2)/3}$, so $Q\to0$ uniformly over the permitted
dimensions $d$.

If $m_b>z^{-b}$ for some $b\geq1$, apply \eqref{eq:fixed-screen} to
$\lfloor z^{-b}\rfloor$ such faces, all with thresholds at most $\sqrt b$.
By \eqref{eq:fixed-S-compare}, the probability of $S\in A$ is then at
most $C(Q^b+B)\leq C(Q+B)$. This also bounds every outward crossing.

Otherwise $m_b\leq z^{-b}$ for every $b$. The necessary condition for an outward crossing is
$W_l\leq\sqrt b+k\delta$ on these coordinates. Complete the retained array with fresh bounded rows in the omitted positions. Its sum has law $S$, and the added value in every coordinate has absolute value at most $k\delta$. The retained event therefore has probability at most
\[
 F_{m_b}(\sqrt b+2k\delta)+CB.
\]
Since $\log m_b\leq by$ and $k\delta\sqrt{\log(em_b)}=O(1)$, \eqref{eq:fixed-screen} and $r\leq1$ give
\[
 m_bz^bF_{m_b}(\sqrt b+2k\delta)
 \leq Cz^{br}\exp\{c_ab\sqrt y+Cb\log(e+y)\}
 \leq CQ^b,\qquad
 m_bz^bB\leq B.
\]
Multiplying by $\lambda_0^h/(b!h!)$ from \eqref{eq:fixed-mark-weight} and summing $b,h$ proves
\begin{equation}\label{eq:fixed-directional-pre}
 \sup_A\mathbb P\{\mathbf1_{\{T_n\in A\}}\ne\mathbf1_{\{S\in A\}}\}
 \leq C(B+Q),\qquad
 d_{\rm rect}(S,G_R)\leq CB.
\end{equation}
Because
$y=(q-2)\log n/2-(q/2)\log \log(end)-\log K\leq(q-2)\log n/2$,
\[
 Q\leq Cn^{-(q-2)r/2}\log^{qr/2}(end)\log^C(en)
 \exp\left\{\frac{\sqrt{a(q-2)}}{1-a}\sqrt{\log n}\right\}.
\]

\subsubsection{A dimension-independent bound for the rare-observation term}
For $0<r\leq1$, let
\[
 R_*(n)=n^{-(q-2)r/2}\log^{C_0}(en)
 \exp\left\{\frac{\sqrt{a(q-2)}}{1-a}\sqrt{\log n}\right\},
\]
with $C_0$ sufficiently large. First consider the faces with $t_j\leq1$. Choose $D_1$ sufficiently
large that, if either sign group has more than $n^{D_1}$ such faces,
\eqref{eq:fixed-screen} at $t=1$ bounds the Gaussian probability of the
whole rectangle by $O(n^{-2})$. The probability of an outward crossing
is then $O(B+n^{-2})$. Otherwise there are at most $2n^{D_1}$ such
faces, so the logarithm of their number is $O(\log(en))$. Applying
\eqref{eq:fixed-directional-pre} to this subrectangle gives
$O(B+R_*(n))$.

For the remaining selected thresholds, $1<t_j\leq\sqrt b$ implies
$b\geq2$. Repeating the preceding calculation on this group gives
$O(B+Q^2)$. Moreover,
\[
 \log(end)\leq Cn^{(q-2)/(3q)},\quad z\leq Cn^{-(q-2)/3},\quad
 Q^2\leq n^{-2(q-2)r/3+o(1)}=o(n^{-(q-2)r/2}).
\]
If $r>1/q$, the strict inequality $(q-2)r/2>(q-2)/(2q)$ makes $R_*(n)=o(n^{-(q-2)/(2q)})$, so this correction is absorbed by $B$. If $r\leq1/q$, it is the term $E_{q,a}(n)$ in the theorem. This proves \eqref{eq:fixed-sample-upper}. The estimates bound the probability that the two statistics disagree
on each subrectangle of the same sample. A union bound therefore
combines the subrectangles, even when their coordinates are dependent.

\subsubsection{Nearly antipodal pairs: a lower bound for the sample-size rate}
A lower-bound construction in fewer than $d$ coordinates extends to
dimension $d$ by appending independent symmetric signs and leaving those
coordinates unconstrained in the test rectangle. This preserves the
moment bound, the correlation condition, and the discrepancy.
Let $h=n^{-(q-2)/(2q)}$ and $b=n^{1/q}$. Independently of a symmetric sign $U$, let $V$ have masses $\lambda/(2n)$ at $b$ and $-b$, and remaining mass at zero. Choose $\lambda\in(0,1/16)$ sufficiently small that $(1+\lambda^{1/q})^q<K$. The standardized pair
\[
 (X_+,X_-)=(U+V,-U+V)/\sigma,\qquad \sigma^2=1+\lambda h^2,
\]
has centered coordinates with unit variance, $q$th absolute moments
below $K$, and negative correlation $(\lambda h^2-1)/(1+\lambda h^2)$.
For the lower orthant at zero, the sample event reduces to $hJ\leq-|W|$, where $W=n^{-1/2}\sum_iU_i$ and $J=b^{-1}\sum_iV_i$ are independent. For large $n$, $b\geq2$ ensures positive denominators. Berry--Esseen, symmetry of $J$, and $\mathbb E|J|\leq\lambda$ bound the sample probability by
\[
 \mathbb E\{\Phi(h|J|)-1/2\}+O(n^{-1/2})
 \leq\lambda h/\sqrt{2\pi}+O(n^{-1/2}).
\]
The matched Gaussian orthant has probability $\pi^{-1}\arctan(\sqrt\lambda h)\geq\sqrt\lambda h/(2\pi)$. Our choice of $\lambda$ and $n^{-1/2}=o(h)$ give a discrepancy of at least $c_{q,K}h$.

For the logarithmic dimension factor, use $m=\lfloor d/2\rfloor$ independent copies of the pair, let $t=\Phi^{-1}(1-1/m)\geq1$ for large $m$, and constrain each pair by $T_+\leq t$ and $T_-\geq-t$. If $P,Q$ are the one-pair sample and Gaussian probabilities, respectively, then
\[
 P\geq\mathbb P(W\leq t-h|J|),\qquad
 Q=\mathbb E\Phi(\sigma t-h\sqrt\lambda|Z|),\quad Z\sim N(0,1).
\]
The first inequality uses $|\pm U_i+V_i|\geq1$, so that every unscaled denominator is at least $\sqrt n$.
Assume $\log(ed)\leq c_0n^{\min\{(q-2)/q,1/3\}}$ with $c_0>0$ small. Then $t\leq Cn^{1/6}$ and $th\leq C\sqrt{c_0}$. For the symmetric sign sum, we use the following normal-tail
comparison, as in related lower-bound arguments \citep{DL2021}:
uniformly for $y\leq t$,
\begin{equation}\label{eq:fixed-sign-tail}
 \mathbb P(W>y)\leq\{1+C(1+t)n^{-1/2}\}\bar\Phi(y).
\end{equation}
For $y<0$ this follows from Berry--Esseen and $\bar\Phi(y)\geq1/2$. For nonnegative lattice points $z\leq\sqrt{n/2}$, Stirling's inequalities and the binary-entropy bound imply
\[
 \mathbb P(W=z)\leq \frac2{\sqrt n}\phi(z)
 \frac{\exp\{C/n-z^4/(12n)\}}{\sqrt{1-z^2/n}}
 \leq\frac2{\sqrt n}(1+C'/n)\phi(z).
\]
The second inequality follows by using $-\frac12\log(1-u)\leq u$ for $0\leq u\leq1/2$ and maximizing $z^2-z^4/12$. An upper Riemann sum adds at most $2\phi(y)/\sqrt n$ to $\bar\Phi(y)$, and the tail above $\sqrt{n/2}$ is exponentially negligible. Mills' inequalities prove \eqref{eq:fixed-sign-tail}.

The number $N$ of nonzero marks has mean $\lambda$ and satisfies $\mathbb Ee^{vN}\leq\exp\{\lambda(e^v-1)\}$ for fixed $v>0$, while $|J|\leq N$. The normal hazard bound on $(-\infty,t]$ gives $\mathbb E\bar\Phi(t-h|J|)\leq C\bar\Phi(t)$. After averaging, the error contributed by
\eqref{eq:fixed-sign-tail} is therefore $O(\phi(t)n^{-1/2})$. Taylor's inequality
\[
 |\Phi(t+z)-\Phi(t)-z\phi(t)|
 \leq C\phi(t)(t+|z|)z^2e^{t|z|}
\]
yields
\[
 P-Q\geq h\phi(t)\{\sqrt{2\lambda/\pi}-\lambda-O(th)-O(n^{-1/q})\}
 \geq c_{q,K}h\phi(t).
\]
Choose $c_0$ first and $n$ large afterwards. Also $1-Q\leq C/m$. Independence and Mills' inequalities give
$P^m-Q^m\geq mQ^{m-1}(P-Q)\geq c_{q,K}h\sqrt{\log(ed)}$. Bounded $m$ is covered by the two-dimensional construction.

For $2<q<3$, the pair construction already proves the complete lower bound. Indeed, in the range $\log(ed)\leq c_0n^{(q-2)/q}$,
\[
 \frac{n^{-1/2}\log^{3/2}(ed)}
      {n^{-(q-2)/(2q)}\sqrt{\log(ed)}}
 =n^{-1/q}\log(ed)
 \leq c_0n^{(q-3)/q}\leq C.
\]
For larger dimension, use a number of coordinates whose logarithm is
$(c_0/2)n^{(q-2)/q}$ and pad by independent signs. The pair construction in that smaller dimension has a positive constant discrepancy. This proves the claimed lower bound, capped at one, for every $2<q<3$.

\subsubsection{Nonzero skewness and the dimension restriction}
The next bounded example gives the
$n^{-1/2}\log^{3/2}(ed)$ term when $q\geq3$. Its scalar expansion is
valid for every $q>2$.

\begin{lemma}[A bounded law with nonzero third moment]
\label{lem:fixed-skew-lower}
For every $q>2$ and $K>1$, there is a centered bounded law with
variance one and $q$th absolute moment below $K$ such that, for its
product law in dimension $d$, the self-normalized vector has rectangle
error at least
\[
 c_{q,K}n^{-1/2}\log^{3/2}(ed)
\]
whenever $L_0\leq\log(ed)\leq c_0n^{1/3}$, where $L_0$ is sufficiently large and $c_0>0$ is sufficiently small. In particular, uniform convergence over the row class requires $\log d=o(n^{1/3})$.
\end{lemma}
\begin{proof}

Choose a fixed $\theta<1/2$ sufficiently close to $1/2$ that
\[
 X=\frac{I-\theta}{\sqrt{\theta(1-\theta)}},\qquad I\sim\operatorname{Bernoulli}(\theta),
\]
has $\mathbb E|X|^q<K$, and use independent coordinates. Put $\nu=(1-2\theta)/\sqrt{\theta(1-\theta)}>0$. Then $\mathbb EX^3=\nu$ and $X^2=1+\nu X$. If $S=n^{-1/2}\sum_iX_i$, the original statistic exceeds $x\geq0$ exactly when
\[
 S>y=x\sqrt{1+\nu^2x^2/(4n)}+\nu x^2/(2\sqrt n).
\]
Uniformly for $1\leq x\leq Cn^{1/6}$,
\begin{equation}\label{eq:fixed-skew-tail}
 \log\frac{\mathbb P(T>x)}{\bar\Phi(x)}
 =-\frac{\nu x^3}{3\sqrt n}
 +O\!\left(\frac{1+x}{\sqrt n}+\frac{1+x^4}{n}\right).
\end{equation}
The leading cubic term is the classical skewness correction for a
Studentized mean \citep{WANG2005,GSS2022}. For completeness, we derive
the required expansion for this two-point law. Let $K_0(v)=\log\mathbb Ee^{vX}$, put $z=y/\sqrt n$, and solve $K_0'(v)=z$. Analyticity of its cumulant generating function gives
\[
 v=z-\nu z^2/2+O(z^3),\quad
 vz-K_0(v)=z^2/2-\nu z^3/6+O(z^4),\quad
 k=v\sqrt{nK_0''(v)}=y+O(y^3/n).
\]
Under exponential tilting, Berry--Esseen applies uniformly for small $v$. The function $u\mapsto e^{-ku}\mathbf1_{\{u>0\}}$ has total variation two, so its tilted expectation differs by $O(n^{-1/2})$ from $M(k)=e^{k^2/2}\bar\Phi(k)\asymp(1+k)^{-1}$. Returning to the original measure gives
\[
 \mathbb P(S>y)=e^{-n(vz-K_0(v))}M(k)\{1+O((1+k)/\sqrt n)\}.
\]
Since $(\log M)'=O(1)$ on $[0,\infty)$, subtracting $\log\bar\Phi(y)$ gives $\nu y^3/(6\sqrt n)+O((1+y)/\sqrt n+(1+y^4)/n)$. Finally,
$y=x+\nu x^2/(2\sqrt n)+O(x^3/n)$ and $-(\log\bar\Phi)'(x)=x+O((1+x)^{-1})$ prove \eqref{eq:fixed-skew-tail}.

Take $x=\Phi^{-1}(1-1/d)$. When $\log(ed)\leq c_0n^{1/3}$ and $x$ is above a sufficiently large fixed constant, \eqref{eq:fixed-skew-tail} gives a tail deficit at least $c_{q,K}\bar\Phi(x)x^3/\sqrt n$. Independence makes the common-threshold discrepancy at least $c_{q,K}\log^{3/2}(ed)/\sqrt n$. For bounded $x$, the lower bound from the nearly antipodal pair is
larger. If $\log(ed)>c_0n^{1/3}$, use a smaller dimension whose
logarithm is of order $c_0n^{1/3}$, and append independent signs. The
discrepancy is then bounded below by a positive constant.
\end{proof}

The larger of the pair lower bound and Lemma~\ref{lem:fixed-skew-lower} proves \eqref{eq:fixed-universal-lower} for $q\geq3$. The preceding pair argument already proves it for $2<q<3$.

\subsubsection{A rare-observation lower bound for fixed and vanishing gaps}

\begin{proposition}\label{prop:rare-mark-unified}
Fix $q>2$, $K>1$, $0\leq a_0<1$ and $M>0$. There are constants $\eta,B>0$, depending only on $q,K,a_0,M$, such
that, for all sufficiently large $n$, with
\begin{gather*}
 A=B\sqrt n\,\log n,\qquad
 \pi=\eta A^{-q},\qquad v=A^2\pi(1-\pi),\\
 m_n=\left\lceil\frac{M\log n}{n\pi}\right\rceil,\qquad
 t=1-(\log n)^{-1/2},
\end{gather*}
one has, uniformly for $0\leq a\leq a_0$,
\begin{equation}\label{eq:rare-mark-unified}
 \Delta_{n,q,K}(m_n,a)\geq
 \left[
 \frac12\overline\Phi\left(
 \frac{\sqrt{a+v}\sqrt{2\log m_n}-t\sqrt{1+v}}{\sqrt{1-a}}
 \right)-Cn^{-M}\right]_+.
\end{equation}
If $a\log n\to\infty$ and
$M>(q-2)a_0/[2(1-a_0)]+3$, then
\begin{equation}\label{eq:rare-mark-logarithm}
 \log\Delta_{n,q,K}(m_n,a)
 \geq-\frac{(q-2)a\log n}{2(1-a)}
 +\frac{\sqrt{(q-2)a\log n}}{1-a}
 -C\{1+a\log\log n+\log(2+a\log n)\}.
\end{equation}
In particular,
$\log m_n=(q-2)\log n/2+(q+1)\log\log n+O(1)$.
\end{proposition}

\begin{proof}
Take independent symmetric signs $C_i$, independent signs $V_{ij}$ with mean $\sqrt{1-a}$, and independent Bernoulli$(\pi)$ variables $I_{ij}$, all three families being independent. Put
\[
 Y_{ij}=C_iV_{ij},\qquad
 X_{ij}=\frac{Y_{ij}+A(I_{ij}-\pi)}{\sqrt{1+v}}.
\]
Each $Y_{ij}$ is a symmetric sign, and distinct coordinates of $Y_i$
have covariance $1-a$. Thus $X$ is centered with unit marginal variances and off-diagonal covariance
$(1-a)/(1+v)\leq1-a$. Moreover,
\[
 \|X_{ij}\|_q\leq1+A\|I_{ij}-\pi\|_q
 \leq1+2\eta^{1/q}.
\]
Choose $\eta$ so this is less than $K^{1/q}$.

The stated expansion of $\log m_n$ follows from
$A=B\sqrt n\log n$. Hoeffding's inequality and a union bound over these polynomially many coordinates give
\[
 \max_{j\leq m_n}\left|\sum_iY_{ij}\right|
 \leq D\sqrt{n\log n}
\]
except on an event of probability at most $Cn^{-M}$, for a fixed
$D=D(q,M)$. Choose $B>4(D+1)$ and then take $n$ large with these constants fixed.

Suppose a column has $k\geq1$ marks. Its unscaled numerator is at least
$A(k-n\pi)-D\sqrt{n\log n}$, and the Euclidean norm of its observations is at most
$\sqrt n+A\sqrt k+A\sqrt n\,\pi$. Its ratio is therefore at least
\[
 \frac{k-n\pi-D/(B\sqrt{\log n})}
      {\sqrt k+1/(B\log n)+\sqrt n\,\pi}>t.
\]
For the last inequality, use
$k-t\sqrt k\geq1-t=(\log n)^{-1/2}$ and
$n\pi+\sqrt n\,\pi=o((\log n)^{-1/2})$. The numerator is positive for all large $n$.
The probability that the whole array has no mark is at most
$e^{-nm_n\pi}\leq n^{-M}$. Hence $\mathbb P(T_{n,j}\leq t\text{ for all }j)\leq Cn^{-M}$.

The comparison Gaussian can be written as
\[
 G_j=\frac{\sqrt{1-a}\,Z_0+\sqrt{a+v}\,Z_j}{\sqrt{1+v}},
\]
where all the $Z_j$ are independent standard normals. Since
$\mathbb P(\max_{j\leq m_n}Z_j\leq\sqrt{2\log m_n})\geq1/2$, intersecting this event with the appropriate independent half-space for $Z_0$ proves \eqref{eq:rare-mark-unified}.

To obtain \eqref{eq:rare-mark-logarithm}, note that $v\log n=o(1)$ and expand the square of
\[
 y_n=\frac{\sqrt{a+v}\sqrt{2\log m_n}-t\sqrt{1+v}}{\sqrt{1-a}}.
\]
As $a\log n\to\infty$, this argument is eventually positive. The negative quadratic contribution to $-\tfrac12y_n^2$ is bounded below by
\[
 -\frac{(q-2)a\log n}{2(1-a)}-C(1+a\log\log n),
\]
the positive cross term is bounded below by
$\sqrt{(q-2)a\log n}/(1-a)-C$, and the remaining term is bounded below by a constant. Mills' lower bound contributes an additional term of magnitude at most
$C\log(2+a\log n)$ to the logarithm. The assumption on $M$ makes the subtracted sample probability negligible compared with this Gaussian lower bound.
\end{proof}

For fixed $0<a\leq a_c$, Proposition~\ref{prop:rare-mark-unified} gives \eqref{eq:fixed-rare-lower}. The upper bound \eqref{eq:fixed-sample-upper} gives the matching logarithmic expansion \eqref{eq:fixed-endpoint-log}, because $(q-2)r/2\leq(q-2)/(2q)$ and the rare-observation term dominates $B_q(n,m_n)$, including when
$r=1/q$.
Finally,
\[
 B_q(n,d)\leq2n^{-(q-2)/(2q)}\log^2(end).
\]
If $a>a_c$, the upper bound attains exponent $(q-2)/(2q)$, and the two-dimensional pair excludes every larger exponent. If $a\leq a_c$, the lower bound at $m_n$ rules out the endpoint exponent $(q-2)r/2$ with every fixed logarithmic factor. For every smaller positive exponent, the extra power of $n$ bounds
the factor $\exp\{C\sqrt{\log n}\}$. This completes the proof of Theorem~\ref{thm:fixed-sample-rate}.

\subsection{Convergence throughout \texorpdfstring{$\log d=o(n^{1/3})$}{log d = o(n\string^(1/3))}}
\label{sec:fixed-cubic}
To extend convergence to $\log d=o(n^{1/3})$, we combine Gaussian
boundary counts with replacement of one observation at a time. Put
\[
 w=n^{-1/6}\sqrt{\log(ed)},\qquad F(w)=w[1+\log_+(1/w)]^3.
\]
\begin{theorem}\label{thm:fixed-complete-envelope}
For fixed $q\geq3$, $K>1$ and $a>0$,
\begin{equation}\label{eq:fixed-complete-envelope}
 \Delta_{n,q,K}(d,a)\leq1\wedge C\left[
 \min\{F(w),B_q(n,d)\}+E_{q,a}(n)\right].
\end{equation}
Together with Lemma~\ref{lem:fixed-skew-lower}, this yields
\eqref{eq:main-cubic-criterion} for fixed $a>0$.
\end{theorem}
The Gaussian estimates and the comparison for bounded observations
below are uniform over deterministic sets of deleted rows.

\subsubsection{Gaussian boundary counts and smooth approximations}
Let $g=-\log\bar\Phi$. Mills' inequalities give
$0<g'(x)\leq C(1+x_+)$, $g'(x)>x$ for $x>0$, and $0<g''(x)\leq1$. The latter bound follows also from
$\operatorname{Var}(Z\mid Z>x)=1-g''(x)$. Gaussian vectors in this subsection have unit marginal variances and may be singular.
For finite thresholds $b_j$, set
$z_j(x)=g(x_j)-g(b_j)$ and $Y(x)=\max_jz_j(x)$.
Lemma~\ref{lem:score-strip} gives the dimension-free strip bound
\begin{equation}\label{eq:cubic-score-ac}
 \mathbb P\{s<Y(G)\leq t\}\leq t-s,\qquad s<t.
\end{equation}
We next count coordinates whose transformed values are close to their
transformed thresholds. The expectation retains the full rectangle
event even when only a subset of the constraints is counted.

\begin{lemma}[Gaussian boundary counts at a fixed correlation gap]\label{lem:cubic-face-count}
Let $\mathcal C$ be a Gaussian polyhedron. Select a group of its constraints
$G_j\leq b_j$ such that distinct selected coordinates have correlations at most
$1-a$, where $a>0$. For
\[
 0<\alpha<\min\{1,a/(1-a)\},
\]
and a sufficiently large $\lambda=\lambda(a,\alpha)$,
\begin{equation}\label{eq:cubic-fractional}
 \mathbb E\left[\mathbf1_{\mathcal C}
 \left(\sum_j e^{-\lambda[g(b_j)-g(G_j)]}\right)^\alpha\right]
 \leq C_{a,\alpha}.
\end{equation}
When $a=1$, decrease the gap slightly before choosing $\alpha$.

If all selected thresholds satisfy $b_j\geq4$ and
\begin{equation}\label{eq:cubic-high-counting}
 \#\{j:g(b_j)\leq t\}\leq2e^{D_1t},\qquad t>0,
\end{equation}
then \eqref{eq:cubic-fractional} holds with exponent one for a
sufficiently large $\lambda=\lambda(a,D_1)$, with a bound depending
only on $a,D_1$. More precisely, if
\[
 N_y=\#\{j:0\leq g(b_j)-g(G_j)\leq y\},
\]
then
\begin{equation}\label{eq:cubic-high-strip}
 \mathbb E[\mathbf1_{\mathcal C}N_y]
 \leq\mathcal H(y):=
 \begin{cases}
 Cy,&0<y\leq y_0,\\
 C_2e^{C_2y},&y>y_0.
 \end{cases}
\end{equation}
Here $y_0>0$, $C$, and $C_2$ depend only on $a,D_1$. These
conclusions persist when selected threshold scores are increased.
They apply separately to the two sign groups, without a restriction
on correlations across the groups.
\end{lemma}
\begin{proof}
For the count bounds, represent $G_j=\langle v_j,Z\rangle$ with $\|v_j\|=1$. To prove \eqref{eq:cubic-fractional}, fix $y>0$ and split the faces at $B=4\sqrt{1+y}$. If $m$ thresholds are at most $B$, Slepian's inequality \citep{Slepian1962} and splitting the independent Gaussian maximum at $u=\sqrt{2(1-\xi)\log m}$ give
\[
 \mathbb P(\mathcal C)\leq \Phi(u)^m+
 \Phi\left(\frac{B-\sqrt a\,u}{\sqrt{1-a}}\right)
 \leq Ce^{Cy}m^{-\alpha},
\]
where $\xi>0$ is chosen so that $(1-\xi)^2a/(1-a)>\alpha$. Thus the expected $\alpha$th power of the count for thresholds at
most $B$ is bounded by $Ce^{Cy}$.
For a constraint with $b_j>B$ and score gap at most $y$, $t_j=b_j-\langle v_j,Z\rangle$ obeys $0\leq t_j\leq2y/b_j$. Put $c_j=8y/(ab_j)$ and translate this strip intersected with $\mathcal C$ by $c_jv_j$. Writing $\rho=1-a$ and
$H_j(z)=c_j(\langle v_j,z\rangle-b_j)-c_j^2/2$, the image of the $j$th strip satisfies
\[
 H_j-H_k\geq\tfrac12(c_k-\rho c_j)^2+
 c_j^2\{(1-\rho^2)/2-t_j/c_j\}>0.
\]
Thus all images are disjoint. Their Gaussian density ratio is at least
$\exp(-c_jb_j-c_j^2/2)\geq e^{-C_ay}$. This bounds the first moment of the count for high thresholds by
$e^{C_ay}$. Subadditivity and $m^\alpha\leq m$ for integer $m\geq0$ now give
$\mathbb E[\mathbf1_{\mathcal C}N_y^\alpha]\leq Ce^{Cy}$. Partitioning the score gaps into successive intervals and summing
with exponential weights proves \eqref{eq:cubic-fractional}.

Under \eqref{eq:cubic-high-counting}, faces with $g(b_j)\leq2y$ number at most $2e^{2D_1y}$. For every remaining constraint with score gap at most $y$,
$g(G_j)\geq g(b_j)/2$ implies $G_j>b_j/2$, since
$g(b_j/2)<g(b_j)/2$ for $b_j\geq4$. The same translation gives an expected count at most $e^{C_ay}$. To obtain the linear bound for $y\leq a/8$, use $c_j=1/b_j$, $\vartheta=a/8$, and map a strip point $x$ to
\[
 x+c_ju(x)v_j,\qquad
 u(x)=1+\vartheta\left(1-\frac{b_j-\langle v_j,x\rangle}{2yc_j}\right).
\]
The map is injective with Jacobian $1+\vartheta/(2y)$. Its image is again strictly separated by the functions $H_j$: the relevant lower bound is
$c_j^2[u-1/2-t_j/c_j-\rho^2u^2/2]>0$ for $1\leq u\leq1+\vartheta$. The density ratio is bounded below by a positive constant. Disjointness therefore gives $\mathbb E[\mathbf1_{\mathcal C}N_y]\leq C_ay$. This proves \eqref{eq:cubic-high-strip}. Integrating against
$\lambda e^{-\lambda y}\,dy$ then proves
\eqref{eq:cubic-fractional} with exponent one.
\end{proof}

Define the flat cutoff
\[
 h(t)=\frac{e^{-1/t}}{e^{-1/t}+e^{-1/(1-t)}}\quad(0<t<1),
 \qquad h(t)=0\ (t\leq0),\quad h(t)=1\ (t\geq1).
\]
For score vector $z$, let
$A(z)=1+\sum_je^{\lambda z_j}$,
$B_*(z)=\sum_je^{\lambda z_j}h(z_j)$, and
$\mathcal R(z)=A(z)/(1+B_*(z))$.
Choosing $\lambda$ so that the preceding weighted bounds hold also at $\lambda/2$ gives
\begin{equation}\label{eq:cubic-gaussian-ratio}
 \mathbb E\mathcal R(z(G))^\alpha\leq C.
\end{equation}
On $Y\leq1$, apply \eqref{eq:cubic-fractional} after increasing scores of all thresholds by one. On $k\leq Y<k+1$, $k\geq1$, one has $B_*\geq e^{\lambda k}$ and
$A-1-B_*\leq e^{\lambda/2}\sum_je^{\lambda z_j/2}$.
Apply the weighted bound on the orthant shifted by $k+1$ and sum the resulting geometric sequence. The constant terms contribute at most one because the slices are
disjoint. The same proof gives \eqref{eq:cubic-gaussian-ratio} with exponent one under the high-threshold hypotheses, at every upward threshold shift.

For $N\asymp n$ independent Gaussian rows, let $T^G$ be the
self-normalized vector and $Z=N^{-1/2}\sum_iG_i$. Gaussian norm concentration and Cauchy--Schwarz give
\begin{equation}\label{eq:cubic-G-calibration}
 \mathbb E\min\{1,\max_j|g(T_j^G)-g(Z_j)|\}
 \leq C\frac{\log^{3/2}(ed)}{\sqrt N}=:u_N.
\end{equation}
To verify this estimate, put
$D_N=\max_j|\|G_{\cdot j}\|_2/\sqrt N-1|$. The required moment bounds are
$\|D_N\|_2\leq C\sqrt{\log(ed)/N}$ and
$\|(1+\max_j|Z_j|)^2\|_2\leq C\log(ed)$. On $D_N\leq1/2$ the score distortion is bounded by $C(1+\max_j|Z_j|)^2D_N$; on its complement, the score difference capped at one is at most $2D_N$.

For a width $0<\delta<1/8$, set
\[
 B_t(z)=\sum_je^{\lambda z_j}h(z_j/\delta+t),\quad t=0,1,
 \qquad\chi(t)=1-h(t-1),
\]
\[
 H^+(z)=1-h\left(\frac{2MB_0(z)}{M+A(z)}\right),\qquad
 H^-(z)=\chi(A(z)/M)[1-h(B_1(z))],\quad M\geq4.
\]
These functions satisfy
$H^-\leq\mathbf1_{\{Y\leq0\}}\leq H^+$. The expectation of either
function at $T^G$ differs from the Gaussian rectangle probability by at
most
\begin{equation}\label{eq:cubic-bracket-calibration}
 C(\delta+M^{-\alpha}+u_N/\delta).
\end{equation}
To check this, the scores of $T^G$ and $Z$ differ by at most
$\delta$ except on an event of probability at most $u_N/\delta$. The Gaussian strip $|Y(Z)|\leq2\delta$ has probability $O(\delta)$. Far inside the orthant, the lower bracket equals one unless $A$ exceeds a fixed multiple of $M$, controlled by \eqref{eq:cubic-fractional}. Far outside, positivity of the upper bracket forces $B_0\geq1$ and $A\geq MB_0$. At thresholds shifted up by $2\delta$, this implies $\mathcal R(z(Z)-2\delta\mathbf1)\geq cM$, controlled by \eqref{eq:cubic-gaussian-ratio}. For an ordinary Gaussian vector, the same argument omits the term
$u_N/\delta$. Applying \eqref{eq:cubic-G-calibration} at width one,
with shift $1/4$, also gives
$\mathbb P(\mathcal R(z(T^G))>cM)\leq C(M^{-\alpha}+u_N)$.

\subsubsection{Replacing bounded observations and solving the error recursion}
\label{subsec:cubic-bounded}
Fix an ambient scale $\Lambda\geq\log(en)$ with $\log(ed)\leq\Lambda$. Let
$w_0=n^{-1/6}\sqrt\Lambda$ and $u=(w_0)^3$; below $w_0$ is sufficiently small. Consider independent centered rows with unit marginal variances and common covariance $R$ having the within-sign gap $a$. Each row is either bounded coordinatewise by $B_0n^{1/3}$ with
uniformly bounded marginal third absolute moments, or is a centered
Gaussian row with covariance $R$. For every $n/2\leq m\leq n$, the
self-normalized vector based on these $m$ rows satisfies
\begin{equation}\label{eq:cubic-bounded-bound}
 d_{\rm rect}(T_m,G_R)\leq Cw_0\log^2(e/w_0).
\end{equation}
All constants are uniform over subcollections and deterministic deleted-row sets, using the same ambient $\Lambda$.

Choose $\alpha$ as in \eqref{eq:cubic-fractional}, the corresponding fixed $\lambda$, and
$\delta=D_0w_0\log^2(e/w_0)$, $M=\max\{4,\delta^{-1/\alpha}\}$.
Increase $D_0$ and then reduce the permitted $w_0$ so that $M^{-\alpha}=\delta$, $\delta<1/8$, and $\log(eM)\leq C_\alpha \log(e/w_0)$.
For a retained sample of $N\asymp n$ rows, write
\[
 Q_j=\sum_iX_{ij}^2,\quad \tau_j=\frac{\sum_iX_{ij}}{\sqrt{Q_j}},\quad
 s_j=1+|\tau_j|,\quad \xi_j=\frac{s_jy_j}{\sqrt{Q_j}},\quad
 x_j(v)=\frac{\tau_j+v/s_j}{\sqrt{1+(v/s_j)^2}}.
\]
Insertion of a row $y$ gives exactly $x_j(\xi_j)$.
Except on an event of probability $O(n^{-20})$, the remaining sample
satisfies $Q_j\geq N/2$ for all $j$ and
$\max_j|\tau_j|\leq C\sqrt\Lambda$. We verify this using only
independence of the rows and the stated moment bounds. Truncating each square at a sufficiently large fixed constant, chosen using the common third absolute moment bound, leaves mean at least $3/4$. Hoeffding's inequality controls the lower tail of $Q_j$. Writing
$S_j=\sum_iX_{ij}$ for the remaining numerator,
$e^{z-z^2/2}\leq1+z+Cz^2$ implies
$\mathbb Ee^{tS_j-t^2Q_j/2}\leq e^{CNt^2}$. On a shell $q_0\leq Q_j<2q_0$, $q_0\geq N/2$, take $t=cx/\sqrt{q_0}$ to bound the event $S_j>x\sqrt{Q_j}$ by $e^{-cx^2}$. Sum dyadic shells up to index $\lceil x^2\rceil$, and use $\mathbb EQ_j=N$ for the remaining tail. A sign change and a union bound finish the claim. These estimates are uniform over the mixtures of rows considered
here.

On this event for the remaining sample, an incoming bounded row has
$\|\xi\|_\infty\leq CB_0w_0\leq c\delta/\log^2(eM)$. For an incoming Gaussian row, the same bound holds on
$\{\max_j|y_j|\leq B_1n^{1/3}\}$. The complementary probability is
at most $2d e^{-cn^{2/3}}$, which is negligible since
$\Lambda=O(n^{1/3})$. The first two Taylor coefficients and the fixed-order moments of the
Gaussian coordinate maximum grow at most polynomially in $n$.
Cauchy--Schwarz therefore bounds their contribution on the exceptional
Gaussian event by $O(n^{-20})$.

We first fix the weights used to bound the derivatives. Set
$m_0=M^{-2}$ and, at the retained score vector, define
\[
 a_j=e^{\lambda z_j}/A,\qquad
 b_j=\frac{e^{\lambda z_j}[h(z_j/\delta+t)+m_0]}{B_t+m_0A}
\]
for $t=0,1$. Each vector has total mass at most one and is independent
of the incoming row. For any such nonnegative vector $b$, let
$\operatorname{Diag}_r(b)$ be the order-$r$ tensor with entry $b_j$
when all indices equal $j$, and zero otherwise. Tensor products, index
permutations, and diagonal lifts preserve the bound of one on total
mass. Here a diagonal lift adds indices constrained to equal existing
ones.

We next estimate the derivatives along the insertion path
$v=\theta\xi$, $0\leq\theta\leq1$. The first three derivatives of
$g(x_j(v_j))$ are uniformly bounded for $|v_j|\leq1$, by the insertion
formula, $s_j=1+|\tau_j|$, $g'(x)\leq C(1+x_+)$, and the bounds on
$g''$ and $g'''$. Direct differentiation of the flat cutoff gives
\begin{equation}\label{eq:cubic-flat-estimate}
 |h^{(k)}(t+v)|\leq C_k\log^{2k}(eM)[h(t)+m_0],\qquad
 |v|\leq c/\log^2(eM),\quad0\leq k\leq3.
\end{equation}
Indeed, $|h^{(k)}|\leq C_kh\log^{2k}(e/h)$ and
$|(h+m_0)'|/(h+m_0)\leq C\log^2(eM)$; integrate the second bound
before applying the first. Consequently, along the path, $A(v)$ and
$B_t(v)+m_0A(v)$ are comparable to their retained values. Here the
argument $v$ means evaluation at the inserted score vector. The
derivatives of $A$ are bounded entrywise by nonnegative diagonal
tensors of mass $CA$; the order-$k$ derivatives of $B_t$ have a bound
of mass $C_k\delta^{-k}\log^{2k}(eM)(B_t+m_0A)$.

Applying these bounds to $H^+$ and $H^-$ produces products and
permutations of the normalized tensors above. For $H^+$, a nonzero
derivative requires $2MB_0(v)/(M+A(v))\leq1$, which gives
$M(B_0+m_0A)/(M+A)\leq C$. For $H^-$, a nonzero product derivative
requires $A(v)\leq2M$ and $B_1(v)\leq1$, so $B_1+m_0A\leq C$.
Terms differentiating only $\chi(A/M)$ have total mass $O(1)$.

The factor $\delta^{-3}$ is needed only on a small event. If the
maximum score of the remaining sample is less than $-C\delta$, the
upper bracket is constant. If that maximum is greater than $C\delta$
and the upper bracket varies, then $B_0(v)\geq1$ and
$A(v)\geq MB_0(v)$. After shifting all threshold scores up by
$C\delta$, the preceding comparison of the sums implies
$\mathcal R>cM$. Outside the boundary strip, the lower bracket either
vanishes or reduces to the cutoff $\chi(A/M)$.

It follows that a nonnegative tensor determined by the remaining
sample bounds the absolute third derivative entrywise along the whole
insertion path. Its total mass is at most
\begin{equation}\label{eq:cubic-third-majorant}
 C\{1+\delta^{-3}\log^6(eM)\mathbf1_{\mathcal V}\},
\end{equation}
where $\mathcal V$ is the union of the $C\delta$ boundary strip for
the retained scores and the event $\mathcal R>cM$ at the shifted
thresholds. Both events depend only on the remaining sample.
The linear and quadratic Taylor terms cancel conditionally: incoming bounded and Gaussian rows have identical means and covariances. In original coordinates their polynomial is exactly
\[
 F(\tau)+\sum_jF_j(\tau)\frac{y_j}{\sqrt{Q_j}}
 +\frac12\sum_{jk}F_{jk}(\tau)\frac{y_jy_k}{\sqrt{Q_jQ_k}}
 -\frac12\sum_j\tau_jF_j(\tau)\frac{y_j^2}{Q_j}.
\]
Thus the denominator's quadratic contribution is included. Conditioning on the remaining sample, H\"older's inequality bounds
each mixed third absolute moment by the marginal third absolute moment
bound. The expected Taylor remainder for one replacement is therefore at
most $Cn^{-3/2}\Lambda^{3/2}$ times the expectation of
\eqref{eq:cubic-third-majorant} over the remaining sample.

We can compare the event involving $\mathcal R$ directly, without
using a bound for the rectangle distance. Majorize
$\mathbf1_{\{\mathcal R\geq cM\}}$ by $\chi(cM/\mathcal R)$ and apply the same derivative calculation at width one. Where a derivative is nonzero,
$M(1+B_*+m_0A)/A\leq C_c$. The absolute third derivative is bounded
entrywise by a tensor determined by the remaining sample, with total
mass $C_c\log^6(eM)$. Summing over all row replacements gives error $Cu\log^6(eM)$.
The Gaussian endpoint is bounded by \eqref{eq:cubic-gaussian-ratio} and \eqref{eq:cubic-G-calibration}. Thus
\begin{equation}\label{eq:cubic-ratio-transfer}
 \mathbb P(\mathcal R(z(T_m))>cM)\leq C(M^{-\alpha}+u\log^6(eM)).
\end{equation}
This also holds with $M'$ in place of $M$ whenever $M'$ is bounded
by a fixed power of $n$ and $w_0\log^2(eM')$ is sufficiently small. Under the high-threshold hypotheses, exponent one replaces $\alpha$.

Keep $n,d,\Lambda$ and the bound on each row fixed, and let $D_m$ be
the supremum of the rectangle errors over the auxiliary class with
exactly $m$ rows. The boundary strip in
\eqref{eq:cubic-third-majorant} has probability at most
$C\delta+2D_{m-1}$, by comparison with two shifted Gaussian rectangles
and \eqref{eq:cubic-score-ac}. Summing the replacement errors and using
\eqref{eq:cubic-ratio-transfer} and
\eqref{eq:cubic-bracket-calibration} gives, uniformly for
$\lceil n/4\rceil+1\leq m\leq n$,
\[
 D_m\leq C[\delta+M^{-\alpha}+u/\delta+u+
 u\delta^{-3}\log^6(eM)(\delta+D_{m-1}+M^{-\alpha}+u\log^6(eM))]+O(n^{-10}).
\]
Since $u/\delta\leq C\delta$, $u\log^6(eM)\leq C\delta$, and
$u\delta^{-3}\log^6(eM)\leq C_\alpha D_0^{-3}$, choose $D_0$ so the coefficient of $D_{m-1}$ is at most $1/2$. Starting from $D_{\lceil n/4\rceil}\leq1$ gives $D_m\leq C\delta+O(2^{-n/4})$ for $n/2\leq m\leq n$. This proves \eqref{eq:cubic-bounded-bound}.

\subsubsection{Covariance comparison for the truncated sample}
Suppose correlation matrices $R,R'$ both have within-sign gap at least $a/2$, and
$d_\Sigma=\max_{jk}|R_{jk}-R'_{jk}|\leq Cn^{-1/3}$. For the ambient $\Lambda$ above,
\begin{equation}\label{eq:cubic-covariance-comparison}
 d_{\rm rect}(G_R,G_{R'})\leq Cw_0\log^2(e/w_0).
\end{equation}
Use the same brackets with $\delta=D_0w_0\log^2(e/w_0)$ and $M=\delta^{-1/\alpha}$. The derivative calculation at second order, now directly in score coordinates, gives
\[
 \|D_z^2H^\pm\|_1+\|D_zH^\pm\|_1
 \leq C\{1+\delta^{-2}\log^4(eM)\mathbf1_{\mathcal E}\},
\]
where $\mathcal E$ is the union of a $C\delta$ score strip and a shifted fixed-width ratio event. Along the Gaussian covariance interpolation, \eqref{eq:cubic-score-ac} and \eqref{eq:cubic-gaussian-ratio} give
$\mathbb P(\mathcal E)\leq C(\delta+M^{-\alpha})$. Since $g'\leq C(1+|x|)$ and $|g''|\leq1$, the chain rule yields
\[
 \mathbb E\|D_x^2[H^\pm(z(G))]\|_1
 \leq C\Lambda\{1+\delta^{-2}\log^4(eM)(\delta+M^{-\alpha})\}+O(n^{-10}).
\]
To justify the tail term, restrict first to $\max_j|G_j|\leq C_0\sqrt\Lambda$. Gaussian coordinate tails, $\Lambda\geq\log(en)$, and the polynomial bounds on $M,\delta^{-1}$ absorb the complement after taking $C_0$ large. Gaussian integration by parts gives the interpolation derivative
$\frac12\sum_{jk}(R'_{jk}-R_{jk})\mathbb E\partial_{jk}F(G_t)$. The convex interpolation preserves the gap. An independent Gaussian perturbation justifies the formula for singular endpoints. Compare each endpoint with its rectangle probability using
\eqref{eq:cubic-bracket-calibration}, omitting the term for Gaussian
self-normalization. Because
$d_\Sigma\Lambda\leq C(w_0)^2$ and $d_\Sigma\Lambda \log^4(eM)\leq C\delta^2$, the total error is $O(\delta)$, proving \eqref{eq:cubic-covariance-comparison}.

\paragraph{Application to the truncated sample.}
For the original observations with $q\geq3$, use the cutoff $\tau=2(Kn)^{1/q}$ and the notation $Y,\mu$ of \Cref{sec:fixed-rates}; write $U$ for the self-normalized vector formed from the centered
truncated rows $Y_i-\mu$. Let $R'$ be the correlation matrix of $Y_1$. Hölder's inequality gives, for any $j,k,l$,
\[
 \mathbb E[|X_jX_k|;|X_l|>\tau]
 \leq\|X_j\|_q\|X_k\|_q\mathbb P(|X_l|>\tau)^{1-2/q}.
\]
Hence $\max_{jk}|R_{jk}-R'_{jk}|\leq C\tau^{2-q}\leq Cn^{-1/3}$, and the truncated gap is at least $a/2$ for large $n$. Rescaling each coordinate of $Y-\mu$ to variance one leaves $U$
unchanged. The standardized observations are bounded by $Cn^{1/3}$
and have uniformly bounded third absolute moments, so
\eqref{eq:cubic-bounded-bound} applies.

For $d\leq n$, the sample-rate proof already gives
$C(n^{-(q-2)/(2q)}\log^2(en)+E_{q,a}(n))$, and $n^{-(q-2)/(2q)}\log^2(en)\leq CF(w)$ for large $n$. The latter follows from $(q-2)/(2q)\geq1/6$, $w\geq cn^{-1/6}$, and monotonicity of $w\log^3(e/w)$ on a sufficiently small fixed interval. If $w$ is larger than a small constant, the asserted bound follows by increasing $C$.
Thus assume $d>n$, put $\Lambda=\log(ed)$, and set
\[
 v_q=n^{-(q-2)/(2q)}\sqrt{\log(ed)}\leq w,\qquad D=C_0w\log^2(e/w),\qquad u=w^3.
\]
For $w$ sufficiently small, the truncated statistic has rectangle
error at most $D$ relative to $G_{R'}$. The same bound holds for each
deterministic choice of at least $n/2$ retained rows. Equation \eqref{eq:cubic-covariance-comparison} adds at most $CD$ to
the error. The same constants apply to every subset of coordinates,
with the original $\log(ed)$ retained in the bound.

\subsubsection{Boundary counts for the truncated statistic at high thresholds}
First assume that all selected signed thresholds satisfy
$b_j\geq A_0\sqrt{\log(en)}$ and that
\eqref{eq:cubic-high-counting} holds in each sign group. We remove these
restrictions at the end of the proof. Write
$\mathcal C_t(x)=\{\max_j z_j(x)\leq t\}$ and
$N_s(x;t)=\#\{j:t-s<z_j(x)\leq t\}$.
For the Gaussian vector with correlation matrix $R'$, \eqref{eq:cubic-high-strip} gives
$\mathbb E[\mathbf1_{\mathcal C_t}N_s]\leq\mathcal H(s)$ at every upward shift $t\geq0$.
To obtain the corresponding bound for a statistic $V$ based on the
remaining truncated observations, fix $M\geq1$. Independently for each
constraint, replace $z_j\leq t$ by $z_j\leq t-s$ with probability
$1/M$. Average the difference between
the probabilities of the original rectangle and the tightened rectangle. This difference is
$\mathbb E[\mathbf1_{\mathcal C_t}(1-(1-1/M)^{N_s})]$ and is at least
$(1-e^{-1})M^{-1}\mathbb E[\mathbf1_{\mathcal C_t}\min(N_s,M)]$. Each of these rectangle probabilities differs from its Gaussian
counterpart by at most $D$, while their Gaussian difference is at most
$\mathcal H(s)/M$. Therefore
\begin{equation}\label{eq:cubic-QH1}
 \mathbb E[\mathbf1_{\mathcal C_t(V)}\min\{N_s(V;t),M\}]
 \leq C[\mathcal H(s)+MD].
\end{equation}
A nonpositive tightened threshold score gives the empty rectangle;
limits from below handle atoms. In particular, with
$m_s=\max\{1,\mathcal H(s)/D\}$, the event
$\mathcal C_t(U),N_s(U;t)>m_s$ has probability $O(D)$.

We also need a bound that holds simultaneously for all large strip
widths. Put $z'_j=z_j-t$,
$A_t=1+\sum e^{\lambda z'_j}$ and $B_t^*=\sum e^{\lambda z'_j}h(z'_j)$. The first-moment form of \eqref{eq:cubic-gaussian-ratio}, and the fixed-width comparison \eqref{eq:cubic-ratio-transfer}, give
\begin{equation}\label{eq:cubic-QH2}
 \mathbb P\{A_t(U)/(1+B_t^*(U))>D^{-1}\}
 \leq C[D+u\log^6(e/D)]=O(D).
\end{equation}
The hypotheses of that comparison hold because $D^{-1}$ is bounded
by a fixed power of $n$ and $w\log^2(e/D)$ is small. On $\mathcal C_t(U)$ one has $B_t^*=0$, so outside the exceptional event, simultaneously for all $s\geq0$,
\begin{equation}\label{eq:cubic-QH3}
 \mathcal C_t(U)\quad\Longrightarrow\quad
 \#\{j:z_j(U)\geq-s\}\leq D^{-1}e^{\lambda(t+s)}.
\end{equation}
Only $t=0,2$ is needed in this implication.

\subsubsection{Restoring large observations after deleting their rows}
Except on an event of probability $O(n^{-10})$, every centered
truncated column has squared norm at least $cn$,
$\|U\|_\infty\leq C\sqrt{\log(ed)}$, and each original column has at
most $C\log(ed)$ entries exceeding $\tau$ in absolute value. Use the
bounds for the remaining columns in the proof of
\eqref{eq:cubic-bounded-bound} and the binomial bound with mean at most
$\lambda_0=2^{-q}$.
Deleting any $k\leq C\log(ed)$ bounded rows keeps the squared norm above $c'n$, because $k\tau^2/n\leq Cv_q^2$ and $w$ is small. The truncated statistic changes by at most
$C(k\tau/\sqrt n+\sqrt{\log(ed)}\,k\tau^2/n)\leq Ckn^{-(q-2)/(2q)}$, and its scores change by at most
\[
 h_k=Ckv_q.
\]

For a coordinate $j$ that crosses its threshold, delete the $k$ rows
where $|X_{ij}|>\tau$, and write $V=U^{(-I)}$. In the remaining column, the original observations were not truncated
and differ from their centered versions only by $\mu_j$; their statistics differ by $O(n^{-(q-2)/(2q)})$ and their scores by $O(v_q)$. Put $v_* =\sqrt{n/\log(ed)}$. If the largest removed entry is in $(b/2,b]$ with $\tau<b\leq2v_*$, set $y=b\sqrt{\log(ed)/n}\leq2$. Reinserting all $k$ entries changes the statistic by at most
$C(kb/\sqrt n+\sqrt{\log(ed)}\,kb^2/n)\leq Ckb/\sqrt n$. The scores of the original statistic, the statistic after deletion,
and the full truncated statistic therefore differ by at most $s=Cky$.
Choose the constant so that $h_k\leq s$. Every intermediate ratio is
$O(\sqrt{\log(ed)})$.
If a deleted entry exceeds $v_*$, set $s=Ck$. The score of the
original statistic is at most the score after deletion plus $s$. This follows by inserting one row at a time and using
\[
 g\left(\frac{x+v}{\sqrt{1+v^2}}\right)\leq g(x)+C.
\]
For $x\leq0$ the new ratio is at most one; for $x>0$ it is at most $\sqrt{x^2+1}$, and the hazard bound proves the claim.
For a fixed selected coordinate, summing over the possible deleted
row indices gives the probability bounds
\begin{equation}\label{eq:cubic-QH4}
 \begin{cases}
 Cv_q^q y^{-q}\lambda_0^{k-1}/(k-1)!,&\text{a moderate bin},\\
 Cv_q^q\lambda_0^{k-1}/(k-1)!,&\text{at least one entry above }v_*.
 \end{cases}
\end{equation}
At one deleted position the observation exceeds $b/2$ or $v_*$ in
absolute value; at every other deleted position it exceeds $\tau$. After dropping the requirement that no further entry in this column
exceeds $\tau$, these conditions involve only deleted rows and are
independent of the remaining sample.

\paragraph{Outward crossings.}
The strip $|\max_jz_j(U)|\leq CD$ has probability $O(D)$ by
\eqref{eq:cubic-score-ac} and the comparison for bounded rows. Outside this
strip, an outward crossing from $U$ to the original statistic requires
a selected coordinate with $k\geq1$ and $s\geq cD$. Choose a fixed $s_0>0$ with $8s_0\leq y_0$. On
$\mathcal C_0(U)$ impose \eqref{eq:cubic-QH3}, and also impose
$N_b(U;0)\leq m_b$ at dyadic widths from $c'D$ through $4s_0$,
including the next dyadic endpoint. Choose $c'$ sufficiently small.
There are $O(\log(e/w))$ such count bounds. A union bound shows that
they fail with total probability $O(D\log(e/w))$.
Fix $k$ and an interval for the magnitude of the largest deleted
entry. The remaining sample satisfies $\mathcal C_{h_k}(V)$, and the
selected score after deletion is greater than $-s$, because the
original statistic strictly violates its upper constraint. Every
candidate with score greater than $-s$ after deletion had score greater
than $-s-h_k$ in $U$.
For $cD\leq s\leq s_0$, let $b(k,s)$ be the smallest width in
this dyadic grid strictly larger than $s+h_k$. Define the deterministic
upper bound
\[
 \operatorname{cap}_{k,s}=
 \begin{cases}
 m_{b(k,s)},&cD\leq s\leq s_0,\\
 D^{-1}e^{\lambda(s+h_k)},&s>s_0.
 \end{cases}
\]
It depends only on $k$, the magnitude bin through $s$, and the fixed
parameters and scales already chosen. It does not depend on the
observations in either the deleted or the remaining rows. The count
bounds for $U$ imply that the number of candidates is at most
$\operatorname{cap}_{k,s}$. Keep this condition on the remaining sample
and discard the other conditions involving the full sample. For fixed
$I$, the sum of the candidate indicators is at most
$\mathbf1_{\mathcal C_{h_k}(V)}
\min\{N_{s+h_k}(V;h_k),\operatorname{cap}_{k,s}\}$.
Here $h_k\leq s$, and, in the first case, $b(k,s)\leq4s\leq y_0$.
Thus \eqref{eq:cubic-QH1} bounds its expectation by $O(s)$ for
$cD\leq s\leq s_0$ and by $O(e^{Cs})$ otherwise. The event and the count in this expectation depend only on the
remaining sample. We can therefore multiply by the probability bounds
for the deleted rows before summing over candidate coordinates.
If the selected coordinate has its largest deleted entry below
$2v_*$ and $Cky\geq cD$, there are at most $C(\log(e/w)+\log(ek))$ relevant dyadic bins for this $k$. At small widths,
$kv_q^qy^{1-q}\leq Ck^qv_q^qD^{1-q}$; at larger widths, $y^{-q}\leq Ck^q$ and $e^{Cs}\leq e^{Ck}$. Multiplying by \eqref{eq:cubic-QH4} and summing the terms with their
factorial denominators gives
$O((v_q^qD^{1-q}+v_q^q)\log(e/w))=O(D\log(e/w))$, because $v_q\leq w\leq CD$. Entries above $v_*$ contribute $O(v_q^q)$.

\paragraph{Inward crossings near the boundary.}
Suppose $CD<\max_jz_j(U)\leq2$ and the original statistic is inside
the rectangle. Select a constraint attaining the maximum score of $U$. Impose \eqref{eq:cubic-QH3} on $\mathcal C_2(U)$ and, at the $O(\log(e/w))$ dyadic shifts $b$ from $c'D$ to a fixed
constant, the bounds $N_{C'b}(U;b)\leq m_{C'b}$ on $\mathcal C_b(U)$. Their total failure probability is $O(D\log(e/w))$.
If the largest deleted entry is below $2v_*$, the maximum score of
$U$ is at most $s=Cky$. For small $s$, let $b$ be the next dyadic value above $s$. The
remaining sample satisfies $\mathcal C_{b+h_k}(V)$, and the selected
score is greater than $-h_k$. Any candidate with score greater than
$-h_k$ after deletion had full-sample score in $[-2h_k,b]$. The number
of such candidates is therefore at most $m_{C'b}$. Equation
\eqref{eq:cubic-QH1} bounds the corresponding expected count by $O(s)$.
For larger $s$, use $\mathcal C_{\min(2,s)+h_k}(V)$ and the count bound
$D^{-1}e^{\lambda(2+2h_k)}$ from \eqref{eq:cubic-QH3}; this gives
$O(e^{Cs})$. If a deleted entry exceeds $v_*$, use $\mathcal C_{2+h_k}(V)$ and the
same count bound. The expectation is $O(e^{C(2+h_k)})$. The selected
score after deletion is greater than $-h_k$ because it attained a
positive maximum before deletion. Multiplying by
\eqref{eq:cubic-QH4} and summing again gives $O(D\log(e/w))$.

\paragraph{Inward crossings far outside.}
If $\max_jz_j(U)>2$, impose \eqref{eq:cubic-QH2} at $t=1$. Since $B_1^*\geq1$ and $A_1/(1+B_1^*)\leq D^{-1}$, the normalized exponential weights
$\omega_j=e^{\lambda z_j}/(1+\sum_le^{\lambda z_l})$ on faces with $z_j(U)>1$ have sum at least $cD$. On an inward crossing, every one of these scores must decrease by
more than one. The crossing indicator is therefore bounded by
$CD^{-1}$ times the sum of the weights of these crossings. For fixed $I$, $\omega_j(U)\leq e^{2\lambda h_k}\omega_j(V)$. Keep
this weight and the restrictions on the deleted rows, and discard the
other restrictions. The weights of $V$ sum to at most one. When the
largest deleted entry is below $2v_*$, we have $Cky\geq1$, so at most
$C(1+\log(ek))$ intervals remain and $y^{-q}\leq Ck^q$. Entries above
$v_*$ are handled together. Summing the factorially decaying probability
bounds gives $O(D^{-1}v_q^q)=O(D)$, since $q\geq3$ and $v_q\leq w$.
Thus, on the high-threshold constraints satisfying
\eqref{eq:cubic-high-counting}, the original and truncated statistics
disagree with probability $O(D\log(e/w))$.

\subsubsection{Combining low and high thresholds}
Split each sign group at $b_j=A_0\sqrt{\log(en)}$. If there are more than $n^{D_2}$ low faces, retain $\lceil n^{D_2}\rceil$ of them. For $D_2$ sufficiently large in terms of $a,A_0$, the Slepian
bound \eqref{eq:fixed-screen}, at threshold $A_0\sqrt{\log(en)}$,
makes their Gaussian probability $O(n^{-3})$. If $a=1$, apply the
bound at a fixed smaller positive gap. Their sample mass is $O(n^{-3}+n^{-(q-2)/(2q)}\log^2(en)+E_{q,a}(n))$ by Theorem~\ref{thm:fixed-sample-rate}. Since the full rectangle is contained in this subrectangle, the same
bound applies to its probability.
Otherwise the low subrectangle has polynomial dimension. Sort the high thresholds in each group. If $k\leq e^{D_1g(b_k)}$ fails, retain the shortest initial segment
for which it fails. Minimality gives the counting hypothesis \eqref{eq:cubic-high-counting} with factor two on that prefix. Since $\log k>D_1g(b_k)\geq D_1b_k^2/2$ and $b_k\geq A_0\sqrt{\log(en)}$, Slepian's inequality bounds its Gaussian mass by $O(n^{-3})$ when $A_0,D_1$ are sufficiently large. The preceding bounds for high-threshold crossings, bounded rows, and
covariance replacement bound the probability of this initial-segment
event under the original sample by $O(D\log(e/w)+n^{-3})$. The full rectangle is contained in this prefix event, so it has the
same probability bound. Upward score shifts preserve the counting condition.

It remains to combine the high thresholds for which the counting
condition holds with the polynomially many low thresholds. On the low group, the proof of \eqref{eq:fixed-sample-upper} controls the disagreement of the original statistic and the ordinary centered truncated sum $S$ by $C(n^{-(q-2)/(2q)}\log^2(en)+E_{q,a}(n))$. On the same observations, Bernstein's inequality and the truncated moment bound give, except on an event of probability $O(n^{-3})$,
\begin{align*}
 \max_j|S_j|&\leq C\sqrt{\log(en)},\\
 \max_j\left|n^{-1}\sum_i(Y_{ij}-\mu_j)^2-1\right|
 &\leq C\left\{\sqrt{\frac{(1+\tau^{(4-q)_+})\log(en)}{n}}
                 +n^{-(q-2)/q}\log(en)\right\}.
\end{align*}
Since $n^{-1/2}\sqrt{1+\tau^{(4-q)_+}}\leq Cn^{-(q-2)/(2q)}$, the truncated self-normalized vector satisfies $\|U-S\|_\infty\leq Cn^{-(q-2)/(2q)}\log(en)$ there. Expand and contract the rectangle formed by the low thresholds by
this deterministic amount. The comparison for the ordinary sum and
Gaussian anti-concentration then bound the probability that $U$ and
$S$ disagree on this subrectangle by
$O(n^{-(q-2)/(2q)}\log^2(en))$. One-sided limits treat atoms and empty inner rectangles.
A union bound over the two groups now bounds the probability that the
original and truncated self-normalized statistics disagree on the full
rectangle by
$C(D\log(e/w)+n^{-(q-2)/(2q)}\log^2(en)+E_{q,a}(n))$, without any independence assumption between groups. The comparison from $U$ to $G_{R'}$ in
\eqref{eq:cubic-bounded-bound}, followed by that from $G_{R'}$ to $G_R$
in \eqref{eq:cubic-covariance-comparison}, adds $O(D)$ to the error. Finally $n^{-(q-2)/(2q)}\log^2(en)\leq CF(w)$ and $D\log(e/w)=O(F(w))$. Hence
$\Delta_{n,q,K}(d,a)\leq1\wedge C[F(w)+E_{q,a}(n)]$.
Combining with \eqref{eq:fixed-sample-upper} proves \eqref{eq:fixed-complete-envelope}. Since $E_{q,a}(n)\to0$ for fixed $a>0$, sufficiency follows when $\log d=o(n^{1/3})$; necessity follows from the bounded example with nonzero third moment
in Lemma~\ref{lem:fixed-skew-lower}.

\subsection{The leading rate for a vanishing correlation gap}
\label{sec:shrinking-rate}
Fix $q\geq3$ and $K>1$, and suppose
\[
 a=a_n\longrightarrow0,\qquad A_n=a_n\log n\longrightarrow\infty.
\]

\begin{theorem}[The leading logarithmic rate when the correlation gap vanishes]
\label{thm:shrinking-coefficient-new}
Fix $0<\beta<1/3$ and $C_0>0$. For every $0<\zeta<\frac{q-2}{2}$, for all sufficiently large $n$,
\begin{equation}\label{eq:shrinking-upper-new}
 \sup_{2\leq d\leq\exp(C_0n^\beta)}\Delta_{n,q,K}(d,a_n)
 \leq C_{q,K,\beta,C_0,\zeta}\exp\{-(\frac{q-2}{2}-\zeta)A_n\}.
\end{equation}
For every fixed $c>0$, uniformly over $c n^{(q-2)/2}\leq d\leq\exp(C_0n^\beta)$,
\begin{equation}\label{eq:shrinking-equivalence-new}
 \frac{\log\Delta_{n,q,K}(d,a_n)}{A_n}\longrightarrow-\frac{q-2}{2}.
\end{equation}
\end{theorem}

\subsubsection{Selecting low and high thresholds}

As before, split the finite signed constraints of the rectangle into
two sign groups. The correlation condition holds separately within
each group. Assume first that the rectangle has Gaussian probability at least
\[
 \delta_0=\exp\{-(\frac{q-2}{2}-\zeta)A_n\}.
\]
In each group select the first $m_0=\lfloor n^{(q-2-\zeta)/2}\rfloor$ faces in increasing order of threshold, or all faces if fewer are present. We call these the low-threshold constraints and the remaining ones
the high-threshold constraints. There are at most
$2n^{(q-2-\zeta)/2}$ low-threshold constraints in total.

Let $B$ be the lowest threshold in a nonempty high-threshold group.
Apply Slepian's inequality \citep{Slepian1962} to the $m_0$ selected
low-threshold constraints in the same sign group, comparing with the
equicorrelated Gaussian vector. For every $u$, this gives
\[
 \delta_0\leq
 \Phi\!\left(\frac{B-\sqrt a\,u}{\sqrt{1-a}}\right)+\Phi(u)^{m_0}.
\]
Choose $u=\{2\log m_0-6\log\log m_0\}^{1/2}$. By Mills' inequalities,
\[
 \Phi(u)^{m_0}\leq\exp\{-c(\log m_0)^2\}=o(\delta_0).
\]
Moreover,
\[
 \frac{\sqrt a\,u}{\sqrt{2A_n}}\longrightarrow\sqrt{\frac{q-2}{2}-\zeta/2},
 \qquad
 \frac{\Phi^{-1}(\delta_0/2)}{\sqrt{2A_n}}
 \longrightarrow-\sqrt{\frac{q-2}{2}-\zeta}.
\]
It follows that $B\geq c_{q,\zeta}\sqrt{A_n}$ for large $n$.

For $h\geq B$, let $k$ be the number of high thresholds at most $h$ in one group, and let $M_k$ be the maximum over their Gaussian coordinates. When $k\geq2$, Sudakov--Fernique comparison gives $\mathbb E M_k\geq c\sqrt{a\log k}$. Since $\mathbb P(M_k\leq h)\geq\delta_0$, Gaussian concentration gives
\[
 \mathbb E M_k\leq h+\sqrt{2\log(1/\delta_0)}\leq C_{q,\zeta}h.
\]
For the rest of this subsection, use the smooth approximation to
$x_+^2$ given by
\[
 g(x)=\frac14\bigl(x+\sqrt{1+x^2}\bigr)^2,
\]
in place of $-\log\bar\Phi(x)$. It satisfies $g(h)\geq h^2$ for
$h\geq0$. Thus each high set $I$ satisfies, for a fixed $D_0=D_0(q,\zeta)$,
\begin{equation}\label{eq:high-count-new}
 \#\{j\in I:g(b_j)\leq t\}
 \leq2\exp\{D_0(1+t)/a\},\qquad t\geq0.
\end{equation}
The same inequality remains valid when every selected threshold score is increased by a common nonnegative amount.

\subsubsection{A high-threshold Gaussian moment estimate}

\begin{lemma}[Fractional moments of the high-threshold count]
\label{lem:high-fractional-new}
Let $C=\{x:\langle v_j,x\rangle\leq b_j\}$ be a Gaussian
polyhedron in its finite-dimensional factor space, where the $v_j$ are
unit vectors, and let $\gamma$ be standard Gaussian measure on that
space. Select a subset $I$ of the constraints, while retaining all
constraints in the event $C$. Suppose the selected vectors have pairwise inner products at most $1-a$ and their thresholds satisfy \eqref{eq:high-count-new}. Fix $H>0$ and take $a$ sufficiently small that $\alpha=Ha\leq1$. There is $C_H=C(H,D_0)$, independent of dimension, $a$, and $\gamma(C)$, such that
\begin{equation}\label{eq:high-count-moment-new}
 \mathbb E\bigl[\mathbf1_C N_y^\alpha\bigr]\leq C_H e^{C_Hy},
 \quad N_y=\#\{j\in I:g(b_j)-g(\langle v_j,Z\rangle)\leq y\},\quad y\geq0.
\end{equation}
For some $\lambda=C_H/a$, both $u=\lambda$ and $u=\lambda/2$ satisfy
\begin{equation}\label{eq:high-weight-moment-new}
 \mathbb E\left[\mathbf1_C
 \left(\sum_{j\in I}e^{-u\{g(b_j)-g(\langle v_j,Z\rangle)\}}\right)^\alpha\right]
 \leq C_H.
\end{equation}
The same estimates hold after increasing the selected threshold
scores or adding further Gaussian constraints to $C$.
\end{lemma}
\begin{proof}
Write $\pi=\gamma(C)$. The case $\pi=0$ is immediate. Define $h_s(z)=\mathbb P(z+\sqrt sZ\in C)$ and
\[
 V_C(x)=\int_0^1
 \mathbb E\left|\nabla\log h_s\bigl((1-s)x+\sqrt{s(1-s)}Z\bigr)\right|^2\,ds.
\]
We first prove
\begin{equation}\label{eq:energy-exp-new}
 \int_C e^{V_C(x)/8}\,d\gamma(x)\leq\pi^{1/2}.
\end{equation}
Let $P$ be Wiener measure and $Q=P(\,\cdot\mid B_1\in C)$. For
$t<1$, the density of $Q$ relative to $P$ on the Brownian filtration
is $L_t=h_{1-t}(B_t)/\pi$. With
$u_t=\nabla\log h_{1-t}(B_t)$, Girsanov's theorem and It\^o's formula
give, under $Q$,
\[
 \log\frac{h_{1-t}(B_t)}\pi=M_t+\frac12A_t,
 \qquad M_t=\int_0^t u_r\,dW_r,\quad A_t=\int_0^t|u_r|^2\,dr.
\]
Here $W$ is Brownian motion under $Q$. This is the F\"ollmer drift
representation \citep{LEHEC2013}. To justify the terminal identities,
set
\[
 \sigma_m=(1-1/m)\wedge\inf\{t:A_t\geq m\},\qquad m\geq2.
\]
The stopped stochastic integral has mean zero under $Q$, so
\[
 \mathbb E_Q A_{\sigma_m}
 =2\mathbb E_Q\log L_{\sigma_m}
 =2\mathbb E_P[L_{\sigma_m}\log L_{\sigma_m}].
\]
The stopping times increase to one, since $A_t$ is finite for every
$t<1$. Also $0\leq L_{\sigma_m}\leq\pi^{-1}$ and
$L_{\sigma_m}\to\mathbf1_C(B_1)/\pi$ under $P$. The function
$x\log x$, with $0\log0=0$, is bounded on $[0,\pi^{-1}]$.
Dominated convergence for the entropy terms and monotone convergence
for $A_{\sigma_m}$ therefore give
$\mathbb E_QA_1=2\log(1/\pi)<\infty$. Consequently $M_t$ converges
in $L^2(Q)$, and the terminal identity is
$M_1+A_1/2=\log(1/\pi)$, $Q$-almost surely.

Apply Fatou's lemma to the nonnegative exponential supermartingale
$\exp(-M_t/2-A_t/8)$. The terminal identity gives
$\mathbb E_Qe^{A_1/8}\leq\pi^{-1/2}$. Conditional on $B_1=x\in C$,
the path has the ordinary Brownian bridge law. Hence
$\mathbb E_Q(A_1\mid B_1=x)=V_C(x)$ for $\gamma$-almost every
$x\in C$. Conditional Jensen's inequality, followed by integration
against $\gamma$ on $C$, proves \eqref{eq:energy-exp-new}.

Fix $x\in C$ with $V_C(x)<\infty$, and put $w=1+y$,
$B_y=4\sqrt w$. The pointwise count bound below is automatic when
$V_C(x)=\infty$. The selected thresholds at most $B_y$ contribute at most $2\exp(CD_0w/a)$ indices by \eqref{eq:high-count-new}. For $l\geq0$, let $k_l$ be the number of active indices with
\[
 2^lB_y\leq b_j<2^{l+1}B_y,
 \qquad g(b_j)-g(x_j)\leq y,\quad x_j=\langle v_j,x\rangle.
\]
The score-gap condition implies $b_j/2\leq x_j\leq b_j$ and
$0\leq b_j-x_j\leq2y/b_j$. Fix $U\geq1$ and use
\[
 J_l(U)=[e^{-U}4^{-l}/16,\,4^{-l}/8].
\]
For every selected index in this band and every $s\in J_l(U)$,
\[
 (b_j-x_j)/\sqrt s+\sqrt s\,x_j\leq C_U\sqrt w.
\]
Gaussian differentiation gives
\[
 s\nabla\log h_s(z)=\mathbb E(W\mid W\in C)-z,
 \qquad W\sim N(z,sI).
\]
The conditional mean belongs to the convex set $C$. Hence its
distance from $z$ is at least the positive amount by which $z$ violates
any selected defining inequality. Applying this to the preceding
gradient identity shows that the contribution of $J_l(U)$ to $V_C(x)$
is at least
\[
 \int_{J_l(U)}\frac{ds}{s}
 \mathbb E\bigl(\sqrt{1-s}M_l-C_U\sqrt w\bigr)_+^2
 \geq (U+\log2)\bigl(c\mathbb EM_l-C_U\sqrt w\bigr)_+^2,
\]
where $M_l$ is the maximum of the selected Gaussian projections in that band. Sudakov--Fernique gives $\mathbb EM_l\geq c\sqrt{a\log k_l}$ for $k_l\geq2$. Consequently,
\begin{equation}\label{eq:band-energy-new}
 \log(k_l\vee1)\leq\frac{CV_C(x)}{Ua}+\frac{C_Uw}{a}.
\end{equation}
The coefficient $C$ in the first term is absolute; this follows from $(r-t)_+^2\geq r^2/2-t^2$ after rescaling constants.

To count occupied bands, use instead the disjoint intervals $[4^{-l}/16,4^{-l}/8]$. In an occupied band, a single standard normal projection lying in $[K\sqrt w,(K+1)\sqrt w]$, for a sufficiently large numerical $K$, makes the same integrand at least $ce^{-Cw}$ in expectation. Summing over disjoint intervals gives
\[
 \#\{l:k_l>0\}\leq Ce^{Cw}V_C(x).
\]
Combining this with \eqref{eq:band-energy-new} and the low-threshold count yields
\[
 \alpha\log(1+N_y)\leq C_{U,H,D_0}w
 +\alpha\log(1+V_C(x))+(CH/U)V_C(x).
\]
Choose $U$ so that $CH/U\leq1/16$. Since $\alpha\leq1$, $\alpha\log(1+v)\leq v/16+C$ for $v\geq0$. Exponentiating and applying \eqref{eq:energy-exp-new} proves \eqref{eq:high-count-moment-new}.

On $C$, every selected score gap is nonnegative, and
\[
 \sum_{j\in I}e^{-u\{g(b_j)-g(x_j)\}}
 \leq\sum_{k\geq0}e^{-uk}N_{k+1}.
\]
Subadditivity of the $\alpha$th power and \eqref{eq:high-count-moment-new} prove \eqref{eq:high-weight-moment-new} once $\alpha\lambda/2>C_H+1$. Taking $\lambda=C'_H/a$ ensures this for both values of $u$. Using
score intervals of length one keeps the geometric-series bound
independent of $1/\alpha$. The other constraints of $C$ remain in each event, and increasing the
threshold scores preserves \eqref{eq:high-count-new}. This proves the final assertion.
\end{proof}

\subsubsection{Joint comparison of truncated and untruncated coordinates}

\begin{lemma}[Joint comparison for low and high thresholds]
\label{lem:mixed-comparison-new}
Suppose $\log(2nd)\leq C_0n^\beta$, $0<\beta<1/3$, and $a\to0$, $a\log n\to\infty$. Fix $e>0$. Assume there are at most $C n^{(q-2)/2-e}$ low signed faces. For the high-threshold constraints, set $z_j=g(G_j)-g(b_j)$ and
$Y=\max_jz_j$, where $G$ is the corresponding signed-coordinate
version of $G_R$. Suppose that for each fixed $H>0$, eventually $\alpha=Ha\leq1$ and, for some $\lambda=C_H/a$, at $u=\lambda,\lambda/2$ and every $r\geq0$,
\begin{equation}\label{eq:mixed-hypothesis-new}
 \mathbb E\left[\mathbf1_{\{Y\leq r\}}
 \left(1+\sum_{j\,\mathrm{high}}e^{u(z_j-r)}\right)^\alpha\right]\leq C_H.
\end{equation}
There are $\varepsilon,c>0$ depending only on $q,K,\beta,C_0,e$,
and independent of $H$, such that the discrepancy for this rectangle
between $T_n$ and $G_R$ is at most
\begin{equation}\label{eq:mixed-result-new}
 C_H\{n^{-c}+e^{-H\varepsilon a\log n}\}.
\end{equation}
The bound holds for all sufficiently large $n$, with the starting
index allowed to depend on $H$ and the given sequence $(a_n)$.
An empty low or high set is simply omitted.
\end{lemma}
\begin{proof}
\emph{Truncation and the Gaussian coupling.}
Choose $0<r_0<\min(e/(4q),1/16)$ and put $\tau=n^{1/2-r_0}$. For a low signed coordinate define $Y_l=X_l\mathbf1_{\{|X_l|\leq\tau\}}$, $\mu_l=\mathbb EY_l$, and $W_l=n^{-1/2}\sum_i(Y_{il}-\mu_l)$. The probability that truncation changes any observation in a low
coordinate is at most $Cn^{-e+qr_0}$. Furthermore,
\[
 \sqrt n\max_l|\mu_l|\leq Cn^{-\frac{q-2}{2}+(q-1)r_0},\qquad
 \mathbb EY_l^4\leq C(1+\tau^{(4-q)_+}).
\]
Bernstein's inequality and a union bound over the polynomially many low columns give, except on an event of probability $O(n^{-3})$,
\begin{align*}
 \max_l|W_l|&\leq C\sqrt{\log n},\\
 \max_l\left|n^{-1}\sum_iY_{il}^2-1\right|
 &\leq C\left\{\tau^{2-q}
 +\sqrt{(1+\tau^{(4-q)_+})\log n/n}
 +n^{-2r_0}\log n\right\}.
\end{align*}
For some $c_0>0$, every original self-normalized low coordinate
differs from its corresponding $W_l$ by at most $n^{-c_0}$, except on
an event of probability $O(n^{-c_0})$.

Use as the comparison row the joint vector $(Y_l-\mu_l:l\ \mathrm{low};X_j:j\ \mathrm{high})$, and let its Gaussian replacement have exactly the full joint covariance. Couple its ordinary Gaussian sum $(W^G,G_{\mathrm{high}})$ with the original Gaussian $(G_{\mathrm{low}},G_{\mathrm{high}})$ using the covariance of $(Y_l-\mu_l,X_l,X_j)$. The high coordinates in this coupling coincide; each low difference has standard deviation at most $C\tau^{1-q/2}$. After decreasing $c_0$, the maximum absolute difference is at most
$n^{-c_0}$ except on an event of probability $O(n^{-3})$. Gaussian anti-concentration bounds a boundary strip of width $b$
for either Gaussian low-coordinate vector by $Cb\sqrt{\log n}$. This joint coupling includes an upper and lower face of the same original coordinate when they fall in different sets.

\emph{Smooth brackets.}
Let $h$ be an increasing $C^3$ function equal to zero on $(-\infty,0]$ and one on $[1,\infty)$, and set $\chi(u)=1-h(u-1)$. Choose a small $\varepsilon>0$ below and put $\delta=n^{-\varepsilon}$, $M=n^\varepsilon$. For $N_\ell$ low faces set
\[
 D(W)=\theta^{-1}\log\sum_{l\,\mathrm{low}}e^{\theta(W_l-b_l)},
 \qquad \theta=\delta^{-1}\log(eN_\ell).
\]
The difference between $D(W)$ and $\max_l(W_l-b_l)$ lies in
$[0,\delta]$. The tests
\[
 P^-=1-h((D+2\delta)/\delta),\qquad
 P^+=1-h((D-2\delta)/\delta)
\]
bound the indicator of the original low-threshold event from below and
above, respectively, whenever
every original low-coordinate ratio differs from $W_l$ by at most
$n^{-c_0}\leq\delta$. The preceding estimates show that these
approximation bounds fail with probability $O(n^{-c_0})$.
For high ratios $x$, put
\[
 z_j(x)=g(x_j)-g(b_j),\quad
 A(x)=1+\sum_j e^{\lambda z_j(x)},\quad
 B_t(x)=\sum_j e^{\lambda z_j(x)}h(z_j(x)/\delta+t),\quad t=0,1,
\]
\[
 H^-(x)=\chi(A(x)/M)\{1-h(B_1(x))\},\qquad
 H^+(x)=1-h\!\left(\frac{2MB_0(x)}{M+A(x)}\right).
\]
These functions bound the indicator of the high-threshold event
from below and above: $B_0=0$ when all scores are nonpositive, whereas
$B_1>1$ when a score is positive. Thus $P^-H^-$ and $P^+H^+$ bracket
the full rectangle event on the same event of simultaneous
approximation for the low coordinates.

\emph{Deleted-row estimates.}
For each fixed high coordinate $j$, the remaining observations are
independent and centered, with unit variance and a uniformly bounded
third absolute moment. They may come from either the original
rows or their matching Gaussian rows. Their squared norm $Q_j$ and
self-normalized ratio $t_j$ satisfy
\[
 \mathbb P(Q_j<N/2)\leq Ce^{-cN},\qquad
 \mathbb P(|t_j|>x)\leq Ce^{-cx^2}+Ce^{-cN},\quad N=n-1.
\]
The first estimate follows by truncating the squares at a fixed
constant and applying Hoeffding's inequality. For the second, let
$S_j$ be the sum of the remaining observations.
The bound $e^{u-u^2/2}\leq1+u+Cu^2$ gives $\mathbb Ee^{hS_j-h^2Q_j/2}\leq e^{CNh^2}$. On each dyadic shell $b\leq Q_j<2b$, $b\geq N/2$, choose $h=c_1x/\sqrt b$ and apply Chernoff to either sign. Sum $O(1+x^2)$ shells; Markov's inequality for $Q_j$ bounds the remaining shells by $Ce^{-cx^2}$. A union bound therefore gives an event, determined by the remaining
sample and of probability $1-O(n^{-20})$, on which
\[
 Q_j\geq(n-1)/2,\qquad \max_j|t_j|\leq C\sqrt{\log(2nd)}.
\]

\emph{Inserting one observation in the high coordinates.}
On this event, for an incoming high row $y$, set
\[
 v_j=(1+|t_j|)y_j/\sqrt{Q_j},\quad
 w_j=e^{\lambda z_j(t_j)}/A(t),\quad
 Q_* =\sum_jw_j|v_j|^3.
\]
Conditionally on the remaining sample, these weights are
deterministic, independent of the incoming row, and satisfy
$\sum_jw_j\leq1$. The exact inserted ratio and its exponential score increment are
\[
 x_j(v)=\frac{t_j+v/(1+|t_j|)}{\sqrt{1+\{v/(1+|t_j|)\}^2}},\qquad
 r_j(v)=e^{\lambda\{g(x_j(v))-g(t_j)\}}.
\]
Globally, $g(x_j(v))-g(t_j)\leq3/2$, because $0\leq g(x)-x_+^2\leq1/2$ and $(x_j(v))_+^2\leq(t_j)_+^2+1$. On $|v|\leq1$, the first three derivatives of the score increment are bounded absolutely, uniformly in $t_j$. This follows by direct differentiation of the exact ratio, using $g'(x)\leq C(1+|x|)$ and bounded $g''$, $g'''$.
These bounds are uniform in $j$, the retained values $t_j$, and the
thresholds. Fix $H$ and the sequence $(a_n)$, and choose a constant
$C_h\geq1$, depending only on the fixed cutoff $h$, large enough that
\[
 L_n=C_h e^{3\lambda/2}(1+\lambda)^3,\qquad
 r_n=\frac{\log L_n}{\log n}\longrightarrow0.
\]
The convergence follows from $\lambda=C_H/a=o(\log n)$. Each of the
two functions $r_j(v)$ and $r_j(v)h(z_j(x_j(v))/\delta+t)$ is bounded
globally by $L_n$, and its derivative of order $1\leq k\leq3$ on
$[-1,1]$ is bounded by $L_n\delta^{-k}$. The same $L_n$ works for
every law in the row class, dimension, threshold vector, retained
sample, replacement step, and coordinate.

Normalize the updated sums $A'$ and $B_t'$ by the retained value
$A=A(t)$:
\[
 A'/A=1/A+\sum_jw_jr_j(v_j),\qquad
 B_t'/A=\sum_jw_jr_j(v_j)h(z_j(x_j(v_j))/\delta+t).
\]
Their linear terms, quadratic terms, and remainders after the
quadratic terms are bounded respectively by
\[
 L_n\delta^{-1}Q_*^{1/3},\quad
 L_n^2\delta^{-2}Q_*^{2/3},\quad
 L_n^3\delta^{-3}Q_*.
\]
Indeed, Taylor's theorem gives the scalar remainder bounds for
$|v_j|\leq1$. For $|v_j|>1$, bound the function globally and its
coefficients at zero; each term is then bounded by the corresponding
power of $|v_j|$. Since $L_n\geq1$ and $\delta\leq1$, increasing
$C_h$ makes these bounds at order $k$ at most
$L_n^k\delta^{-k}|v_j|^k$. Weighted H\"older's inequality gives
$\sum_jw_j|v_j|^k\leq Q_*^{k/3}$ for $k=1,2,3$, proving the
displayed estimates. The same argument bounds the full first
increment by $L_n\delta^{-1}Q_*^{1/3}$.

At zero, the normalized coordinates belong to $\{1\}\times[0,1]$.
On $1/2\leq u\leq3/2$, $|v|\leq2$, each outer bracket, viewed as a
function of $(u,v)$, has derivatives of order $i\leq3$ bounded by
$C_iM^i$. For $H^+$, use $A/(M+Au)\leq2$. For $H^-$, a nonzero
value or derivative requires $Au\leq2M$, hence $A\leq4M$.

There are consequently homogeneous polynomials $H_1,H_2$ of
degrees one and two in the incoming row, with $H_0=H(t)$, such that,
for $i=1,2,3$,
\begin{equation}\label{eq:mixed-Taylor-new}
 \left|H(x(v))-\sum_{j=0}^{i-1}H_j\right|
 \leq C_iL_n^i n^{2i\varepsilon}Q_*^{i/3}.
\end{equation}
To verify this for every incoming row, put
$u_*=L_n\delta^{-1}Q_*^{1/3}$. If $u_*$ is below a sufficiently
small absolute constant, the full inner increment stays in the
indicated box. Taylor expansion of the outer bracket, together with
the inner remainder bounds, gives $C_iM^iu_*^i$. Otherwise the
bracket is bounded by one and its polynomial terms satisfy
$|H_j|\leq C_jM^ju_*^j$, $j=1,2$. Because $M\geq1$ and $u_*$ is
then bounded below, their sum is again at most $C_iM^iu_*^i$ at the
required order. This proves \eqref{eq:mixed-Taylor-new}; in particular,
its factor $L_n^i=n^{ir_n}$ is uniform over all the objects listed
above. Conditionally on the remaining sample,
\begin{equation}\label{eq:mixed-Q-moment-new}
 \mathbb EQ_*\leq C(\log(2nd)/n)^{3/2}
\end{equation}
for either incoming law, using only its marginal third absolute moments.

\emph{The low test and joint moment cancellation.}
For either low-threshold test, write $P_0,P_1,P_2$ for its Taylor
terms through degree two in the incoming row. Its softmax weights are
$p_l=e^{\theta(W_l-b_l)}/\sum_{k\,\mathrm{low}} e^{\theta(W_k-b_k)}$. Along the
insertion path, they are bounded by a constant times their retained
values whenever $\theta\max_l|y_l|/\sqrt n\leq1$. For original truncated rows this holds if $\varepsilon<r_0/2$; for Gaussian rows its failure probability is exponentially small by a union bound over polynomially many low faces. Differentiating $p_l$ via
$\partial_kp_l=\theta p_l(\mathbf1_{k=l}-p_k)$ gives, at each order
$i\leq3$, a nonnegative tensor determined by the remaining sample
that bounds the derivatives along this path. Its total mass is at
most $C_i n^{i\varepsilon}\log^i(en)$, because
$\log(eN_\ell)\leq C\log(en)$. Conditional H\"older's inequality,
including repeated coordinates, therefore gives
\[
 \|P_i\|_{3/i}\leq C_i n^{i\varepsilon}\log^i(en)n^{-i/2}
 \quad(i=1,2),\qquad
 \mathbb E|P-P_0-P_1-P_2|
 \leq Cn^{3\varepsilon}\log^3(en)n^{-3/2}.
\]
These expectations and norms are conditional on the remaining sample;
the constants are uniform over its permitted values and all low
thresholds.
Cauchy--Schwarz and moment bounds for the Gaussian maximum show that
the exceptional Gaussian event contributes an exponentially small error,
including to the quadratic Taylor term.

The polynomial of total degree at most two for the product is
\[
 P_0H_0+P_1H_0+P_0H_1+P_2H_0+P_0H_2+P_1H_1.
\]
Its conditional expectations are identical for the original row and the jointly matching Gaussian row, including the low--high cross terms. The remainder is exactly
\[
 (P-P_0-P_1-P_2)H+P_0(H-H_0-H_1-H_2)
 +P_1(H-H_0-H_1)+P_2(H-H_0).
\]
Put $\widetilde L_n=L_n+C\log(en)$, so that
$\log\widetilde L_n/\log n\to0$ for fixed $H$. The preceding
estimates, \eqref{eq:mixed-Taylor-new},
\eqref{eq:mixed-Q-moment-new}, and conditional H\"older's inequality
bound the expected absolute remainder for one row by
$C\widetilde L_n^3 n^{6\varepsilon}(\log(2nd)/n)^{3/2}+O(n^{-10})$.
On the exceptional event for the remaining sample, use boundedness of
the test functions. Summing over all row replacements gives total error
\[
 C\widetilde L_n^3 n^{-(1-3\beta)/2+6\varepsilon}+O(n^{-9}).
\]
Choose once and for all
\[
 0<\varepsilon<\min\{r_0/2,c_0/2,(1-3\beta)/24\}.
\]
This choice is independent of $H$. Since $\widetilde L_n=n^{o(1)}$
for each fixed $H$, the replacement error is bounded by a fixed
negative power of $n$ whose exponent is independent of $H$. Only the
starting index may depend on $H$ and $(a_n)$.

\emph{Comparison with the Gaussian rectangle probability.}
Let $G$ be the normalized Gaussian sum over the high coordinates and
$T^G$ the vector of ratios formed from the same Gaussian rows.
Gaussian and chi-square tail bounds give, except on an event of
probability $O(n^{-3})$,
\[
 \max_{j\,\mathrm{high}}|g(T_j^G)-g(G_j)|
 \leq C\log^{3/2}(2nd)/\sqrt n\leq\delta.
\]
For $Y=\max_{j\,\mathrm{high}}\{g(G_j)-g(b_j)\}$, every interval
of width $v$ has probability at most $Cv$. Indeed,
$\phi(t)/(g'(t)\overline\Phi(t))$ is bounded on $\mathbb R$. Apply the
argument of Lemma~\ref{lem:score-strip} to the shifted threshold scores,
conditioning on the maximizing coordinate; regularization covers
singular covariance matrices.

If $Y<-2\delta$ and $H^-(T^G)\ne1$, then $A(G)>Me^{-\lambda\delta}$. Hypothesis \eqref{eq:mixed-hypothesis-new} and Markov's inequality give probability at most $C_HM^{-\alpha}$. To control $H^+$, define
\[
 B_*(z)=\sum_j e^{\lambda z_j}h(z_j),\qquad
 R(z)=\frac{1+\sum_j e^{\lambda z_j}}{1+B_*(z)}.
\]
Uniformly in upward shifts $r\geq0$,
\begin{equation}\label{eq:R-moment-new}
 \mathbb ER(z-r)^\alpha\leq C_H.
\end{equation}
On $\max_j(z_j-r)\leq1$, use \eqref{eq:mixed-hypothesis-new} after one further upward shift. On the slice $k\leq\max_j(z_j-r)<k+1$, $k\geq1$, the denominator is at least $e^{\lambda k}$, whereas
\[
 A(z-r)-1-B_*(z-r)\leq e^{\lambda/2}\sum_j e^{\lambda(z_j-r)/2}.
\]
Shift the latter scores by $k+1$ and apply the hypothesis at $\lambda/2$. Subadditivity produces a geometric series with ratio controlled by $\alpha\lambda$, a fixed positive constant; this proves \eqref{eq:R-moment-new}.

If $Y>2\delta$ and $H^+(T^G)>0$, then $B_0(T^G)\geq1$ and $A(T^G)>MB_0(T^G)$. On the good Gaussian coupling event, the ordinary scores shifted by $2\delta$ satisfy
\[
 A(z-2\delta)\geq e^{-3\lambda\delta}A(T^G),\qquad
 B_*(z-2\delta)\leq B_0(T^G).
\]
Consequently $R(z-2\delta)\geq cM$, and \eqref{eq:R-moment-new} bounds this probability by $C_HM^{-\alpha}$. The remaining boundary strip has probability at most $C\delta$. For the low coordinates, the joint coupling above bounds the error
at the Gaussian endpoint by $O(\delta\sqrt{\log n}+n^{-3})$. Products of $[0,1]$-valued factors differ by at most the sum of their individual differences. Since $M^{-\alpha}=e^{-H\varepsilon a\log n}$, adding the truncation
error, the replacement error, and the two errors at the Gaussian
endpoint proves \eqref{eq:mixed-result-new}.
\end{proof}

\subsubsection{Completion of the upper rate}

For rectangles of Gaussian probability at least $\delta_0$, the
threshold selection satisfies the low-coordinate bound with
$e=\zeta/2$ and the high-coordinate count bound
\eqref{eq:high-count-new}. Lemma~\ref{lem:high-fractional-new}, applied to each sign group and combined by subadditivity, proves \eqref{eq:mixed-hypothesis-new}. Choose $\varepsilon$ as in Lemma~\ref{lem:mixed-comparison-new}, and only then choose a fixed $H$ with $H\varepsilon>\frac{q-2}{2}$. Since $A_n=o(\log n)$, both terms in \eqref{eq:mixed-result-new} are $o(\delta_0)$.

For a rectangle with Gaussian probability smaller than $\delta_0$, enlarge its coordinate intervals continuously until that probability equals $\delta_0$. Such an enlargement exists, and Gaussian marginal continuity gives continuity even for singular covariance. Inclusion controls the positive sample-minus-Gaussian error by the comparison error for the enlarged rectangle plus $\delta_0$; the opposite error is at most the original Gaussian mass. This proves \eqref{eq:shrinking-upper-new} uniformly over all rectangles.

\subsubsection{A lower-bound construction in polynomial dimension}
Fix $c>0$, set $d_0=\lfloor cn^{(q-2)/2}\rfloor$ and $\rho=1-a$, and put
\begin{gather*}
 t=\sqrt{2a\log(2d_0)},\qquad b=t+4,\\
 \eta=\eta_0b^{-q},\qquad
 \pi_0=\frac{\eta}{nd_0},\qquad v=nb^2\pi_0(1-\pi_0).
\end{gather*}
Take mutually independent symmetric signs $U_i$, standard normals $E_{ij}$ and Bernoulli$(\pi_0)$ marks, and define
\[
 X_{ij}=\frac{\sqrt\rho\,U_i+\sqrt a\,E_{ij}
                    +\sqrt n\,b(I_{ij}-\pi_0)}{\sqrt{1+v}}.
\]
This row is centered, with unit marginal variances and off-diagonal
covariances $\rho/(1+v)\leq1-a$. For fixed $q$, the $L^q$ norm of
$\sqrt\rho\,U_i+\sqrt a\,E_{ij}$ tends to one as $a\to0$. The rare
perturbation has $q$th moment at most
$C\eta_0n^{(q-2)/2}/d_0$. Thus, for every $K>1$, a sufficiently small
fixed $\eta_0=\eta_0(q,K,c)>0$ ensures the required moment bound for
all large $n$.

Write
\[
 W=n^{-1/2}\sum_iU_i,\qquad
 H_j=n^{-1/2}\sum_iE_{ij},\qquad
 N_j=\sum_iI_{ij},\qquad e_0=\eta b/d_0.
\]
At zero thresholds only the numerator signs matter. Thus the sample probability is
$\mathbb P(M_1\leq e_0)$, where
\[
 M_1=\max_{j\leq d_0}\{\sqrt\rho\,W+\sqrt a\,H_j+bN_j\}.
\]
Let $M_0$ be the same maximum without the marks. Exactly one mark occurs with probability at least $\eta/2$. Conditional on its column $j$, impose
\[
 -t-2\leq\sqrt\rho\,W\leq-t-1,\qquad
 0\leq H_j\leq1,\qquad
 \max_{k\ne j}H_k\leq\sqrt{2\log(2d_0)}.
\]
These three events are independent, and their joint probability is at least
\[
 c\exp\{-(t+2)^2/(2\rho)\}-Cn^{-1/2}
\]
by the scalar Berry--Esseen bound and a Gaussian union bound. The error is negligible because $t^2=O(a\log n)=o(\log n)$. On the stated event,
$M_0\leq-1$, whereas the marked coordinate entering $M_1$ is at least two.
Marks only increase the maximum, and $e_0\to0$, so
\begin{equation}\label{eq:one-mark-crossing-new}
 \mathbb P(M_0\leq e_0)-\mathbb P(M_1\leq e_0)
 \geq c\eta\exp\{-(t+2)^2/(2\rho)\}.
\end{equation}

Let $H_{\max}=\max_{j\leq d_0}H_j$. The probability of the zero-threshold lower orthant under the
comparison Gaussian is
\[
 \mathbb E\Phi(-\sqrt{a+v}\,H_{\max}/\sqrt\rho).
\]
Condition on $H_{\max}$ and apply Berry--Esseen to $W$. The bounded normal density and
$\sqrt{a+v}-\sqrt a\leq\sqrt v$ show that its difference from
$\mathbb P(M_0\leq e_0)$ is at most
\[
 C\{n^{-1/2}+e_0+\sqrt v\,\mathbb E|H_{\max}|\}
 \leq C\{n^{-1/2}+d_0^{-1/2}\sqrt{\log(2d_0)}\}.
\]
This is polynomially small, whereas \eqref{eq:one-mark-crossing-new} is $n^{-o(1)}$. Appending independent symmetric signs gives, for $d\geq d_0$,
\begin{equation}\label{eq:shrinking-lower-new}
 \begin{aligned}
 \Delta_{n,q,K}(d,a)
 &\geq C^{-1}(1+\sqrt{a\log n})^{-q}\\
 &\qquad{}\times\exp\left\{-\frac{(q-2)a\log n}{2(1-a)}
                  -C\sqrt{a\log n}-C\right\}\\
 &=\exp\{-[\tfrac{q-2}{2}+o(1)]a\log n\}.
 \end{aligned}
\end{equation}
Together with \eqref{eq:shrinking-upper-new}, followed by $\zeta\downarrow0$,
this proves \eqref{eq:shrinking-equivalence-new}. Because all thresholds are zero, the same lower-bound construction
applies to the ordinary normalized sum.

For the multiplicative correction, Proposition~\ref{prop:rare-mark-unified}
at $m_n\asymp n^{(q-2)/2}(\log n)^{q+1}$ gives
\begin{equation}\label{eq:shrinking-correction-new}
 \begin{aligned}
 \log\Delta_{n,q,K}(m_n,a_n)
 &\geq-\frac{(q-2)a_n\log n}{2(1-a_n)}
       +\frac{\sqrt{(q-2)a_n\log n}}{1-a_n}\\
 &\qquad{}-C\{1+a_n\log\log n+\log(2+a_n\log n)\}.
 \end{aligned}
\end{equation}
The last term is $o(\sqrt{a_n\log n})$. For example, taking
$a_n=(\log n)^{-1/2}$ yields
\[
 \log\Delta_{n,q,K}(m_n,a_n)
 \geq-\frac{q-2}{2}\sqrt{\log n}
        +\sqrt{q-2}\,(\log n)^{1/4}-O(\log\log n).
\]
Consequently, for every fixed $B_0<\infty$,
\[
 \frac{\Delta_{n,q,K}(m_n,a_n)}
 {\exp\{-(q-2)a_n\log n/2\}\log^{B_0}(enm_n)}
 \longrightarrow\infty.
\]
Thus the leading logarithmic rate cannot, uniformly over such gap
sequences, be converted to a multiplicative estimate using only a fixed
power of a logarithm.

\begin{proof}[Proof of Theorem~\ref{thm:main-shrinking}]
Theorem~\ref{thm:shrinking-coefficient-new} proves
\eqref{eq:main-leading-gap} and the sufficiency of
$a_n\log n\to\infty$ in \eqref{eq:main-exponential-criterion} along
subsequences on which $a_n\to0$. Along subsequences where the gaps stay bounded below, use
Theorem~\ref{thm:fixed-complete-envelope} at a fixed smaller positive gap.
If $a_n\log n$ fails to diverge, take a subsequence on which it is bounded
and apply Theorem~\ref{thm:main-polynomial}. This proves the converse.
For \eqref{eq:main-cubic-criterion}, take any fixed positive gap
below $\liminf_n a_n$ in Theorem~\ref{thm:fixed-complete-envelope}
to obtain sufficiency. Necessity follows from the independent skewed
coordinates of Lemma~\ref{lem:fixed-skew-lower}, which are admissible
even at $a=1$.
\end{proof}

\subsection{A counterexample under arbitrary dependence}\label{sec:qcs-proof}

\begin{proof}[Proof of Proposition~\ref{prop:qcs-counterexample}]
Set $d=d_n$, choose $b=2$, and take $\lambda>0$ such that
$\lambda b^{2+\tau}<K-1$. Let $U_i$ be independent symmetric signs and
let $B_{ij}$ be independent Bernoulli variables with success probability
$p=\lambda/(nd)$, independent of the signs. Set
\[
 Y_{ij}=U_i+b\sqrt n(B_{ij}-p),\qquad
 v=b^2np(1-p),\qquad X_{ij}=Y_{ij}/\sqrt{1+v}.
\]
Each coordinate is centered and has variance one. The common
off-diagonal correlation is $(1+v)^{-1}$. Also
\[
 \E|X_{ij}|^{2+\tau}\longrightarrow1+\lambda b^{2+\tau}<K.
\]
Indeed, on $B_{ij}=0$ the unscaled moment tends to one, while the
$B_{ij}=1$ contribution is
$p(b\sqrt n)^{2+\tau}(1+o(1))\to\lambda b^{2+\tau}$.
The common scale does not change the self-normalized ratios.
For the finitely many initial $n$ at which these bounds do not yet hold,
use independent symmetric signs.

Write $W_n=n^{-1/2}\sum_iU_i$ and $N_j=\sum_iB_{ij}$.
Call a Bernoulli success a mark. The total number of marks
converges in distribution to a Poisson variable with mean
$\lambda$, and
\[
 \Pp(\max_jN_j\ge2)
 \le d\binom n2p^2\longrightarrow0.
\]
Thus, with probability tending to one, the largest $N_j$ is zero if
there is no mark and one otherwise. The signs are independent of all
marks. Since the denominators are positive for all sufficiently large
$n$, the event $\max_jT_{n,j}\le0$ is exactly
\[
 W_n+b\max_jN_j-b\lambda/d\le0.
\]
The scalar CLT gives
\begin{equation}\label{eq:qcs-cdf-limit}
 \Pp(\max_jT_{n,j}\le0)
 \longrightarrow \tfrac12e^{-\lambda}
                  +(1-e^{-\lambda})\Phi(-b).
\end{equation}

On the event that each coordinate receives at most one mark, the
empirical second moment of the unscaled observations $Y_{ij}$ is
$1+o(1)$ in an unmarked column and $1+b^2+o(1)$ in a marked column,
uniformly over all columns. To check this,
expand
\[
 \frac1n\sum_iY_{ij}^2
 =1+\frac{2b}{\sqrt n}\sum_iU_i(B_{ij}-p)
       +b^2\sum_i(B_{ij}-p)^2.
\]
The middle term is uniformly $o(1)$ when $N_j\le1$, and the last term
is $b^2N_j+o(1)$. With probability tending to one there is an unmarked
coordinate. If a mark is present, the bounds $|T_{n,j}|\le t$ for all $j$
require $W_n$ to be within $o(1)$ of both $[-t,t]$ and
$[-b-t\sqrt{1+b^2},-b+t\sqrt{1+b^2}]$.
For $t=1/2$ and $b=2$ these intervals are disjoint. Consequently,
\begin{equation}\label{eq:qcs-abs-cdf-limit}
 \Pp(\max_j|T_{n,j}|\le1/2)
 \longrightarrow e^{-\lambda}\{2\Phi(1/2)-1\}.
\end{equation}

For independent standard normals $Z_0,Z_1,\ldots,Z_d$, the matching
Gaussian vector is
\[
 G_{P,j}=\frac{Z_0+\sqrt v\,Z_j}{\sqrt{1+v}}.
\]
Here $v\asymp d^{-1}$, so
$\sqrt v\max_j|Z_j|\to0$ in probability. Its maximum converges to $Z_0$
and its maximum absolute value converges to $|Z_0|$. Together with \eqref{eq:qcs-cdf-limit}--\eqref{eq:qcs-abs-cdf-limit},
these limits yield the two strict inequalities in the proposition.
The respective probability ratios converge to
\[
 1+(1-e^{-\lambda})(2\Phi(2)-1),\qquad
 1+(1-e^{-\lambda})\frac{2\Phi(1/2)-1}{2\barPhi(1/2)}.
\]
The independent Gaussian component has variance $v/(1+v)>0$,
which makes the covariance matrix positive definite. One may also index the row
laws only by dimension: take $n_d=\lceil d^{2/\tau}\rceil$ and fix the
constructed law $P_d$ at that sample size. The same limits hold along
this subsequence, so the example does not require changing a fixed
$d$-dimensional law as its sample size varies.
\end{proof}

\section{Shrinking correlation gaps}\label{sec:shrinking-proofs}

We prove the upper and lower bounds in Theorem~\ref{thm:main-critical}.
For the upper bound, we retain the sign of the cubic term in the
Lindeberg expansion and estimate its expectation by Gaussian integration
by parts. For the lower bound, we construct a sample with rare large
observations and compare its maximum with a Gaussian process indexed by
a spherical net. The two bounds determine the power of the correlation
gap; their logarithmic factors differ.

\subsection{Gaussian comparison at a shrinking correlation gap}
\label{sec:higher-boundary}
\label{sec:log-gap-improvement}

Throughout this subsection, $q>3$ and $K<\infty$ are fixed.
Lemma~\ref{lem:reductions} expresses each rectangle event as a system of
$m\le2d$ one-sided inequalities in signed coordinates. These coordinates form
two groups, according to their signs. Correlations within either group
are at most $1-a$; correlations between the groups are unrestricted.
Choose a real parameter
\[
 L\ge \log(2nd),\qquad L=O(n^{1/3}).
\]
The freedom to take $L$ larger than $\log(2nd)$ will ensure that the
estimates below also apply when the dimension is small.

We first replace each observation by an auxiliary pair. Its numerator
has the same joint moments through order three as the original observation
and contains an independent Gaussian component with covariance $R/2$.
The numerator construction is due to \citet[Section~4]{GZ2026}. We choose
the denominator contribution to match the additional moments involving
the denominator that enter the Taylor expansion of the ratio through
order three. Such mixed moments also appear in scalar expansions for general
numerator--denominator pairs; see \citet[Theorem~2.1]{GSS2022}.
These moment identities cancel all terms through degree three in the
first replacement stage.

We then replace the auxiliary pairs by Gaussian numerators and constant
denominator contributions. In this second stage, terms through degree
two still cancel, but the expected cubic term must be estimated. The
independent Gaussian component in the retained rows permits integration
by parts. An estimate for auxiliary Gaussian variables at indices
selected by the soft-maximum weights gives the logarithmic improvement
in the upper bound.

We transform the coordinates with the increasing convex function
\[
 g(t)=\frac{(t+\sqrt{t^2+4})^2}{4}.
\]
It satisfies $0\le g(t)-t_+^2\le 2$, $0<g'(t)\le C(1+t_+)$, and
all derivatives of order at least two are bounded.
For the finite thresholds $b_j$ of the signed coordinates, put
\[
 z_j(t)=g(t_j)-g(b_j),\qquad
 Y(t)=\max_j z_j(t),\qquad
 A_\ell(z)=\ell^{-1}\log\sum_j e^{\ell z_j}.
\]
Lemma~\ref{lem:core-localization}, applied separately to the two groups,
gives
\begin{equation}\label{eq:newcrit-localization}
 \Pp\{Y(G)\le s,\ A_\ell(z(G))>s+u\}
 \le \frac{C_\alpha}{a\ell u},
 \qquad \ell\ge C_\alpha/a,\quad \ell u\ge 2\log 2,
\end{equation}
whenever $\Pp\{Y(G)\le s\}\ge\alpha>0$.
Indeed, the log-sum-exponential over both groups exceeds the larger of
the two group values by at most $\ell^{-1}\log2$.
At level $s$, the threshold in coordinate $j$ is
$g^{-1}(g(b_j)+s)$. Its argument is positive whenever
$\Pp\{Y(G)\le s\}>0$.
Lemma~\ref{lem:core-score-density} also yields
\[
 \Pp\{u<Y(G)\le v\}\le C(v-u),\qquad u<v,
\]
with an absolute constant.

\subsubsection{Smooth upper and lower approximations to a rectangle indicator}

\begin{lemma}[Bounded smooth bracketing functions]
\label{lem:newcrit-brackets}
Suppose $\Pp\{Y(G)\le0\}\in[\varepsilon,1-\varepsilon]$.
For sufficiently small fixed $\tau>0$, let
$\lambda=\Lambda/a$, with $\Lambda$ sufficiently large and
$\lambda\tau\ge16\log2$. There are uniformly bounded smooth functions
$F^-,F^+$ satisfying
\begin{equation}\label{eq:newcrit-brackets}
 F^-(z)\le\1_{\{\max_jz_j\le0\}}\le F^+(z),
\end{equation}
and
\begin{equation}\label{eq:newcrit-calibration}
 \E[F^+(z(G))-F^-(z(G))]
 \le C_\varepsilon\{\tau+(a\lambda\tau)^{-1}
                                      +e^{-\lambda\tau/2}\}.
\end{equation}
Each function is a constant plus a linear combination of two
log-sum-exponential functions applied to convex increasing coordinate
scores. The coefficients are bounded in terms of $\tau$. We may include
an additional constant score in either sum.
\end{lemma}
\begin{proof}
In this proof, write $Y=\max_jz_j$ and $A_\ell=A_\ell(z)$, and evaluate
them at $z(G)$ when taking probabilities or expectations.
Choose $\tau>0$ so that $20C\tau\le\varepsilon/2$, where $C$ is the
absolute constant in the preceding Gaussian interval bound. Put
\[
 D_c(z)=\lambda^{-1}\log\left(e^{\lambda c}
                                     +\sum_j e^{\lambda z_j}\right),
 \qquad
 F^-(z)=\frac{D_{-\tau}(z)-D_{-2\tau}(z)}{\tau}
                     -\frac{e^{-\lambda\tau}}{\lambda\tau}.
\]
The quotient in $F^-$ lies between zero and one. If $Y>0$, then
$\sum_j e^{\lambda z_j}\ge1$, so
$D_{-\tau}-D_{-2\tau}\le e^{-\lambda\tau}/\lambda$.
This proves the first inequality in \eqref{eq:newcrit-brackets}. On
$Y\le-3\tau$ and $A_\lambda\le-5\tau/2$, the quotient is at least
$1-e^{-\lambda\tau/2}/(\lambda\tau)$.
The event $-3\tau<Y\le0$ has probability at most
$C_\varepsilon\tau$. Moreover, $\Pp(Y\le-3\tau)\ge\varepsilon/2$, so
\eqref{eq:newcrit-localization} bounds the probability of
$\{Y\le-3\tau,\ A_\lambda>-5\tau/2\}$ by
$C_\varepsilon/(a\lambda\tau)$.

For the upper bound, choose a convex $C^\infty$ function $u$ with
$u=0$ on $(-\infty,0]$, $0\le u'\le1$, and $u'=1$ on
$[\tau,\infty)$. Define
\[
 v_1(s)=s+u(s),\qquad v_0(s)=s+u(s-\tau),\qquad
 F^+(z)=1-\frac{A_\lambda(v_1(z))-A_\lambda(v_0(z))}{\tau}.
\]
The difference $v_1-v_0$ lies in $[0,\tau]$, vanishes on $s\le0$,
and equals $\tau$ on $s\ge2\tau$. Hence $0\le F^+\le1$ and
$F^+=1$ on $Y\le0$. To bound it outside that set, integrate the
derivative along $v_s=(1-s)v_0+sv_1$. Relative to the weights
proportional to $e^{\lambda z_j}$, the new weights are tilted by the
nondecreasing function $e^{\lambda(v_s(z_j)-z_j)}$.
Tilting by a nondecreasing function decreases the mass of a lower level
set. This follows from the nonnegative covariance of two nondecreasing
functions under the probability measure defined by the original weights.
Thus
\[
 F^+(z)\le
 \frac{\sum_{z_j<2\tau}e^{\lambda z_j}}{\sum_j e^{\lambda z_j}}
 \le \exp\{\lambda(\tau+A_{\lambda/2}(z)/2-Y)\}.
\]
The strip $0<Y\le10\tau$ has probability $O_\varepsilon(\tau)$.
Partition $Y>10\tau$ into $(y,3y/2]$ with
$y=10\tau(3/2)^k$, $k\ge0$. If
$A_{\lambda/2}\le8y/5$, the last bound is at most $e^{-\lambda\tau}$.
On the remaining event, apply \eqref{eq:newcrit-localization} with
level $3y/2$ and excess $y/10$. Summing these probability bounds over $k$ and using $0\le F^+\le1$
controls $\E[F^+(z(G))\mathbf1_{\{Y(G)>0\}}]$. Together with the estimate for $F^-$, this proves
\eqref{eq:newcrit-calibration}.
All nonconstant coordinate scores have the form $v(g(t)-g(b))$, where
$v$ is convex and increasing, $1\le v'\le2$, and its higher derivatives
are bounded in terms of $\tau$.
\end{proof}

\subsubsection{Row replacement with an independent Gaussian component}

Represent one observation by a pair $(y,h)$, where $y$ contributes to
the numerator and $h_j\ge0$ contributes to the squared denominator.
The original pair is $(X,X^2)$, with the square taken coordinatewise;
the final Gaussian pair is $(G,\mathbf1)$, with $G\sim N(0,R)$.
After deleting row $i_0$ from an array of such pairs, define
\[
 Q_j=\sum_{i\ne i_0}h_{ij},\qquad
 t_j=Q_j^{-1/2}\sum_{i\ne i_0}y_{ij},\qquad
 d_j=\sqrt{n/Q_j},\qquad
 \mathcal B=\{\min_j Q_j\ge n/4\}.
\]
An undefined ratio outside $\mathcal B$ may be assigned any value;
there we use only bounded test differences.

Use the numerator construction of \citet[Section~4.1]{GZ2026} and
choose its squared-denominator contribution as follows:
\begin{equation}\label{eq:newcrit-augmented}
 Y_j=2IX_j+G_j/\sqrt2,\qquad
 H_j=4IX_j^2+\frac12,\qquad
 I\sim\operatorname{Bernoulli}(1/8).
\end{equation}
Here $I,X,G$ are independent, and the same indicator $I$ is used for
every coordinate of a row. The independent Gaussian component
$G/\sqrt2$ has covariance $R/2$. The choice of $H$ gives $\E H_j=1$
and $\E[Y_iH_j]=\E[X_iX_j^2]$, the additional identities needed for the
denominator terms in the ratio expansion. All rows are independent,
including in arrays that combine original, auxiliary, and Gaussian
pairs; coordinates within a row need not be independent.

\begin{lemma}[Bounds after deleting one observation]\label{lem:shr-retained}
Uniformly over arrays formed from the three types of pairs above and
over nonempty deterministic sets $\mathcal I$ of signed coordinates,
\begin{equation}\label{eq:newcrit-retained-tails}
 \Pp(\mathcal B^c)\le Cde^{-cn},\qquad
 \Pp\{\mathcal B,\ \max_{j\in\mathcal I}|t_j|>x\}
 \le C|\mathcal I|e^{-cx^2},\qquad x\ge1.
\end{equation}
Consequently, for every fixed $p<\infty$,
\begin{equation}\label{eq:newcrit-retained-moments}
 \E\big[\mathbf1_{\mathcal B}(1+\max_j|t_j|)^p\big]
 \le C_p L^{p/2}.
\end{equation}
\end{lemma}
\begin{proof}
Every denominator contribution is nonnegative, has mean one, and has a
uniformly bounded moment of order $q/2$. Truncation at a sufficiently large fixed constant leaves
the mean at least $3/4$. For each coordinate, the truncated contributions
are independent across rows. Hoeffding's inequality and a union bound
over the coordinates give the first estimate; the constants can be
adjusted to include the finitely many small values of $n$.

For every real $v$,
$e^{v-v^2/2}\le1+v+Cv^2$.  Centering therefore gives, with
$S_j=\sum_{i\ne i_0}y_{ij}$,
\[
 \E\exp\{uS_j-u^2Q_j/2\}\le e^{Cnu^2}.
\]
For an auxiliary row, integrating out the independent Gaussian component
cancels the exponential factor from the constant $1/2$ in $H_j$. The
remaining expectation is $\E\exp\{2uIX_j-2u^2IX_j^2\}$, to which the
same inequality applies. For a Gaussian row, the exponential has
expectation one.
On a shell $b\le Q_j<2b$, where $b\ge n/4$, choose
$u=\pm c_0x/\sqrt b$ with a sufficiently small fixed $c_0$.
The exponential Markov inequality bounds each event
$\pm S_j>x\sqrt{Q_j}$ on this shell by $e^{-cx^2}$.
Summing over the first $1+\lceil x^2\rceil$ dyadic shells gives the
same exponential bound after reducing $c$. Markov's inequality and
$\E Q_j=n-1$ bound the probability of the remaining shells by
$C2^{-x^2}$.
A union bound over $\mathcal I$ proves the second inequality.
Integrating the tail bound, capped at one, proves
\eqref{eq:newcrit-retained-moments}.
\end{proof}

We insert an independent pair $(y,h)$ along the path
\begin{equation}\label{eq:newcrit-insertion}
 x_j(v)=\left(t_j+\frac{vd_jy_j}{\sqrt n}\right)
        \left(1+\frac{v^2d_j^2h_j}{n}\right)^{-1/2},
 \qquad 0\le v\le1.
\end{equation}
The first three derivatives of $x_j(v)$ at $v=0$ are
$d_jy_j/\sqrt n$, $-t_jd_j^2h_j/n$, and
$-3d_j^3y_jh_j/n^{3/2}$.
For either bracketing function, write $f(t)=F^\pm(z(t))$.
In the resulting Taylor expansion, assign degree one to $y$ and degree
two to $h$, as prescribed by their powers of $v$ in
\eqref{eq:newcrit-insertion}. The terms of degree at most two depend only on
$\E y_j$, $\E y_iy_j$, and $\E h_j$.
Conditional on the remaining rows, the expected cubic coefficient for
an original incoming row, excluding the factor $1/(6n^{3/2})$, is
\begin{equation}\label{eq:newcrit-cubic}
\begin{split}
 \mathcal K(t)
 &=\sum_{i,j,k}d_id_jd_k f_{ijk}(t)\kappa_{ijk}
   -3\sum_{i,j}d_id_j^2t_jf_{ij}(t)\kappa_{ijj}
   -3\sum_jd_j^3f_j(t)\kappa_{jjj},\\
 \kappa_{ijk}&=\E[X_iX_jX_k].
\end{split}
\end{equation}
The last two terms arise from differentiating the denominator in
\eqref{eq:newcrit-insertion}; we will estimate their sum.

\begin{lemma}[The third-order remainder]\label{lem:newcrit-remainder}
Put $r=q\wedge4>3$. For either bracketing function from
Lemma~\ref{lem:newcrit-brackets}, expand $v\mapsto f(x(v))$ through
degree three at $v=0$ and evaluate the remainder at $v=1$.
For each realization of the remaining rows in
$\mathcal B\cap\{\max_j|t_j|\le C_0\sqrt L\}$, the conditional expectation
of the absolute remainder is at most
\begin{equation}\label{eq:newcrit-remainder}
 Ce^{C\lambda}(1+\lambda)^C(L/n)^{r/2}+O(n^{-12}).
\end{equation}
The estimate holds for each of the three types of incoming pairs.
When taking the expectation of the Taylor polynomial on $\mathcal B$,
removing the restriction $\max_j|t_j|\le C_0\sqrt L$ changes the result
by $O(n^{-12})$ per row when $C_0$ is sufficiently large.  The complement of $\mathcal B$ contributes at most
$Cde^{-cn}$ to a bounded test difference.
The constants may depend on $q,K,\tau$.
\end{lemma}
\begin{proof}
For an original row put $e_j=C\sqrt{L/n}(1+|X_j|)$.
On the event specified in the lemma, differentiation of
\eqref{eq:newcrit-insertion} and the chain rule give
\[
 \left|\frac{d^s}{dv^s}\psi_j(x_j(v))\right|\le C_s e_j^s,
 \qquad 1\le s\le4,\quad e_j\le1,
\]
for every scalar coordinate score $\psi_j$ in a bracketing function.
Indeed, $|X_j|/\sqrt{Q_j}$ and $\sqrt{X_j^2/Q_j}$ are at most
$Ce_j/\sqrt L$, the derivative of each score is at most $C\sqrt L$,
and each further chain-rule term is bounded by the corresponding
power of $e_j$.  The derivatives at zero satisfy the same polynomial
bounds without restricting $e_j$.
For every value of the incoming observation, Cauchy--Schwarz gives
$(x_j(v)_+)^2\le(t_j{}_+)^2+1$, and hence
$g(x_j(v))-g(t_j)\le3$.
The increasing transformations $v_0,v_1$ multiply this upper bound by
at most two.

For each log-sum-exponential term, let $w_{uj}$ be its weights before
inserting the new observation, and write
\[
 S_u(v)=\sum_jw_{uj}
  \exp\{\lambda[\psi_{uj}(x_j(v))-\psi_{uj}(t_j)]\},
 \qquad M_r=\sum_{u,j}w_{uj}e_j^r.
\]
Here $u$ indexes the finitely many log-sum-exponential terms and
$\psi_{uj}$ is the corresponding coordinate score. Any additional
constant score has zero increment.
H\"older's inequality applied to the weights bounds each degree-$s$
Taylor coefficient of $S_u$, $s\le3$, by
$C(1+\lambda)^s M_r^{s/r}$.
The remainder in the scalar exponential expansion is bounded by
$e^{C\lambda}(1+\lambda)^C e_j^r$:
for $e_j\le1$, use the fourth derivative and $e_j^4\le e_j^r$;
for $e_j>1$, use the bound on the score increase and the Taylor polynomial at zero.
The bounds are uniform over $0\le v\le1$.

If $e^{C\lambda}(1+\lambda)^C M_r^{1/r}$ is sufficiently small,
every $S_u(v)$ stays above $1/2$.  Expand its logarithm and retain
total degree at most three.  The remainder is at most
$e^{C'\lambda}(1+\lambda)^{C'}M_r$.
If this smallness condition fails, boundedness of the bracketing function
and the coefficient bounds give the same estimate after enlarging $C'$.
Boundedness is used only for the combined bracketing function; its
individual log-sum-exponential terms need not be bounded.
Independence of the incoming row gives
$\E[M_r\mid\text{retained rows}]\le C(L/n)^{r/2}$.

For a Gaussian row first restrict to
$\max_j|G_j|\le C_2\sqrt L$ and use
$e_j=C\sqrt{L/n}(1+|G_j|)$.
The increase of a score is $O(L/\sqrt n)=o(1)$.
For an auxiliary row use
$e_j=C\sqrt{L/n}(1+|X_j|+|G_j|)$.
First insert $2IX_j$ and its square into the numerator and squared
denominator, respectively. The preceding positive-part bound then applies.  The additional $1/2$ in the denominator decreases
the positive part, while the Gaussian numerator changes the coordinate
by $O(\sqrt{L/n})$.
The additional score increase is $O(1+L/\sqrt n)$, uniformly in $X_j$
and $v$.  The local derivative bounds follow from
$|Y_j|\le C(|X_j|+|G_j|)$ and $H_j\le C(1+X_j^2)$.
On $\{\max_j|G_j|>C_2\sqrt L\}$, use boundedness of the bracketing
function and apply H\"older's inequality with exponent $r/3>1$ to the
degree-three Taylor polynomial.
Gaussian tails make the contribution $O(n^{-12})$ after increasing $C_2$.

Finally, for $1\le s\le3$, differentiation of the log-sum-exponential terms gives
\begin{equation}\label{eq:shr-absolute-derivatives}
 \sum_{j_1,\ldots,j_s}|f_{j_1\cdots j_s}(t)|
 \le C(1+\lambda)^{s-1}(1+\|t\|_\infty)^s.
\end{equation}
Lemma~\ref{lem:shr-retained}, the displayed derivative bound, and the
bounded absolute moments of the incoming row through order three give
the claimed estimate on $\{\max_j|t_j|>C_0\sqrt L\}$.
On $\mathcal B^c$, we use boundedness of the test difference without
expanding it.
\end{proof}

\begin{lemma}[Two stages of moment matching]
\label{lem:newcrit-reserve}
Replace the original pairs one at a time by
\eqref{eq:newcrit-augmented}, and then replace these auxiliary pairs
one at a time by $(G,\mathbf1)$. In the first stage, conditional on the
remaining rows, the expectations of all terms through weighted degree
three agree for the outgoing and incoming pairs. In the second stage,
the expectations still agree through degree two, but the uncancelled
cubic coefficient is \eqref{eq:newcrit-cubic}.
Conditional on the indicators and original observations in the rows
that remain after deletion, the statistic in the second stage has the form
\begin{equation}\label{eq:newcrit-fixed-gaussian}
 t_j=u_j+\theta d_jZ_j,\qquad Z\sim N(0,R),\qquad
 \theta^2=\frac{k+(n-1-k)/2}{n}\ge\frac{n-1}{2n},
\end{equation}
where $k$ is the number of remaining Gaussian pairs.
Under this conditioning, $u_j$, $d_j$, and $\mathbf1_{\mathcal B}$
are fixed, and the remaining randomness in $t$ is Gaussian.  For $n\ge2$, $\theta\ge1/2$, and on $\mathcal B$,
$d_j\le2$.
If $\lambda=o(\log n)$ and $L=O(n^{1/3})$, the total remainder is bounded by
\begin{equation}\label{eq:shr-two-sweep-remainder}
 Cn e^{C\lambda}(1+\lambda)^C(L/n)^{r/2}+o(1),
 \qquad r=q\wedge4.
\end{equation}
\end{lemma}
\begin{proof}
Since $4\E I=1/2$ and $8\E I=1$, independence and centering give
\[
 \E Y_j=0,\quad \E H_j=1,\quad
 \E Y_iY_j=R_{ij},\quad
 \E Y_iY_jY_k=\kappa_{ijk},\quad
 \E Y_iH_j=\kappa_{ijj}.
\]
These identities make the expectations of all terms of weighted degree
at most three agree for the original and auxiliary pairs in the
expansion along \eqref{eq:newcrit-insertion}. For the Gaussian pair,
the moments through degree two still agree, while both the mixed third
moments of the numerator and the mixed numerator--denominator moments
vanish.
Lemma~\ref{lem:newcrit-remainder} applies throughout both stages.
In the second stage, the Gaussian numerator components of the remaining rows
have covariance $[k+(n-1-k)/2]R$, while the denominators depend only
on the retained indicators and original variables.
This gives \eqref{eq:newcrit-fixed-gaussian}.
There are $2n$ replacements. Summing the remainder bounds and the
contributions from the exceptional events proves
\eqref{eq:shr-two-sweep-remainder}.
Under the stated assumptions, the leading term is bounded by
$n^{1-r/3+o(1)}=o(1)$.
A finite-dimensional factor representation proves the same assertions
when $R$ is singular.
\end{proof}

\subsubsection{Bounding the cubic term by Gaussian integration by parts}
\label{subsec:newcrit-twofold}
Gaussian covariance and entropy estimates for soft-maximum weights
appear in \citet[Lemmas~11--12]{GZ2026}. The reduction of a weighted
Gaussian matrix to a tilted mean and a signed trace is given in
\citet[Section~H.2]{CCKKSUP2026}, the supplement to \citet{CCKK2026}.
Here the coordinate scores are nonlinear and convex. Conditioning on
one Gaussian coordinate fixes the numerator of its soft-maximum weight
and makes the remaining tilt log-concave. The proof uses this
conditional structure to obtain a bound uniform in the dimension of
the Gaussian representation space. Write $\|\cdot\|_*$ for the trace norm and
$\|\cdot\|_{\rm op}$ for the operator norm.

\begin{lemma}[A third-derivative bound uniform in $\dim V$]\label{lem:newcrit-twofold}
Let $v_1,\ldots,v_m$ be unit vectors in a finite-dimensional real Hilbert
space $V$, and let $Z$ be standard Gaussian in $V$.
Put $Z_j=\langle v_j,Z\rangle$ and
$t_j=u_j+\theta d_jZ_j$, where $u_j,d_j\in\mathbb R$ and $\theta>0$.
Let $h_j$ be $C^3$ convex functions with bounded first derivatives
and polynomially growing third derivatives, and put
\[
 F(t)=\lambda^{-1}\log\sum_j e^{\lambda h_j(t_j)},\qquad \lambda>0.
\]
For arbitrary linear operators $T_i:V\to V$, set
$A=\max_i\{|d_i|\,\|h_i'\|_\infty\|T_i\|_{\rm op}\}$.
Then
\begin{equation}\label{eq:newcrit-twofold}
 \left|\E_Z\sum_{i,j,k}d_id_jd_kF_{ijk}(t)
                  \langle v_j,T_iv_k\rangle\right|
 \le CA\theta^{-2}\log(2m).
\end{equation}
The constant is independent of the smoothing parameter $\lambda$,
the dimension of $V$, and the intercepts $u_j$. A constant score may be
included by taking $d_i=0$.
\end{lemma}

\begin{proof}
Only the symmetric part of $T_i$ contributes to the sum in
\eqref{eq:newcrit-twofold}, and
symmetrization does not increase its operator norm.  Write $r=\dim V$, and let
$w_i=e^{\lambda h_i(t_i)}/\sum_je^{\lambda h_j(t_j)}$.
Twice integrating by parts in any orthonormal basis of $V$ gives, for each fixed $i$,
\begin{equation}\label{eq:newcrit-two-ibp-identity}
 \theta^2\mathbb E_Z\sum_{j,k}d_jd_kF_{ijk}(t)
               \langle v_j,T_iv_k\rangle
 =\mathbb E_Z\left[F_i(t)\operatorname{tr}\{T_i(ZZ^\top-I_V)\}\right].
\end{equation}
Since $F_i=w_ih_i'(t_i)$ and $h_i'(t_i)$ depends only on $Z_i$, define
\[
 p_i(x)=\mathbb E[w_i\mid Z_i=x],\qquad
 Q_i(x)=\mathbb E[w_i(ZZ^\top-I_V)\mid Z_i=x].
\]
The duality between the trace norm and the operator norm reduces the
claim to
\begin{equation}\label{eq:newcrit-trace-bound}
 \sum_i\mathbb E\|Q_i(Z_i)\|_*\leq C\log(2m).
\end{equation}

Fix $i$ and $x$. Conditional on $Z_i=x$, the vector $Z$ is Gaussian
with mean $xv_i$ and identity covariance on $v_i^\perp$. Tilt this law
by the density $w_i/p_i(x)$, and denote the mean and covariance of the
tilted law in $V$ by $\mu_i(x)$ and $C_i(x)$.
The numerator $e^{\lambda h_i(t_i)}$ is constant on this affine space.
The logarithm of the denominator is convex: the log-sum-exponential
function is convex and nondecreasing in each argument, and each
$h_j(u_j+\theta d_j\langle v_j,Z\rangle)$ is convex. Translate the
affine space to $v_i^\perp$. The negative log-density of the tilted law
then has Hessian at least the identity. Applying the Brascamp--Lieb variance inequality
\citep{BL1976,BK2020} to linear functions gives
\[
 0\preceq C_i(x)\preceq I_V,
 \qquad C_i(x)v_i=0.
\]
This also holds when $v_i^\perp$ is zero-dimensional. Consequently,
\[
 Q_i=p_i\{C_i-I_V+\mu_i\mu_i^\top\},\qquad
 \|Q_i\|_*\leq p_i\{r-\operatorname{tr}C_i+\|\mu_i\|^2\}.
\]
To remove the dependence on $r=\dim V$, note that $\sum_iw_i=1$ implies
$\sum_i\mathbb E Q_i(Z_i)=\mathbb E(ZZ^\top-I_V)=0$, so taking traces gives the exact identity
\begin{equation}\label{eq:newcrit-rank-cancellation}
 \sum_i\mathbb E[p_i(Z_i)\{r-\operatorname{tr}C_i(Z_i)\}]
 =\sum_i\mathbb E[p_i(Z_i)\|\mu_i(Z_i)\|^2].
\end{equation}

Since $0<w_i\le1$, the tilted conditional law has relative entropy at
most $\log(1/p_i(x))$ with respect to the conditional Gaussian law.
The entropy variational inequality, applied to
$\langle b,Z-xv_i\rangle$, $b\in v_i^\perp$, yields
\[
 \langle b,\mu_i(x)-xv_i\rangle-\|b\|^2/2
 \leq\log(1/p_i(x)).
\]
Optimizing over $b$ gives
$\|\mu_i(x)-xv_i\|^2\leq2\log(1/p_i(x))$.
Since this displacement is orthogonal to $v_i$,
\begin{align*}
 \sum_i\mathbb E[p_i(Z_i)\|\mu_i(Z_i)\|^2]
 &\leq\mathbb E\sum_iw_iZ_i^2
       +2\sum_i\mathbb E[p_i(Z_i)\log(1/p_i(Z_i))]\\
 &\leq\mathbb E\max_iZ_i^2
       +2\sum_i\bar p_i\log(1/\bar p_i)
 \leq C\log(2m),
\end{align*}
where $\bar p_i=\mathbb Ew_i$ and $\sum_i\bar p_i=1$. The second
inequality uses the concavity of $p\log(1/p)$. A union bound for the
Gaussian tails controls the maximum without requiring independent
coordinates.
Combine this estimate with \eqref{eq:newcrit-rank-cancellation}
to prove \eqref{eq:newcrit-trace-bound}.
Applying the trace-norm inequality conditionally in
\eqref{eq:newcrit-two-ibp-identity} bounds the left side of
\eqref{eq:newcrit-twofold} by
\[
 \theta^{-2}\sum_i
 \E\!\left[|d_ih_i'(t_i)|\,\|T_i\|_{\rm op}\|Q_i(Z_i)\|_*\right]
 \le A\theta^{-2}\sum_i\E\|Q_i(Z_i)\|_*.
\]
This proves \eqref{eq:newcrit-twofold}.
Polynomial growth of $h_j'''$ also implies polynomial growth of $h_j''$.
Together with the bounded first derivatives, these conditions justify
both integrations by parts by Gaussian integrability, or by compact
cutoffs followed by dominated convergence.
\end{proof}

\begin{remark}\label{rem:newcrit-operator-realization}
The proof gives the more precise bound
\[
 2A\theta^{-2}\left\{
 \E\sum_iw_iZ_i^2+
 2\sum_i\E[p_i(Z_i)\log(1/p_i(Z_i))]
 \right\}.
\]
The factor $\log(2m)$ bounds both the weighted second moment and the
sum of the terms $\E[p_i(Z_i)\log(1/p_i(Z_i))]$. The latter is bounded
by the entropy of the probability vector $(\E w_i)_i$.
This refinement alone does not enlarge the dimension range proved below.

For centered unit-variance variables $X_1,\ldots,X_m$, take $V$ to be
their span in $L^2(P)$ and set $v_j=X_j$. Given a bounded real random
variable $b_i$, define $T_i$ by multiplication by $b_i$ followed by
orthogonal projection onto $V$. Then $T_i$ is self-adjoint,
$\|T_i\|_{\rm op}\leq\|b_i\|_\infty$, and
$\langle v_j,T_iv_k\rangle=\mathbb E[b_iX_jX_k]$.
No independence or centering assumption on $b_i$ is needed. For a
fixed finite linear combination of log-sum-exponential functions, apply
the lemma term by term and multiply the bound by the sum of the absolute
coefficients.
\end{remark}

To estimate the cubic term, modify $g$ smoothly and convexly so that it
agrees with $g$ up to $C_0\sqrt L$, is affine to the right of
$C_0\sqrt L+1$, and has slope at most $C\sqrt L$.
To construct the extension, multiply $g''$ on this unit interval
by a smooth cutoff between zero and one, equal to one near the left
endpoint and zero beyond the right endpoint, and integrate twice with
matching value and slope.  Leave every threshold constant $g(b_j)$
unchanged.  The resulting scores $h_j(t)=v(g(t)-g(b_j))$ satisfy
\begin{equation}\label{eq:newcrit-modified-derivatives}
 \max_j\|h_j'\|_\infty\le C\sqrt L,\qquad
 \max_j\|h_j^{(s)}\|_\infty\le C_sL^{s/2},\quad 2\le s\le4.
\end{equation}
The original and modified cubic expressions agree on
$\max_j|t_j|\le C_0\sqrt L$.
Their absolute derivative bounds and Lemma~\ref{lem:shr-retained}
show that, after multiplication by $n^{-3/2}$ and summation over the
rows, the expected difference is $o(1)$ when $C_0$ is sufficiently large.
The initial and final tests and the Taylor remainder still use the
original score. Below, $\E_0 U=\E[\mathbf1_{\mathcal B}U]$ denotes
expectation restricted to $\mathcal B$, without conditioning on that event.

We also need two identities obtained by Gaussian integration by parts.
For $U_{ij}\in L^2(P)$, realize jointly
$Z_k=W(X_k)$ and $W_{ij}=W(U_{ij})$ in an isonormal Gaussian process.
Condition on the indicators and original observations in the remaining
rows, apply Gaussian integration by parts, and then average over those
variables. This gives
\begin{equation}\label{eq:newcrit-one-ibp}
 \E_0\sum_{i,j,k}d_id_jd_k f_{ijk}(t)\E[U_{ij}X_k]
 =\theta^{-1}\E_0\sum_{i,j}d_id_jf_{ij}(t)W_{ij}.
\end{equation}
For each convex log-sum-exponential function,
\begin{equation}\label{eq:shr-gibbs-hessian}
 F_{ij}=\mathbf1_{\{i=j\}}w_ih_i''
       +\lambda(\mathbf1_{\{i=j\}}w_i-w_iw_j)h_i'h_j',
 \qquad \sum_{i,j}|F_{ij}|\le C\lambda L.
\end{equation}
Here and in the subsequent application $\lambda=\Lambda/a\ge1$.

When $q\ge4$, put $W_j'=W(X_j^2)$.  Differentiating the explicit
factor $t_j$ as well as $f_j(t)$ gives
\begin{equation}\label{eq:newcrit-denominator-identity}
 \E_0\left\{\sum_{i,j}d_id_j^2t_jf_{ij}(t)\kappa_{ijj}
                 +\sum_jd_j^3f_j(t)\kappa_{jjj}\right\}
 =\theta^{-1}\E_0\sum_jd_j^2t_jf_j(t)W_j'.
\end{equation}
For each log-sum-exponential term, $F_j=w_jh_j'$ and $\sum_jw_j=1$.
Apply this identity term by term to $f$. Under
\eqref{eq:newcrit-fixed-gaussian}, the absolute value of the right side
is at most $CL^{3/2}$, since
\[
 \|\mathbf1_{\mathcal B}\max_j|t_j|\|_2\le C\sqrt L,\qquad
 \|\max_j|W_j'|\|_2\le C\sqrt L.
\]
These $L^2$ estimates include expectation over the remaining rows;
no bound uniform over the conditional intercepts is needed.
The indices of $W$ are elements of $L^2(P)$ and need not be centered.

\begin{remark}\label{rem:shr-zero-cubic}
If $\kappa_{ijk}=0$ for every triple, the direct replacement from
$(X,X^2)$ to $(G,\mathbf1)$ has no cubic term.
Lemmas~\ref{lem:newcrit-brackets} and \ref{lem:newcrit-remainder}
then give convergence throughout $\log(2nd)=O(n^{1/3})$ when
$a\log n\to\infty$, for this smaller class and $q>3$.
Here the assumption requires all mixed third moments to vanish;
vanishing marginal third moments alone is insufficient.
\end{remark}

\subsubsection{The number of scores near the Gaussian maximum}

Gaussian perturbations and near-maximizing sets are studied by
\citet{CHAT2009}. Chatterjee's entropy bound
\citep[Theorem~12.4]{CHAT2014}, restated in
\citet[Theorem~3.2]{MS2021}, concerns a centered Gaussian process and
a level defined by its expected maximum. The sparsification theorem
of \citet{DNOS2026} approximates the supremum by a smaller shifted
process. The next lemma bounds the number of indices near the
realized maximum, uniformly over arbitrary affine intercepts.

\begin{lemma}[Number of indices near a Gaussian maximum]
\label{lem:newlog-nearmax}
Suppose $a\downarrow0$, $L\ge1$, $C_0,c_0>0$, and $N\le e^{C_0L}$.
Let
\[
 X_i=\langle v_i,Z\rangle+b_i,\qquad
 \|v_i\|\leq1,\qquad \|v_i-v_j\|\geq c_0\sqrt a\quad(i\ne j),
\]
where $Z$ is standard Gaussian in a finite-dimensional real Hilbert space and the $b_i$ are arbitrary deterministic intercepts. Fix $P\geq0$, $A>0$, and $B>0$. For $0\leq t\leq Aa\sqrt L\,\log^P(1/a)$, set
\[
 \mathcal A_t=\{i:X_i\geq\max_jX_j-t\}.
\]
There are constants $C,c>0$, depending only on $C_0,c_0,P,A,B$, such that
\begin{equation}\label{eq:newlog-nearmax}
 \mathbb P\{\log|\mathcal A_t|>C\frac{L\log\log(1/a)}{\log(1/a)}\}
 \leq C\log^{-B}(1/a)+Ce^{-cL}.
\end{equation}
The estimate is uniform over this range of $t$, over the intercepts,
and over deterministic nonempty subsets of the index set. An empty
subset contributes no indices.
\end{lemma}
\begin{proof}
By monotonicity, take $t=Aa\sqrt L\,\log^P(1/a)>0$.
If $L\le \log(1/a)/(8C_0)$, two elements in $\mathcal A_t$ would give a pair
with $|X_i-X_j|\le t$.
The variance of this difference is at least $c_0^2a$.
The one-dimensional Gaussian density bound and a union bound imply
\[
 \Pp\{|\mathcal A_t|\ge2\}
 \le CN^2t/\sqrt a
 \le Ca^{1/4}\log^{P+1/2}(1/a)=o(\log^{-B}(1/a)).
\]
On the complementary event there is a unique index in $\mathcal A_t$.
For the remaining case, rescale $L$ by a fixed constant so that
$N\le e^L$ and $L\ge c\log(1/a)$ for a fixed $c>0$.
For a deterministic nonempty set $S$, let
\[
 W(S)=\mathbb E'\max_{i\in S}\langle v_i,Z'\rangle,
\]
where $Z'$ is an independent standard Gaussian and $\mathbb E'$ denotes
expectation over $Z'$ only. For $u>0$, restricting the perturbed maximum
to $\mathcal A_t$ gives
\[
 \mathbb E\max_i\{X_i+u\langle v_i,Z'\rangle\}
 \geq\mathbb E\max_iX_i-t+u\mathbb EW(\mathcal A_t).
\]
The expectation on the left equals that of the maximum with its Gaussian
part multiplied by $\sqrt{1+u^2}$. For arbitrary intercepts, its increase over $\mathbb E\max_iX_i$ is at most
\[
 (\sqrt{1+u^2}-1)\mathbb E\max_i|\langle v_i,Z\rangle|
 \leq Cu^2\sqrt L.
\]
Take $u=(t/\sqrt L)^{1/2}$ to obtain
\begin{equation}\label{eq:newlog-width}
 \mathbb EW(\mathcal A_t)\leq C(t\sqrt L)^{1/2}.
\end{equation}

Choose $R=\log^{P+B+5}(1/a)$ and
$k=C_*L\log\log(1/a)/{\log(1/a)}$, with $C_*$ sufficiently large and
fixed. If $|\mathcal A_t|>e^k$, either a ball of radius $R\sqrt a$
centered at some $v_i$ contains the vectors corresponding to at least
$e^{k/2}$ indices in $\mathcal A_t$, or a greedy packing selects at least
$e^{k/2}$ such vectors at mutual distances greater than $R\sqrt a$.
In the second case, Sudakov--Fernique comparison gives
\[
 W(\mathcal A_t)\geq cR\sqrt{ak}.
\]
By \eqref{eq:newlog-width}, the probability of this alternative is at most
\[
 C\left(\frac{t\sqrt L}{R^2ak}\right)^{1/2}
 \leq\frac{C\log^{(P+1)/2}(1/a)}{R\sqrt{\log\log(1/a)}}
 \leq C\log^{-B}(1/a).
\]

In the first case, fix the center $v_j$ of such a ball and put
\[
 B_j=\{i:\|v_i-v_j\|\leq R\sqrt a\},\quad
 w_i=\frac{v_i-v_j}{R\sqrt a},\quad e_i=\frac{b_i}{R\sqrt a},\quad
 t'=\frac{t}{R\sqrt a},\quad \beta=\frac{k}{8t'}.
\]
The $w_i$ have norms at most one and mutual distances at least $c_0/R$. Consider
\[
 f(z)=\beta^{-1}\log\sum_{i\in B_j}e^{\beta(\langle w_i,z\rangle+e_i)}.
\]
Let $\mu_i$ be the soft-maximum weights and $\Sigma_\mu$ the covariance
matrix of the vectors $w_i$ under these weights. Then
$\|\nabla f\|\le1$ and $\nabla^2f=\beta\Sigma_\mu$.
If the ball contains at least $e^{k/2}$ indices in $\mathcal A_t$, their
rescaled scores are within $t'$ of the maximum over $B_j$. Therefore
\[
 \max_i\mu_i\leq e^{\beta t'-k/2}=e^{-3k/8}.
\]
For every affine $m$-plane, at most $(C_1R)^m$ of the points $w_i$ lie
within distance $c_0/(4R)$ of it. Their orthogonal projections onto the
plane lie in a ball of fixed radius and are separated by at least
$c_0/(2R)$, so a disjoint-ball packing bound applies. Take
\[
 m=\left\lfloor\frac{k}{16\log(C_1R)}\right\rfloor.
\]
Since $L\ge c\log(1/a)$, taking $C_*$ sufficiently large ensures
$m\ge1$ and $m\ge c_{P,B}C_*L/\log(1/a)$, including the integer part.
Also $k\ge cC_*\log\log(1/a)\to\infty$.
When the ball contains at least $e^{k/2}$ indices in $\mathcal A_t$,
every such tube has total weight at most
$e^{-k/4}=o(1)$. This event is impossible if the ambient dimension is at most $m$.
Otherwise, consider the affine plane through the weighted mean and
parallel to the first $m$ eigendirections of $\Sigma_\mu$. Writing the
eigenvalues of $\Sigma_\mu$ in decreasing order as $\lambda_i$, we obtain
\[
 \sum_{i>m}\lambda_i\geq c_1/R^2.
\]
Put $b=c_1/(mR^2)$, reducing $c_1$ if necessary. If $\lambda_m\geq b$, the first $m$ eigenvalues yield
\[
 \log\det(I+u\beta\Sigma_\mu)\geq m\log(1+u\beta b).
\]
If $\lambda_m<b$, the same conclusion follows from the lower bound on
$\sum_{i>m}\lambda_i$ and the chord inequality
$\log(1+u\beta x)\geq (x/b)\log(1+u\beta b)$ for $0\leq x\leq b$.
Now take $u=\sqrt L$. Then
\[
 u\beta b\geq\frac{c\log(C_1R)}{\sqrt a\,\log^P(1/a) R},
 \qquad \log(1+u\beta b)\geq \log(1/a)/3
\]
for sufficiently small $a$. Thus, on this event, the logarithm of the
determinant is at least $c_{P,B}C_*L$.

The map $z\mapsto z+u\nabla f(z)$ is injective by strong monotonicity and has positive definite derivative. Gaussian change of variables over its image gives
\[
 \mathbb E\left[\det(I+u\nabla^2f(Z))
 e^{-u\langle Z,\nabla f(Z)\rangle-u^2\|\nabla f(Z)\|^2/2}\right]\leq1.
\]
A union bound over the at most $N^2$ vectors $w_i$ obtained from all
centers gives $|\langle w_i,Z\rangle|\le C_2\sqrt L$ simultaneously for
these vectors, except on an event of probability at most $Ce^{-cL}$.
When all these bounds hold, the exponential factor is at least
$e^{-(C_2+1/2)L}$. For a fixed center, the probability that all these
bounds hold and its ball contains at least $e^{k/2}$ near-maximizers is
at most
$\exp\{(C_2+1/2-c_{P,B}C_*)L\}$.

Choose $C_*$ sufficiently large and sum over the at most $e^L$ centers.
Combining this estimate with the separated-point bound proves
\eqref{eq:newlog-nearmax}.
\end{proof}

\subsubsection{Auxiliary Gaussian variables at selected indices}

The next estimate uses the modified scores in the expected cubic term.
The modified function $g$ is strictly increasing and convex, with
$\|g''\|_\infty\le C$ and $0<g'\le C\sqrt L$. For each
log-sum-exponential term, the coordinate scores have the form
\[
 h_j(t)=v\{g(t)-g(b_j)\},\qquad 1\leq v'\leq2,
\]
where $v$ is convex with bounded derivatives through order three.
The threshold values $g(b_j)$ are still those of the original function.

\begin{lemma}[Gaussian variables selected by soft-maximum weights]
\label{lem:newlog-selected}
Suppose $a\downarrow0$ and $a\sqrt L\ge1$ eventually.
There are at most $e^{C_0L}$ coordinate indices, divided into at most two
groups. Let $Z$ be standard Gaussian in a finite-dimensional real
Hilbert space, and write $Z_j=\langle v_j,Z\rangle$, where
$\|v_j\|=1$ and $\|v_i-v_j\|\geq\sqrt{2a}$ within each group. Let
\[
 t_j=u_j+\theta d_jZ_j,\qquad 0<d_j\leq2,\quad 0<\theta\leq1,
\]
with arbitrary deterministic $u_j$.
Conditional on $Z$, sample $I$ with probabilities proportional to
$e^{\lambda h_j(t_j)}$. When a second index is needed, sample $I'$ from
the same probabilities, independently of $I$ conditional on $Z$. We may include an additional
constant score, whose derivative is zero.
Let $\eta_{ij}=\langle z_{ij},Z\rangle$ have unit variance and satisfy, for a fixed $P\geq0$,
\begin{equation}\label{eq:newlog-mark-cov}
 |\langle z_{ij},v_k\rangle|\leq a\log^P(1/a)
 \quad\hbox{for all coordinate indices }i,j,k.
\end{equation}
The variables $\eta_{ij}$ may use additional Gaussian directions and
may be mutually dependent. Put
$\rho=\frac{\log\log(1/a)}{\log(1/a)}$, and let $E$ be the event that
every selected coordinate index satisfies $d_jh'_j(t_j)\ge\rho\sqrt L$.
Then
\begin{equation}\label{eq:newlog-selected}
 \mathbb E[|\eta_{II'}|;E]\leq C_P\sqrt{\frac{L\log\log(1/a)}{\log(1/a)}}.
\end{equation}
The expectation is over both the Gaussian vector and the sampled
indices. The same bound holds for one selected index with $\eta_{II}$
in place of $\eta_{II'}$. A selection of the constant score contributes
zero. The constants are uniform over the intercepts and the $d_j$.
\end{lemma}
\begin{proof}
We use the Gumbel-max identity; see \citet[Section~2]{MTM2014}.
For fixed scores $s_j$, independent standard Gumbels $E_j$ give
\[
 \Pp\{i=\argmax_j(s_j+E_j/\lambda)\}
 =\int_{\mathbb R}
   \lambda e^{-\lambda t}e^{\lambda s_i}
   \exp\!\left\{-e^{-\lambda t}\sum_j e^{\lambda s_j}\right\}dt
 =\frac{e^{\lambda s_i}}{\sum_j e^{\lambda s_j}}.
\]
Choose the Gumbel variables independently of $Z$ and condition on
their offsets $e_j=E_j/\lambda$. The sampled index then maximizes
\[
 S_j=h_j(t_j)+e_j,
\]
including the additional constant score, if present. Write
$Y=\max_jS_j$. This function is $C\theta\sqrt L$-Lipschitz in $Z$.
For a median $y$ under the conditional Gaussian law, concentration gives
\begin{equation}\label{eq:newlog-median-event}
 |Y-y|\le C\theta\sqrt{L\log\log(1/a)}
\end{equation}
except on an event of probability at most
$C\left\{\frac{\log\log(1/a)}{\log(1/a)}\right\}^{8}$, after choosing
the constant in \eqref{eq:newlog-median-event} sufficiently large.

For each coordinate index, there is a unique $t_j^*$ satisfying
$h_j(t_j^*)+e_j=y$. Indeed, $h_j+e_j$ is strictly increasing, its range
is unbounded above, and its lower endpoint is below $y$ because
$Y\ge S_j$ almost surely. Set $b_j^*=(t_j^*-u_j)/(\theta d_j)$.

Within each sign group, define the deterministic subset
\[
 J_0=\{j:d_jg'(t_j^*)\geq\rho\sqrt L/4\}.
\]
The elementary inequality
\[
 |g'(t)^2-g'(s)^2|
 \leq2\|g''\|_\infty|g(t)-g(s)|
 \leq2\|g''\|_\infty|h_j(t)-h_j(s)|
\]
shows that, on the event in \eqref{eq:newlog-median-event}, a selected
index satisfying the slope condition belongs to $J_0$. Indeed,
$d_jh'_j\ge\rho\sqrt L$ implies $d_jg'\ge\rho\sqrt L/2$, and
\[
 \frac{\theta\log^2(1/a)}{\sqrt L\,\{\log\log(1/a)\}^{3/2}}
 \le\frac{a\log^2(1/a)}{\{\log\log(1/a)\}^{3/2}}\longrightarrow0.
\]
For $j\in J_0$ with $Z_j>b_j^*$, convexity and $S_j\leq Y$ give
\[
 \theta\rho\sqrt L\,(Z_j-b_j^*)/4\leq Y-y.
\]
For a selected index $i$ satisfying the slope condition and
$Z_i<b_i^*$, convexity at $t_i$ gives
\[
 \theta\rho\sqrt L\,(b_i^*-Z_i)\leq y-Y.
\]
Consequently, at a selected index satisfying the slope condition, the
affine score $Z_i-b_i^*$ is within
$C\frac{\log(1/a)}{\sqrt{\log\log(1/a)}}$ of the maximum over the set
$J_0$ in its sign group.
Put $P_1=\max\{P,1\}$.  The scale assumption implies
\begin{equation}\label{eq:newlog-localization-slack}
 \frac{\log(1/a)}{\sqrt{\log\log(1/a)}}\le C\log(1/a)\le Ca\log^{P_1}(1/a)\sqrt L
.
\end{equation}
All reductions are uniform in the conditioned offsets and intercepts.

For two selections, use two independent arrays of Gumbel offsets and
apply the preceding reduction to each array. Let $\xi$ be standard
normal and independent of all previous variables.
For each candidate pair $(i,j)$, define
\[
 Z^{ij}=Z+(\xi-\eta_{ij})z_{ij}.
\]
For each fixed pair $(i,j)$, $\|z_{ij}\|=1$ implies that $Z^{ij}$ is
standard Gaussian and is independent of $\eta_{ij}$. A union bound gives
$\max_{ij}|\eta_{ij}|,|\xi|\leq C\sqrt L$, except on an event of
probability at most $Ce^{-cL}$. When these bounds hold, replacing $Z$
by any $Z^{ij}$ changes every affine score in the preceding reduction
by at most $C a\log^P(1/a)\sqrt L$. Combining this estimate with
\eqref{eq:newlog-localization-slack} and the concentration events for
the two selections shows the following: if both selected indices
satisfy the slope condition, they belong to the corresponding
near-maximum sets of affine scores at $Z^{ij}$ with tolerance
\[
 t_0=C'a\log^{P_1}(1/a)\sqrt L.
\]
Those sets are contained in the corresponding sets with tolerance
$2t_0$ evaluated at $Z$. Choose $C'$ large enough for both inclusions.
Apply Lemma~\ref{lem:newlog-nearmax} at $Z$ with tolerance $2t_0$ and a
sufficiently large failure exponent. The cardinalities of the four sets
at $Z$ are then bounded by
\[
 N_0=\exp(C_P\frac{L\log\log(1/a)}{\log(1/a)})
\]
except on an event of probability at most
$C\left\{\frac{\log\log(1/a)}{\log(1/a)}\right\}^{4}$. There are two
sign groups for each of the two fixed offset arrays. The set inclusions
give the same cardinality bound at $Z^{ij}$, with tolerance $t_0$,
simultaneously for all candidate pairs.

Fix a pair of sign groups. For each candidate pair $(i,j)$, consider
the event that $|\eta_{ij}|>u$, both candidates lie in their respective
near-maximum sets at $Z^{ij}$, and both sets have size at most $N_0$.
Sum these probabilities over $(i,j)$. For each pair, independence of
$\eta_{ij}$ and $Z^{ij}$ factors out the same standard normal tail.
For fixed offsets and deterministic score parameters, the two sets are
functions of their Gaussian argument and do not depend on the auxiliary
direction $z_{ij}$. Since every $Z^{ij}$ has the same standard Gaussian
distribution as a single auxiliary vector $\widetilde Z$, each remaining
probability can be evaluated at $\widetilde Z$. Linearity of expectation
then turns their sum into the expected number of pairs in the two sets
at $\widetilde Z$, with the count set to zero unless both sets have size
at most $N_0$. This count is at most $N_0^2$. This argument uses only
the common marginal distribution of the $Z^{ij}$; it does not assume
that these vectors are equal or mutually independent.
The original sum of probabilities is therefore bounded by
\[
 \mathbb P(|\xi|>u)e^{2C_P\frac{L\log\log(1/a)}{\log(1/a)}}.
\]
Summing over the four pairs of sign groups and taking
$u=C_3\sqrt{\frac{L\log\log(1/a)}{\log(1/a)}}$ makes this bound
$O(e^{-c\frac{L\log\log(1/a)}{\log(1/a)}})$. Integrating the same tail
bound above $u$ gives an exponentially small contribution to the
expectation. For the bound involving $\eta_{II}$, use one candidate set
instead of a pair.
Finally, $\|\max_{ij}|\eta_{ij}|\|_2\le C\sqrt L$.
Cauchy--Schwarz bounds the contribution of the exceptional events in
the counting, concentration, and Gaussian tail estimates by $C\sqrt L\,\left\{\frac{\log\log(1/a)}{\log(1/a)}\right\}^{2}+o(1)$.
The part of the expectation with $|\eta_{II'}|\le u$ is at most $u$.
This proves \eqref{eq:newlog-selected}, first for fixed offsets and then
after averaging over them.
\end{proof}

\subsubsection{Bounding the contribution of tail moments}

\begin{lemma}[A logarithmic improvement for tail moments]
\label{lem:newlog-tail}
In the second replacement stage of Lemma~\ref{lem:newcrit-reserve}, let
$U_{ij}\in L^2(P)$ vanish on $\{|X_i|\le M\}$ and satisfy
$\max_{ij}\|U_{ij}\|_2\leq V$. Suppose
\[
 CM^{-(q-2)/2}\leq a\log^{p_0}(1/a)
\]
for some fixed $p_0\ge0$, and suppose the scale assumptions of
Lemma~\ref{lem:newlog-selected} hold.
For a convex log-sum-exponential term $F$ in either bracketing function,
\begin{equation}\label{eq:newlog-tail}
 \left|\mathbb E_0\sum_{ij\ell}d_id_jd_\ell F_{ij\ell}(t)
                         \mathbb E[U_{ij}X_\ell]\right|
 \leq C L^{3/2}V\{1+\lambda \left\{\frac{\log\log(1/a)}{\log(1/a)}\right\}^{1/2}\}.
\end{equation}
The constants may depend on $q,K,p_0$, the fixed bracketing width,
and the lower bound on the Gaussian coefficient. They are uniform over
$a,n,d$ and over laws $P$ satisfying the standing assumptions.
\end{lemma}
\begin{proof}
Let $W$ be isonormal on the span in $L^2(P)$ of the original coordinates
and the variables $U_{ij}$. Realize $Z_k=W(X_k)$ and
$W_{ij}=W(U_{ij})$ jointly. When $\|U_{ij}\|_2>0$, set
$\eta_{ij}=W_{ij}/\|U_{ij}\|_2$; the other terms vanish.
Cauchy--Schwarz, H\"older's inequality, and the marginal tail bound give
\[
 |\operatorname{Cov}(\eta_{ij},Z_k)|
 \leq\|X_k\mathbf1_{\{|X_i|>M\}}\|_2
 \le (\E|X_k|^q)^{1/q}\Pp(|X_i|>M)^{(q-2)/(2q)}\le CM^{-(q-2)/2}.
\]
This estimate uses only the marginal $q$th moments and imposes no independence among coordinates.

Condition on the original observations and Bernoulli indicators in the
remaining rows. On $\mathcal B$, we have $d_j\le2$, and the retained
Gaussian component gives $\theta\ge c>0$. Apply Gaussian integration
by parts under this conditioning, then average over the conditioning
variables. The expression inside the absolute value becomes
\[
 \theta^{-1}\mathbb E_0\sum_{ij}d_id_jF_{ij}(t)W_{ij}.
\]
Write $w_j$ for the Gibbs weights and $s_j=d_jh'_j(t_j)$.  Direct differentiation gives
\[
 \sum_{ij}d_id_jF_{ij}W_{ij}
 =\sum_iw_id_i^2h_i''W_{ii}
 +\lambda\sum_iw_is_i^2W_{ii}
 -\lambda\sum_{ij}w_iw_js_is_jW_{ij}.
\]
The expected absolute value of the first term is at most $CVL^{3/2}$,
since $\max_i|h_i''|\le CL$ and
$\mathbb E\max_{ij}|W_{ij}|\le CV\sqrt L$.
For the remaining terms, sample one or two indices according to the
weights $w_j$, and put $\rho=\frac{\log\log(1/a)}{\log(1/a)}$.
If a selected slope is below $\rho\sqrt L$, the product of slopes is
at most $C\rho L$; the expected absolute value, before multiplication
by $\lambda$, is therefore at most $CVL^{3/2}\rho$.
If all selected slopes are at least $\rho\sqrt L$,
Lemma~\ref{lem:newlog-selected}, applied with logarithmic exponent
$p_0$, and $|s_i|\le C\sqrt L$ give the bound
$CVL^{3/2}\left\{\frac{\log\log(1/a)}{\log(1/a)}\right\}^{1/2}$.
A selection of the constant score contributes zero.

Multiply the last two bounds by $\lambda$ and combine them with the
bound for the first term. Using
$\rho\le\left\{\frac{\log\log(1/a)}{\log(1/a)}\right\}^{1/2}$ and the
fixed lower bound on $\theta$ proves \eqref{eq:newlog-tail}.
\end{proof}

\subsubsection{Proof of the upper dimension bound}

\begin{theorem}\label{thm:shr-critical-upper}
Fix $q>3$ and $1\le K<\infty$.  Suppose $a_n\to0$ and
$a_n\log n\to\infty$.  If
\begin{equation}\label{eq:newlog-target}
 \log(2d_n)=o\!\left[
 n^{1/3}a_n^{4/[3(q-2)]}
 \left\{\frac{\log(1/a_n)}{\log\log(1/a_n)}\right\}^{
       2\min\{1,q-3\}/[3(q-2)]}
 \right],
\end{equation}
then $\Delta_{n,q,K}(d_n,a_n)\longrightarrow0$.
Here $\Delta_{n,q,K}$ is the rectangle error defined in
Section~\ref{sec:main}, with the original denominators and the population
covariance.
\end{theorem}
\begin{proof}
Suppress the subscript $n$.  In this proof put
\[
 s_*=\frac2{q-2},\qquad m_q=\frac{\min\{1,q-3\}}{q-2},
\]
and choose
\begin{equation}\label{eq:shr-analytic-L}
 L=\max\{\log(2nd),\ n^{1/3}(\log n)^{-10-4/(q-2)}\}.
\end{equation}
The second term is $o(n^{1/3}a^{4/[3(q-2)]}\left\{\frac{\log(1/a)}{\log\log(1/a)}\right\}^{2m_q/3})$
because $a\log n\to\infty$.
The same target scale dominates $\log n$.
Thus \eqref{eq:newlog-target} implies that $L$ divided by that scale tends
to zero, while $L=n^{1/3+o(1)}$ and $L=O(n^{1/3})$.
Moreover,
\[
 a\sqrt L
 \ge (a\log n)n^{1/6}(\log n)^{-6-2/(q-2)}
 \longrightarrow\infty.
\]
Lemma~\ref{lem:newlog-selected} is therefore applicable.
The dimension of the random vector is unchanged.

Fix $\varepsilon>0$ and a bracketing width $\tau$. We first consider
rectangles whose Gaussian probabilities lie in
$[\varepsilon,1-\varepsilon]$.
Take $\lambda=\Lambda/a$, with $\Lambda$ fixed and sufficiently large,
and use Lemma~\ref{lem:newcrit-reserve}.
It remains to bound the expected cubic coefficient in
\eqref{eq:newcrit-cubic}.
To truncate the first factor in $\kappa_{ijk}=\E[X_iX_jX_k]$, put
$b_i=X_i\mathbf1_{\{|X_i|\le M\}}$.
The operator defined in Remark~\ref{rem:newcrit-operator-realization}
multiplies by $b_i$ and then projects onto the span of the coordinate
variables. Its norm is at most $M$.
Lemma~\ref{lem:newcrit-twofold}, applied conditionally in
\eqref{eq:newcrit-fixed-gaussian}, bounds this contribution by
$CML^{3/2}$, uniformly over replacement steps.

If $q\ge4$, choose
\begin{equation}\label{eq:newlog-parameters-ge4}
 M=a^{-s_*}\left\{\frac{\log\log(1/a)}{\log(1/a)}\right\}^{1/(q-2)},\qquad
 U_{ij}=X_i\mathbf1_{\{|X_i|>M\}}X_j .
\end{equation}
For $q>4$, H\"older's inequality and the marginal tail bound give
\[
 \|U_{ij}\|_2^2
 \le(\E|X_i|^q)^{2/q}(\E|X_j|^q)^{2/q}
      \Pp(|X_i|>M)^{1-4/q}
 \le CM^{4-q}.
\]
For $q=4$ the same conclusion follows directly from Cauchy--Schwarz.
Consequently,
\[
 V:=\max_{i,j}\|U_{ij}\|_2\le CM^{-(q-4)/2},
 \qquad M^{-(q-2)/2}=a\left\{\frac{\log(1/a)}{\log\log(1/a)}\right\}^{1/2}.
\]
Since $\left\{\frac{\log(1/a)}{\log\log(1/a)}\right\}^{1/2}$ is bounded by a fixed power of ${\log(1/a)}$, the assumptions of
Lemma~\ref{lem:newlog-tail} hold after increasing that power to absorb
the fixed constant.  Combining the bounded and tail contributions gives the bound
\[
 CL^{3/2}\{M+\lambda V\left\{\frac{\log\log(1/a)}{\log(1/a)}\right\}^{1/2}+V\}
 \le CL^{3/2}a^{-s_*}\left\{\frac{\log\log(1/a)}{\log(1/a)}\right\}^{1/(q-2)}.
\]
Equation~\eqref{eq:newcrit-denominator-identity} bounds the sum of the
two denominator terms by $CL^{3/2}$, including when $q=4$.

For $3<q<4$, use two truncation levels:
\begin{equation}\label{eq:newlog-parameters-lt4}
\begin{split}
 M&=a^{-s_*}\left\{\frac{\log\log(1/a)}{\log(1/a)}\right\}^{(q-3)/(q-2)},\qquad
 T=a^{-s_*}\left\{\frac{\log(1/a)}{\log\log(1/a)}\right\}^{1/(q-2)},\\
 U_{ij}&=X_i\mathbf1_{\{M<|X_i|\le T\}}
                X_j\mathbf1_{\{|X_j|\le T\}}.
\end{split}
\end{equation}
Then $M<T$.  Cauchy--Schwarz and
$\E[|X_j|^4\mathbf1_{\{|X_j|\le T\}}]\le KT^{4-q}$ give
\[
 V\le CT^{(4-q)/2},\qquad
 M^{-(q-2)/2}=a\left\{\frac{\log(1/a)}{\log\log(1/a)}\right\}^{(q-3)/2}.
\]
Again Lemma~\ref{lem:newlog-tail} applies.
The part of the third-moment tensor not represented by
$\E[b_iX_jX_k]+\E[U_{ij}X_k]$ is bounded entrywise by $CT^{3-q}$.
Each omitted integrand vanishes unless $|X_i|>T$ or $|X_j|>T$.
H\"older's inequality with exponents $q,q,q,q/(q-3)$ and the marginal
tail bound give the stated estimate.
Multiplying this entrywise estimate by the absolute derivative sum in
\eqref{eq:shr-absolute-derivatives} bounds the omitted contribution by
$C\lambda^2T^{3-q}L^{3/2}$.
Substitution of \eqref{eq:newlog-parameters-lt4} gives
\[
 M+\lambda V\left\{\frac{\log\log(1/a)}{\log(1/a)}\right\}^{1/2}+\lambda^2T^{3-q}+V
 \le Ca^{-s_*}\left\{\frac{\log\log(1/a)}{\log(1/a)}\right\}^{(q-3)/(q-2)}.
\]
For the second term on the left, the exponent of
$\log(1/a)/\log\log(1/a)$ is
$(4-q)/(2(q-2))-1/2=-(q-3)/(q-2)$.
The third term has the same logarithmic exponent, and both have exponent
$-2/(q-2)$ in $a$.
The two denominator terms are bounded directly by $C\lambda L^{3/2}$,
using \eqref{eq:shr-gibbs-hessian}, the first-derivative sum, and
Lemma~\ref{lem:shr-retained}.
This is smaller than the preceding bound, since
\[
 \lambda\{a^{-s_*}\left\{\frac{\log\log(1/a)}{\log(1/a)}\right\}^{(q-3)/(q-2)}\}^{-1}
 =O\!\left(a^{(4-q)/(q-2)}\left\{\frac{\log(1/a)}{\log\log(1/a)}\right\}^{(q-3)/(q-2)}\right)\longrightarrow0.
\]

In both ranges,
\begin{equation}\label{eq:newlog-cubic}
 |\E_0\mathcal K|\le CL^{3/2}a^{-s_*}\left\{\frac{\log\log(1/a)}{\log(1/a)}\right\}^{m_q}.
\end{equation}
Multiplying by $n^{-3/2}$ and summing over the rows bounds the total
cubic contribution by
\[
 Cn^{-1/2}L^{3/2}a^{-s_*}\left\{\frac{\log\log(1/a)}{\log(1/a)}\right\}^{m_q}=o(1).
\]
The Taylor remainders from both replacement stages sum to at most
\[
 Ce^{C\lambda}(1+\lambda)^C n^{1-r/2}L^{r/2}+o(1)
 =o(1),\qquad r=q\wedge4>3.
\]
Indeed, $\lambda=o(\log n)$ and $L=n^{1/3+o(1)}$ make the explicit
leading term $n^{1-r/3+o(1)}=o(1)$.
In the final Gaussian array, every squared denominator equals $n$ and the
numerator is $\sum_iG_i$, so the resulting vector has law $N(0,R)$.
Apply the pointwise inequalities and the Gaussian expectation bound in
Lemma~\ref{lem:newcrit-brackets}. Lemma~\ref{lem:reductions} extends the
comparison to all rectangles with an additional error of at most
$\varepsilon$. First let $n\to\infty$ with $\Lambda$, $\tau$, and
$\varepsilon$ fixed. Then let $\Lambda\to\infty$, followed by
$\tau\downarrow0$ and $\varepsilon\downarrow0$.
\end{proof}

\begin{corollary}\label{cor:shr-critical-upper}
Let $q>3$, $K\ge1$, $a_n\to0$, and $a_n\log n\to\infty$.
Then
\[
 \log d_n=O\!\left(n^{1/3}a_n^{4/[3(q-2)]}\right)
 \quad\Longrightarrow\quad
 \Delta_{n,q,K}(d_n,a_n)\longrightarrow0.
\]
\end{corollary}
\begin{proof}
The exponent of $\log(1/a_n)/\log\log(1/a_n)$ in
Theorem~\ref{thm:shr-critical-upper} is positive for every $q>3$.
\end{proof}

\subsection{Lower bounds at a shrinking correlation gap}
\label{sec:spherical-gap-lower}

Fix $q>2$ and $K>1$.  For $0<a_n<1$ eventually, define
\begin{equation}\label{eq:shr-lower-scale}
 D_*=n^{1/3}a_n^{4/[3(q-2)]}
      \{\log(1/a_n)\}^{\min\{1,\,2(q-1)/[3(q-2)]\}}.
\end{equation}
For $\log(1/a_n)=o(\log n)$, the lower bound uses a spherical
construction. Smaller gaps are handled by arranging several rare
observations in one coordinate.
The dependence transition for Gaussian sample extremes studied by
\citet{HR1989} provides background for correlations approaching one.
Here we construct iid row laws, couple their sums with the Gaussian
vector having the same population covariance, and control the empirical
denominators.

\begin{theorem}\label{thm:newlower-sphere}
Suppose $a_n\to0$ and $D_*/\log n\to\infty$.
For every deterministic integer sequence $d_n\ge2$,
\begin{align*}
 \Delta_{n,q,K}(d_n,a_n)\to0
 &\quad\Longrightarrow\quad \log(2d_n)=o(D_*),\\
 \liminf_n\frac{\log(2d_n)}{D_*}>0
 &\quad\Longrightarrow\quad
       \liminf_n\Delta_{n,q,K}(d_n,a_n)>0,\\
 \frac{\log(2d_n)}{D_*}\to\infty
 &\quad\Longrightarrow\quad\Delta_{n,q,K}(d_n,a_n)\to1.
\end{align*}
The constant $K$ is fixed, and the definition of $\Delta_{n,q,K}$ also
applies for $2<q<3$.
\end{theorem}

\subsubsection{Multiple rare observations at a small gap}

\begin{lemma}\label{lem:shr-multiple-hits}
Let $q>2$, $K>1$, $a_n\in[0,1]$, and let $d_n\ge2$ be deterministic
integers.  If
\begin{equation}\label{eq:shr-multiple-hit-regime}
 a_n\log n\longrightarrow0,\qquad
 \frac{\log d_n}{\log n}\longrightarrow\infty,
\end{equation}
then $\Delta_{n,q,K}(d_n,a_n)\to1$.
\end{lemma}
\begin{proof}
Write $a=a_n$ and $d=d_n$.
Choose a fixed $\rho>0$ so small that $(1+2\rho^{1/q})^q<K$.
For sufficiently large $n$, put
\[
\begin{gathered}
 k=\left\lfloor\min\left\{
     \frac{\log d}{(q+2)\log n},\ n^{1/4}\right\}\right\rfloor,
 \qquad m=\left\lfloor e^{(q+2)k\log n}\right\rfloor,\\
 B=\sqrt n,\qquad \pi=\rho B^{-q},\qquad
 v^2=B^2\pi(1-\pi).
\end{gathered}
\]
Both quantities in the minimum defining $k$ tend to infinity.
Thus $k\to\infty$, $k\le n^{1/4}<n$, $m\le d$, and
\begin{equation}\label{eq:shr-multiple-hit-scales}
 n\pi=\rho n^{-(q-2)/2}\to0,\qquad
 v^2\sim\rho n^{-(q-2)/2},\qquad
 \frac{\log m}{k}\sim(q+2)\log n.
\end{equation}
Let $C_i$ be independent symmetric signs.
Independently, let $V_{ij}$ be independent signs with mean
$\sqrt{1-a}$, and put $Y_{ij}=C_iV_{ij}$ for $1\le j\le m$.
Each coordinate is a symmetric sign, and
$\E Y_{ij}Y_{i\ell}=1-a$ for $j\ne\ell$.
Let $I_{ij}$ be independent Bernoulli variables of parameter $\pi$,
independent of all signs, and define
\begin{equation}\label{eq:shr-multiple-hit-law}
 X_{ij}=\frac{Y_{ij}+B(I_{ij}-\pi)}{\sqrt{1+v^2}}.
\end{equation}
The rows are iid and centered, with unit marginal variances and
common off-diagonal correlation
\[
 r_0=\frac{1-a}{1+v^2}\le1-a.
\]
Minkowski's inequality and $\E|I-\pi|^q\le2\pi$ imply
\[
 (\E|X_{ij}|^q)^{1/q}\le1+B(2\pi)^{1/q}\le1+2\rho^{1/q},
\]
so the moment condition holds.

The counts $J_j=\sum_i I_{ij}$ are independent
$\operatorname{Bin}(n,\pi)$ variables.  Since
$\binom nk\ge(n/k)^k$,
\[
 p_k:=\Pp(J_j=k)\ge (n\pi/k)^k(1-\pi)^n.
\]
For large $n$, $\pi\le1/2$ and $k\le n^{1/4}$, whence
\[
 \log p_k\ge
 -k\left\{\left(\frac{q-2}{2}+\frac14\right)\log n-\log\rho\right\}
 -2n\pi.
\]
Together with \eqref{eq:shr-multiple-hit-scales}, this gives
$mp_k\to\infty$.  Therefore
\begin{equation}\label{eq:shr-multiple-hit-exists}
 \Pp\{\text{some }j\le m\text{ has }J_j=k\}
 =1-(1-p_k)^m\longrightarrow1.
\end{equation}
Choose the first coordinate $j_*$ with $J_{j_*}=k$, taking $j_*=1$ if none exists.
This choice depends only on the Bernoulli indicators and is independent
of the sign array.
Conditional on all Bernoulli indicators, the selected sign sequence
remains iid and symmetric. Hence
\[
 S_*:=\sum_iY_{ij_*}=O_{\Pp}(\sqrt n)=O_{\Pp}(B),
 \qquad
 \E[S_*^2\mid(I_{ij})_{i,j}]=n.
\]

The deterministic normalizing factor in \eqref{eq:shr-multiple-hit-law}
cancels from the self-normalized statistic.
On the event in \eqref{eq:shr-multiple-hit-exists}, the unscaled
numerator and squared denominator at $j_*$ are
\[
 N_*=S_*+B(k-n\pi),\qquad
 Q_*=n+2B\sum_iY_{ij_*}(I_{ij_*}-\pi)
       +B^2\{k(1-2\pi)+n\pi^2\}.
\]
The absolute value of the middle sum is at most $k+n\pi$.
Consequently,
\[
 \frac{N_*}{Bk}=1-\frac{n\pi}{k}+O_{\Pp}(k^{-1})\longrightarrow1,
\]
and, uniformly over the sign array and the locations of the Bernoulli successes,
\[
 \frac{Q_*}{B^2k}
 =1+\frac1k-2\pi+\frac{n\pi^2}{k}
       +O\!\left(\frac1B+\frac{n\pi}{Bk}\right)\longrightarrow1.
\]
The denominator is positive on this event for large $n$, and
\begin{equation}\label{eq:shr-multiple-hit-sample}
 \Pp\{\max_{j\le m}T_{n,j}>\tfrac12\sqrt k\}\longrightarrow1.
\end{equation}

The Gaussian comparison vector can be written as
$G_j=\sqrt{r_0}Z_0+\sqrt{1-r_0}Z_j$, where $Z_0,Z_1,\ldots,Z_m$ are
independent standard normal variables.  Since
$1-r_0=(a+v^2)/(1+v^2)\le a+v^2$,
\[
 \frac{\max_{j\le m}|G_j|}{\sqrt k}
 \le \frac{|Z_0|}{\sqrt k}
     +O_{\Pp}\!\left(\sqrt{(a+v^2)\frac{\log m}{k}}\right)
 \longrightarrow0.
\]
Here $a\log n\to0$ and $v^2\log n\to0$.
Thus $\Pp\{\max_{j\le m}G_j\le\tfrac12\sqrt k\}\to1$.
Together with \eqref{eq:shr-multiple-hit-sample}, this shows that the
event $\{x_j\le\tfrac12\sqrt{k},\ j\le m\}$ has probabilities whose
difference tends to one under the two laws.
Append $d-m$ independent symmetric-sign coordinates and ignore them
in this rectangle. Their absolute moments equal one and their
cross-correlations are zero, so both restrictions remain satisfied.
This includes $a_n=0$.
\end{proof}

\subsubsection{The spherical construction}

We first prove the assertions of Theorem~\ref{thm:newlower-sphere}
in the regime
\begin{equation}\label{eq:shr-slow-gap}
 \log(1/a)=o(\log n).
\end{equation}
The assumption \eqref{eq:shr-slow-gap} is needed only for the spherical
construction.
The required binomial expansion, uniform in the success probability,
is given in Lemma~\ref{lem:binomial-varying}.

Choose a deterministic sequence $L$ satisfying
\begin{equation}\label{eq:newlower-L}
 \liminf L/D_*>0,\qquad L\le n^{1/3}.
\end{equation}
By \eqref{eq:shr-slow-gap},
$D_*=n^{1/3+o(1)}$ and $L=n^{1/3+o(1)}$. Let
\[
 r=\lfloor L/\log(1/a)\rfloor,\qquad
 \ell=\begin{cases}r,&2<q\le4,\\L,&q>4,\end{cases}
 \qquad N=\lfloor e^\ell\rfloor.
\]
Thus $r=n^{1/3+o(1)}\to\infty$, $r\le\ell\le L$, and $\ell/r\le C\log(1/a)$.

Let $U$ be a symmetric sign and $Z$ an independent standard normal
variable. Fix $u\in(0,1/2)$ sufficiently small that
\[
 m_q:=\E|\sqrt{1-u^2}\,U+uZ|^q<K.
\]
Such a choice exists by continuity at $u=0$. Put
$\epsilon=\sqrt{2a}/u$, and let $\mathcal V\subset S^{r-1}$ be a maximal
$\epsilon$-separated set in chordal distance. Disjoint Euclidean balls of
radius $\epsilon/2$ give
\[
 |\mathcal V|\le(1+2/\epsilon)^r.
\]
Maximality also makes $\mathcal V$ an $\epsilon$-cover. There is a numerical
$c_0>0$ such that every chordal cap of radius $6\epsilon$ contains at least
$e^{c_0r}$ net points for all sufficiently large $r$ and sufficiently small
$\epsilon$. Here is a uniform verification. The cap of radius $5\epsilon$ is covered by radius-$\epsilon$ caps
whose centers lie within distance $6\epsilon$ of its center. Spherical cap area at chordal radius $s$ is proportional to
\[
 \int_0^{2\arcsin(s/2)}\sin^{r-2}t\,dt.
\]
To compare the areas of caps of radii $5\epsilon$ and $\epsilon$,
use $\sin t\ge2t/\pi_{\rm c}$ and
$2\arcsin(5\epsilon/2)\ge5\epsilon$ in the numerator, and
$\sin t\le t$ and $2\arcsin(\epsilon/2)\le1.1\epsilon$ in the denominator.
Here $\pi_{\rm c}$ denotes the circle constant. The area ratio is at least
$(2/\pi_{\rm c})^{r-2}(5/1.1)^{r-1}$, proving the assertion.
In particular,
\begin{equation}\label{eq:newlower-dim-budget}
 \log(N|\mathcal V|)\le3L
\end{equation}
eventually.

Choose $C_1$ so large that independent standard Gaussian vectors
$Z_b\in\R^r$, $1\le b\le N$, satisfy
$\max_b\|Z_b\|\le C_1\sqrt\ell$ with probability tending to one. This follows
from $\Pp(\|Z_b\|>\sqrt r+t)\le e^{-t^2/2}$ and a union bound.
Choose $c_1>0$ with $c_1^2/2<c_0/2$, and then choose a fixed $C_2$ so large
that
\begin{equation}\label{eq:newlower-C2}
 (c_1/2)\sqrt{C_2}>36C_1/u.
\end{equation}
Finally choose $\omega>0$ sufficiently small that
\[
 \E|\sqrt{1-u^2}\,U+uZ|^q+\omega^{q-2}<K.
\]
Set
\[
 \delta=C_2a^2\ell/r,\qquad s=\sqrt{1-u^2-\delta},\qquad
 B=\omega\delta^{-1/(q-2)},\qquad
 \pi(1-\pi)=\delta/B^2,\quad 0<\pi<1/2.
\]
These quantities are well-defined eventually. Indeed
$\delta\le C a^2{\log(1/a)}\to0$, $B\to\infty$, and the smaller root for $\pi$ exists.
Also
\begin{equation}\label{eq:newlower-slow}
 B,\ \pi^{-1}=n^{o(1)},\qquad B\pi\to0,\qquad
 \pi B^q\to\omega^{q-2}.
\end{equation}

For each row $i$, block $b$, and $v\in\mathcal V$, define
\begin{equation}\label{eq:newlower-law}
 X_{ibv}=sU_{ib}+u\langle v,Z_{ib}\rangle+B(I_{ibv}-\pi),
\end{equation}
where the $U$ are symmetric signs, the $Z$ are standard Gaussian vectors,
and the $I$ are Bernoulli variables of parameter $\pi$. All the underlying
variables are mutually independent, apart from the sharing explicitly
shown in \eqref{eq:newlower-law}. The rows are iid and centered. Since
$s^2+u^2+\delta=1$, each coordinate has variance one. Within a block,
distinct coordinates have covariance
\[
 s^2+u^2\langle v,w\rangle
 =1-\delta-\frac{u^2}{2}\|v-w\|^2\le1-a,
\]
and different blocks are independent.

The moment bound holds for every fixed $q>2$, even if $\E|Z|^q>K$.
Decompose $\E|X_{ibv}|^q$ according to the Bernoulli variable.
Since $B\pi\to0$, dominated convergence shows that the contribution
from $I=0$ tends to $\E|\sqrt{1-u^2}U+uZ|^q$.
Conditional on $I=1$, the $q$th absolute moment of $X_{ibv}/B$ tends to one:
the Gaussian and sign components have uniformly bounded moments of
every fixed order, whereas $B\to\infty$. Equation \eqref{eq:newlower-slow}
therefore gives
\[
 \E|X_{ibv}|^q
 \longrightarrow \E|\sqrt{1-u^2}U+uZ|^q+\omega^{q-2}<K.
\]
The denominator of each coordinate is positive almost surely, since
$u>0$ gives a continuous marginal law. The Gaussian vector with the same covariance is
\begin{equation}\label{eq:newlower-G}
 G_{bv}=sZ_b^0+u\langle v,Z_b\rangle+\sqrt\delta\,Z_{bv}^1,
\end{equation}
where the scalar variables $Z_b^0$ and $Z_{bv}^1$ are standard normal,
the vectors $Z_b$ are standard Gaussian in $\R^r$, and all these factors
are mutually independent. The factors $Z_b^0$ and $Z_b$ are shared
within each block; $Z_{bv}^1$ is specific to one coordinate.

\subsubsection{Coupling at a Gaussian maximizing coordinate}

Conditional on $Z_b$, consider net points in the cap centered at
$Z_b/\|Z_b\|$, with arbitrary center if $Z_b=0$. At least $e^{c_0r}$ such
points lie within distance $6\epsilon$. Conditional on $Z_b$, the
variables $Z_{bv}^1$ at these points remain independent standard normal.
The Gaussian Mills bound therefore shows that the probability that
none exceeds $c_1\sqrt r$ is at most $\exp\{-\exp(c_3r)\}$ for some
$c_3>0$.
The assertion holds in every block with probability tending to one:
$\log N\le L$ whereas $\exp(c_3r)\gg L$.

On $\max_b\|Z_b\|\le C_1\sqrt\ell$, the difference between $u\|Z_b\|$ and the contribution
$u\langle v,Z_b\rangle$ of a point in that cap is at most
\[
 \frac u2(6\epsilon)^2C_1\sqrt\ell
       =\frac{36C_1}{u}a\sqrt\ell.
\]
By \eqref{eq:newlower-C2}, each block therefore has a coordinate whose
Gaussian value is at least
$sZ_b^0+u\|Z_b\|+(c_1/2)\sqrt{\delta r}$.
Since the directional contribution of any coordinate is at most
$u\|Z_b\|$, every maximizer in that block, and consequently every coordinate attaining the global maximum, satisfies
\begin{equation}\label{eq:newlower-winner-private}
 Z_{bv}^1\ge(c_1/2)\sqrt r .
\end{equation}
A union bound using \eqref{eq:newlower-dim-budget} also gives
$\max_{bv}|Z_{bv}^1|\le C\sqrt L$ with probability tending to one.

For each coordinate, use the increasing quantile coupling of
$J_{bv}=\sum_i I_{ibv}\sim\operatorname{Bin}(n,\pi)$ and $Z_{bv}^1$.
Write
\[
 \sigma=\sqrt{n\pi(1-\pi)}=\sqrt{n\delta}/B .
\]
We claim that, eventually, the deterministic quantile map satisfies
\begin{equation}\label{eq:newlower-quantile}
 \frac{J(z)-n\pi}{\sigma}\ge z+\frac{z^2}{12\sigma}
 \quad\text{for every }\ (c_1/2)\sqrt r\le z\le C\sqrt L.
\end{equation}
Here $L^2/(n\pi)\to0$ by \eqref{eq:newlower-slow} and $L\le n^{1/3}$.
Apply Lemma~\ref{lem:binomial-varying} at
$x=z+z^2/(12\sigma)$. Expanding the logarithm of the Gaussian tail gives
\[
 \log\frac{\Pp\{J\ge n\pi+\sigma z+z^2/12\}}{\barPhi(z)}
   =\left(\frac1{12}+o(1)\right)\frac{z^3}{\sigma}+o(1)
\]
uniformly on the displayed interval. Its right-hand side tends uniformly
to infinity. Indeed
\begin{equation}\label{eq:newlower-strength}
 \frac{Br\sqrt\ell}{\sqrt n}
 \asymp
 \frac{a^{-2/(q-2)}r^{1+1/(q-2)}
                   \ell^{1/2-1/(q-2)}}{\sqrt n}
 \asymp (L/D_*)^{3/2},
\end{equation}
whereas
\[
 \frac{r^{3/2}}{\sigma}
 =\frac{Br\sqrt\ell}{\sqrt n}\frac{r}{\sqrt{C_2}\,a\ell},
 \qquad
 \frac{r}{a\ell}\ge\frac{c}{a\log(1/a)}\longrightarrow\infty.
\]
To verify the last equivalence in \eqref{eq:newlower-strength}, consider
the two ranges of $q$.
For $2<q\le4$, $\ell=r\asymp L/\log(1/a)$ and the total exponent of $r$ is
$3/2$, which gives the exponent one on $\log(1/a)$ in $D_*$. For $q>4$,
$\ell=L$ and $r\asymp L/\log(1/a)$, so that exponent is
$\frac23(1+\frac1{q-2})$. Thus
\eqref{eq:newlower-L} gives the required uniform divergence.
The binomial tail at the integer ceiling is strictly larger than
$\barPhi(z)$. Equivalently, its CDF immediately below that integer is
strictly smaller than $\Phi(z)$, which proves
\eqref{eq:newlower-quantile}, including at lattice thresholds.

This coupling has no accumulating approximation error over coordinates:
\eqref{eq:newlower-quantile} is a deterministic inequality for the common
quantile map on the whole indicated interval. Given a coupled count,
sample its Bernoulli array uniformly among the arrangements with that
count. This recovers the original iid Bernoulli law, independently across coordinates.

For the common signs use the binomial quantile inequality
\[
 \left|J_b^0-\frac n2-\frac{\sqrt n}{2}Z_b^0\right|
       \le1+\frac{(Z_b^0)^2}{8},
 \qquad J_b^0\sim\operatorname{Bin}(n,1/2),
\]
under the increasing quantile coupling
\citep[p.~2732]{CP2004}. Sampling each sign array uniformly conditional
on its count recovers its original iid law. A union bound for the
Gaussian tails yields
\[
 \max_b\left|n^{-1/2}\sum_iU_{ib}-Z_b^0\right|
       =O_{\Pp}(L/\sqrt n).
\]
The common Gaussian row vectors can be sampled conditionally on their sums
so that $n^{-1/2}\sum_iZ_{ib}=Z_b$ exactly. The sign, Gaussian, and
Bernoulli couplings use independent collections of underlying variables.
Together they preserve the full row law \eqref{eq:newlower-law}.

At a coordinate attaining the Gaussian maximum,
\eqref{eq:newlower-winner-private} and \eqref{eq:newlower-quantile}
show that $B(J_{bv}-n\pi)/\sqrt n$ exceeds
$\sqrt\delta Z_{bv}^1$ by at least $cBr/\sqrt n$.
The error from coupling the common signs is negligible because
$L/(Br)=O(\log(1/a)/B)\to0$. Consequently, for a constant $c>0$,
\begin{equation}\label{eq:newlower-numerator-shift}
 \Pp\left\{\max_{bv}n^{-1/2}\sum_iX_{ibv}
           \ge\max_{bv}G_{bv}+cBr/\sqrt n\right\}\longrightarrow1.
\end{equation}

\subsubsection{Controlling the empirical denominators}

For a generic coordinate of the spherical construction, write
$H=sU+u\langle v,Z\rangle$ and $Y=B(I-\pi)$. Independence gives
\[
 X^2-1=(H^2-\E H^2)+2HY+(Y^2-\delta),\qquad
 Y^2-\delta=B(1-2\pi)Y.
\]
The first term has uniformly bounded subexponential norm. For the
cross term, we use the bound, valid for every integer $j\ge2$,
\[
 \E|HY|^j\le j!C^j\delta B^{j-2},
\]
since $H$ is uniformly subgaussian, $|Y|\le B$, and $\E Y^2=\delta$.
The last term is centered, has variance at most $B^2\delta$, and is
bounded in absolute value by $B^2$. Applying the scalar Bernstein
inequality separately to these three terms, then taking a union bound
over at most $e^{3L}$ coordinates, gives
\begin{equation}\label{eq:newlower-denominator}
 \max_{bv}\left|n^{-1}\sum_iX_{ibv}^2-1\right|
  =O_{\Pp}\left((1+B\sqrt\delta)\sqrt{L/n}+B^2L/n\right)
  =o_{\Pp}(1).
\end{equation}
The maximum absolute value of the normalized numerators is
$O_{\Pp}(\sqrt L)$. For the Gaussian and sign components, this follows
from their subgaussian tails. Bernstein's inequality bounds the centered Bernoulli contribution by
$O_{\Pp}(\sqrt{\delta L}+BL/\sqrt n)=O_{\Pp}(\sqrt L)$.

By \eqref{eq:newlower-denominator}, all empirical second moments lie
between $1/2$ and $3/2$ with probability tending to one. The map
$x\mapsto x^{-1/2}$ is Lipschitz on this interval.
The sup-norm difference between the normalized sum and the
self-normalized vector, divided by
$b_n=Br/\sqrt n$, is therefore
\[
 O_{\Pp}\left(\frac{\log(1/a)}B+{\log(1/a)}\sqrt\delta+B\log(1/a)\sqrt{L/n}\right)=o_{\Pp}(1).
\]
The three limits follow from
$\delta\asymp a^2\ell/r$, $\ell/r\le C\log(1/a)$, $a \log^{3/2}(1/a)\to0$,
$B=n^{o(1)}$, and $L\le n^{1/3}$. Combining this estimate with
\eqref{eq:newlower-numerator-shift}, and decreasing a fixed constant if
necessary, yields
\begin{equation}\label{eq:newlower-final-coupling}
 \Pp\{\max T_n\ge\max G+c_4b_n\}\longrightarrow1
\end{equation}
for some $c_4>0$.

\paragraph{From the coupling to a rectangle lower bound.}

Let $F$ be the CDF of the maximum in one Gaussian block. That maximum is
a $1$-Lipschitz function of the independent underlying standard Gaussian
variables: each linear form in the maximum has a coefficient vector of norm one.
The additive term $sZ_b^0$ is independent of the other terms and has a
positive density because $s>0$. Thus $F$ is continuous and strictly
increasing. Gaussian isoperimetry implies
\begin{equation}\label{eq:newlower-isoperimetry}
 F(t+x)\ge\Phi\bigl(\Phi^{-1}(F(t))+x\bigr),\qquad x\ge0.
\end{equation}
Indeed the Euclidean $x$-enlargement of the sublevel set at $t$ lies
inside the sublevel set at $t+x$; apply the Gaussian enlargement
inequality to that set \citep[Section~6.4]{POLLARD2024}.

The maximum of the $N$ independent blocks has $\alpha$-quantile
$Q_\alpha=F^{-1}(\alpha^{1/N})$. For each fixed
$0<\eta<1/2$, \eqref{eq:newlower-isoperimetry} gives
\[
 0\le Q_{1-\eta}-Q_\eta
 \le\Phi^{-1}((1-\eta)^{1/N})-\Phi^{-1}(\eta^{1/N})
 \le C_\eta/\sqrt\ell.
\]
The final bound follows from the normal Mills ratio: both upper-tail
probabilities are asymptotic to positive constants times $N^{-1}$.

If $\liminf \sqrt\ell\,b_n>0$, partition
$[Q_{1/4},Q_{3/4}]$ into a bounded number of intervals of lengths at most
$c_4b_n$. One interval contains $\max G$ with probability bounded below by a
positive constant. At its right endpoint $t$,
\eqref{eq:newlower-final-coupling} implies
\[
 \Pp(\max G\le t)-\Pp(\max T_n\le t)
 \ge\Pp(t-c_4b_n<\max G\le t)-o(1),
\]
so the rectangle distance has a positive lower limit.
If $\sqrt\ell\,b_n\to\infty$, use $t=Q_{1-\eta}$.
Eventually $t-c_4b_n<Q_\eta$, and the discrepancy is at least
$1-2\eta-o(1)$. Letting $\eta\downarrow0$ gives convergence to one.
By \eqref{eq:newlower-strength}, the respective hypotheses are precisely
a positive liminf of $L/D_*$ and divergence of $L/D_*$.

For a general dimension sequence, along any subsequence on which
$\log(2d)/D_*$ is bounded below by a positive constant choose
\[
 L=\min\{\tfrac14\log d,\ n^{1/3}\}.
\]
Then \eqref{eq:newlower-L} holds, and \eqref{eq:newlower-dim-budget}
ensures that the construction uses at most $d$ coordinates. If
$\log(2d)/D_*\to\infty$, then $L/D_*\to\infty$ because
$D_*=o(n^{1/3})$.
Add independent symmetric-sign coordinates to reach dimension $d$ and
leave these added coordinates unrestricted in the test rectangle.
The enlarged vector satisfies the same moment and correlation assumptions.
This proves all three assertions along sequences satisfying
\eqref{eq:shr-slow-gap}; the first follows by applying the positive
lower bound on any subsequence on which the dimension ratio stays
bounded away from zero.

\subsubsection{Completion of the lower bound}

\begin{proof}[Proof of Theorem~\ref{thm:newlower-sphere}]
The spherical construction proves the assertions on every subsequence
on which $\log(1/a_n)=o(\log n)$.
On a subsequence with $\log(1/a_n)\ge\varepsilon\log n$ for some
fixed $\varepsilon>0$, we have $a_n\log n\le n^{-\varepsilon}\log n\to0$.
If $\log(2d_n)/D_*$ is bounded below by a positive constant, the
assumption $D_*/\log n\to\infty$ also gives
$\log d_n/\log n\to\infty$.
Lemma~\ref{lem:shr-multiple-hits} then gives discrepancy tending to one.

If the second assertion failed, there would be a subsequence on which
the error tends to zero while the dimension ratio stays bounded away
from zero.
If $\log(1/a_n)/\log n$ has lower limit zero there, pass to a further
subsequence on which it tends to zero and use the spherical construction.
Otherwise the ratio is bounded below by a positive constant, and the
construction in Lemma~\ref{lem:shr-multiple-hits} applies.  Both conclusions contradict
convergence of the error to zero.
For the third assertion, if convergence to one failed, take a subsequence
whose errors are at most $1-\eta$ for some $\eta>0$ and apply the same
dichotomy.  Either construction forces the error to tend to one.
Finally, if the dimension condition in the first assertion failed,
take a subsequence on which its ratio stays bounded away from zero
and apply the second assertion.
\end{proof}

\begin{remark}\label{rem:shr-lower-log-n}
Without the condition $D_*/\log n\to\infty$, the second assertion need
not hold for every dimension sequence.
If this condition fails, then $D_*\le C\log n$ on a subsequence.
Choose $c>0$ with $cC<(q-2)/2$ and put
$d_n=\lceil e^{cD_*}\rceil\ge2$.
The ratio $\log(2d_n)/D_*$ is at least $c$, whereas on this subsequence
$d_n\le2n^{cC}=o(n^{(q-2)/2})$.
For $q\ge3$, Theorem~\ref{thm:polynomial-envelope} gives error tending
to zero.  For $2<q<3$, apply Theorem~\ref{thm:core-sufficiency} with
singleton groups.
Thus the claimed positive lower limit would fail.
This restriction is needed for a lower bound that is uniform over the
stated dimension sequences. It does not characterize convergence in
every smaller dimension range.
\end{remark}
\section{Comparison lemmas and Gaussian boundary estimates}\label{sec:inputs}

This section gives the comparison and smoothing estimates used in the
preceding proofs. The final subsection combines row replacement with
Gaussian boundary estimates to prove a sufficient comparison theorem for
grouped coordinates. This is the sufficient part of the complete
criterion in Theorem~\ref{thm:legacy-grouped}; its auxiliary estimates
also enter the proofs under higher moments.

\subsection{Reduction to one-sided events and comparison inequalities}

\begin{lemma}\label{lem:reductions}
Every rectangle event $\{X\in A\}$ can be written as $\{Y\le b\}$,
where $Y$ consists of at most $2d$ signed coordinates of $X$ and the
inequality is coordinatewise. If the original coordinates form $m$
groups with within-group correlations at most $1-a$, the signed
coordinates form at most $2m$ groups with the same bound.
A uniform error bound $e_\varepsilon$ on rectangles with Gaussian
probability in $[\varepsilon,1-\varepsilon]$ extends to the bound
$e_\varepsilon+\varepsilon$ on all rectangles. These statements allow a
singular Gaussian covariance.

A lower-bound example in $d_0$ coordinates extends to any $d\ge d_0$
while preserving its rectangle discrepancy and marginal moment bound.
\end{lemma}
\begin{proof}
For $A=\prod_j[l_j,u_j]$, retain the constraints $X_j\le u_j$ and
$-X_j\le-l_j$ with finite endpoints. Assign the positive and negative copies of each group to separate
groups. Correlations within either copy are unchanged. Each coordinate of
the Gaussian reference has a continuous distribution, so the boundary of
a fixed rectangle has Gaussian measure zero even if the joint law is
singular. Regularization by an independent Gaussian vector, followed by
diagonal rescaling, therefore passes fixed-rectangle inequalities to
singular limits.

If the Gaussian probability of $A$ is below $\varepsilon$, expand its
intervals continuously until it equals $\varepsilon$ and use inclusion.
If it exceeds $1-\varepsilon$, intersect $A$ with a moving coordinate
half-space until it equals $1-\varepsilon$ and use the reverse inclusion.
The marginal continuity just noted also proves continuity in these
operations. Infinite endpoints follow by monotone limits. Finally,
append independent symmetric signs and use a cylinder rectangle. The
added coordinates have all absolute moments equal to one and correlation
zero with the original coordinates.
\end{proof}

The following bound for sums of bounded vectors follows from
\citet[Theorem~1(4)]{CCKK2026} and
\citet[Proposition~G.1]{CCKKSUP2026}.

\begin{lemma}[Bounded sums]\label{input:cckk}\label{lem:bounded-sums}
Let $W=\sum_i\xi_i$, where the $\xi_i$ are independent and centered,
$\|\xi_i\|_\infty\le b$, and $\Var(W_j)=1$. If
$Z\sim N(0,\Cov W)$, then
\begin{align}
 \distrect(W,Z)
 &\le C\min\left\{1,\ b\E\|Z\|_\infty\log(ed)
 \sqrt{\min\left(1+\log_+\frac1{b\sqrt{\log(ed)}},\log(ed)\right)}\right\}
 \label{eq:external-cckk}\\
 &\le C\min\{1,b\log^2(ed)\}.\label{eq:external-cckk-simple}
\end{align}
No restriction is imposed on the rank of $\Cov W$. If the marginal
variances lie in a fixed compact subset of $(0,\infty)$, diagonal
rescaling gives the same bounds, with constants depending on that subset.
\end{lemma}

We also use the scalar Berry--Esseen inequality and its self-normalized
version due to \citet{BG1996}. In our applications of the latter, the
scalar observations are conditionally iid, with a bounded standardized
third absolute moment. The proofs retain the required correction for
their conditional mean.
For independent centered scalar variables, \citet[Lemma~3.2, (3.5)]{LSW2011}
give
\begin{equation}\label{eq:scalar-exponential-input}
 \Pp\left\{\left|\sum_iY_i\right|
       \ge x\sqrt{\sum_iY_i^2+5\sum_i\E Y_i^2}\right\}
 \le\sqrt2 e^{-x^2/8},\qquad x>0.
\end{equation}
Nazarov's inequality, in the form proved by \citet{NAZPROOF2017},
transfers a coordinatewise coupling bound to rectangle distance:
\begin{equation}\label{eq:gaussian-coupling-transfer}
 \distrect(V,G)\le p+C\varepsilon\sqrt{\log(2d)}
 \quad\text{if}\quad
 \Pp(\|V-G\|_\infty>\varepsilon)\le p,
\end{equation}
when $G$ has unit marginal variances. Indeed, inward and outward shifts
of the finite signed thresholds bracket each rectangle event. This
form also permits singular covariance.

We use Gaussian concentration, Sudakov--Fernique comparison and Slepian's
inequality \citep{Slepian1962}, as well as Bernstein, Hoeffding and
Rosenthal inequalities \citep{Rosenthal1970}. In particular, centered
Gaussian coordinates with variance one and off-diagonal correlations
at most $1-a$ satisfy
\begin{equation}\label{eq:core-sudakov}
 \E\max_{j\le m}G_j\ge c\sqrt{a\log m},\qquad m\ge2.
\end{equation}
Comparison with $\sqrt{1-a}\,Z_0+\sqrt a\,Z_j$ proves this bound.
Concentration then implies
\begin{equation}\label{eq:core-low-count}
 \Pp(\max_{j\le m}G_j\le B)\ge\alpha>0,\quad B\ge1
 \quad\Longrightarrow\quad \log m\le C_\alpha B^2/a.
\end{equation}
For the Gaussian arguments on convex sets, we use the following
consequence of the Brascamp--Lieb variance inequality \citep{BL1976,BK2020}: a density
proportional to $e^{-V}$ with $\nabla^2V\succeq I$ has covariance at most
$I$. It is applied in Gaussian factor space or on a conditional affine
subspace. We also use Gaussian isoperimetry in the form
$\gamma(A^t)\ge\Phi(\Phi^{-1}(\gamma(A))+t)$, $t\ge0$
\citep[Section~6.4]{POLLARD2024}, and the symmetric-binomial quantile
coupling of \citet[p.~2732]{CP2004}:
\begin{equation}\label{eq:tusnady-input}
 \left|J-\frac n2-\frac{\sqrt n}{2}Z\right|
 \le1+\frac{Z^2}{8},\qquad J\sim\operatorname{Bin}(n,1/2),\quad Z\sim N(0,1).
\end{equation}

\subsection{Binomial tails}

We need a binomial tail expansion that remains uniform when the success
probability tends to zero.
The cubic correction is the usual Cram\'er correction; for
self-normalized versions see \citet{WANG2005,SG2017}.

\begin{lemma}\label{lem:binomial-tail}\label{lem:uniform-binomial}\label{lem:binomial-varying}
Let $J\sim\operatorname{Bin}(n,p)$, where $0<p\le1/2$, and put
$\sigma^2=np(1-p)$. If $x\to\infty$ and $x^4/(np)\to0$, then
\begin{equation}\label{eq:uniform-binomial-tail}
 \log\Pp(J\ge np+\sigma x)
 =\log\barPhi(x)+\frac{(1-2p)x^3}{6\sigma}+o(1).
\end{equation}
The integer threshold is rounded upward. The error is uniform over
families on which the two stated limits hold uniformly.
\end{lemma}
\begin{proof}
Set $k=\lceil np+\sigma x\rceil$, $t=k/n$ and $x'=(k-np)/\sigma$.
The assumptions imply $x/\sigma\to0$, $(t-p)/p\to0$, and
$k\wedge(n-k)\to\infty$. Stirling's formula and Taylor expansion give
\begin{align*}
 \Pp(J=k)&=\frac{1+o(1)}{\sqrt{2\pi nt(1-t)}}e^{-nD(t\Vert p)},\\
 nD(t\Vert p)
 &=\frac{(x')^2}{2}-\frac{(1-2p)(x')^3}{6\sigma}
       +O\left(\frac{(x')^4}{np}\right),
\end{align*}
where $D(t\Vert p)=t\log(t/p)+(1-t)\log((1-t)/(1-p))$.
Replacing $x'$ by $x$ changes the exponent by
$O(x/\sigma+x^2/\sigma^2)=o(1)$, and
$nt(1-t)/\sigma^2\to1$.

The successive mass ratios
\[
 \frac{\Pp(J=k+j+1)}{\Pp(J=k+j)}
 =\frac{(n-k-j)p}{(k+j+1)(1-p)}
\]
decrease in $j$. Their logarithms equal $-(1+o(1))x/\sigma$ uniformly for
$0\le j\le A\sigma/x$, for each fixed $A$. A geometric upper bound and
the sum over this initial range show that the tail-to-mass ratio is
$(1+o(1))\sigma/x$, by first letting $n\to\infty$ and then $A\to\infty$.
Now use $\barPhi(x)\sim e^{-x^2/2}/(\sqrt{2\pi}x)$. The argument applies
to every sequence satisfying the assumptions, which also proves the
uniform assertion.
\end{proof}

\begin{corollary}\label{cor:sign-tails}\label{cor:uniform-sign-tail}
For independent symmetric signs $U_i$ and $b_n=o(n^{1/4})$,
\[
 \sup_{0\le x\le b_n}
 \left|\frac{\Pp(n^{-1/2}\sum_iU_i>x)}{\barPhi(x)}-1\right|\to0.
\]
The assertion also holds with $>$ replaced by $\ge$.
\end{corollary}
\begin{proof}
Apply Lemma~\ref{lem:binomial-tail} with $p=1/2$ for divergent $x$.
The cubic term vanishes. Strict and nonstrict thresholds differ by at
most $2/\sqrt n$; the Gaussian hazard bound makes their normal tail ratio
tend to one uniformly in the stated range. On bounded intervals use the
scalar Berry--Esseen inequality. Any violating sequence has a bounded or
a divergent subsequence, so these two arguments prove the assertion.
\end{proof}
\subsection{Gaussian boundary bounds after a change of scale}

Write $\phi$ and $\Phi$ for the standard normal density and distribution
function, and $\barPhi=1-\Phi$. We shall use
\begin{equation}\label{eq:scores}
 g_h(t)=\frac{(t+\sqrt{t^2+h^2})^2}{4},\quad h>0,
 \qquad g_N(t)=-\log\barPhi(t).
\end{equation}
Both functions are strictly increasing transformations of a coordinate;
we refer to their values as scores. The smooth approximation $g_h$ to
$t_+^2$ has useful derivative bounds for self-normalized ratios. The
normal cumulative hazard $g_N$ gives a Gaussian boundary estimate with
an absolute constant when width is measured on the transformed scale.
Bounds for a common width in the original coordinates appear in
\citet[Lemma~1 and Appendix~C.1]{CCKK2026}. The coordinate widths here
depend on their thresholds. We use a weighted version of the density
decomposition and conditional monotonicity argument in
\citet[Section~5, Lemmas~5--6]{CCK2015}.

\begin{lemma}\label{lem:score-strip}\label{lem:core-score-density}
Let $G_1,\ldots,G_m$ be centered Gaussian variables of variance one,
with arbitrary joint covariance. Let $g_j$ be continuous and strictly
increasing, and suppose
\begin{equation}\label{eq:score-domination}
 C_*:=\sup_{j\le m}\sup_{x<y}
 \frac{g_N(y)-g_N(x)}{g_j(y)-g_j(x)}<\infty.
\end{equation}
For arbitrary real $c_j$ and every interval $I$,
\begin{equation}\label{eq:score-strip}
 \Pp\{\max_{j\le m}(g_j(G_j)-c_j)\in I\}\le C_*|I|.
\end{equation}
For differentiable scores with positive derivatives,
$C_*=\sup_{j,t}\phi(t)/(g_j'(t)\barPhi(t))$.
\end{lemma}
\begin{proof}
First suppose the covariance is positive definite. For $b\in\R^m$, put
\[
 p_j(t)=\Pp\{G_k-b_k\le t\ (k\ne j)\mid G_j=b_j+t\}.
\]
Gaussian regression expresses this conditional probability using
thresholds with slopes $1-R_{jk}\ge0$. Thus $p_j$ is nondecreasing, and the disjoint events specifying the
maximizing coordinate give
\begin{equation}\label{eq:weighted-gaussian-faces}
 \sum_j\barPhi(b_j)p_j(0)
 \le\sum_j\int_0^\infty\phi(b_j+t)p_j(t)\,dt\le1.
\end{equation}
If $g_j=g_N$, differentiate the orthant probability with thresholds
$b_j(u)=g_N^{-1}(c_j+u)$ with respect to $u$. The derivative is the
left side of \eqref{eq:weighted-gaussian-faces}, evaluated at $b=b(u)$,
since
$b_j'(u)=\barPhi(b_j(u))/\phi(b_j(u))$. The orthant is empty below
the first admissible level, and its probability tends to zero at that
level. Integration proves \eqref{eq:score-strip} with constant one.
Apply this result to $(G+\varepsilon Z)/\sqrt{1+\varepsilon^2}$ and let
$\varepsilon\downarrow0$ to allow singular covariance. Each transformed
coordinate is atomless, and hence so is its finite maximum.

For general scores, let $\ell_j=g_j(-\infty)$. Condition
\eqref{eq:score-domination} implies $g_j(+\infty)=+\infty$, and
\[
 K_j(v)=
 \begin{cases}
 g_N(g_j^{-1}(v)),&v>\ell_j,\\
 0,&v\le\ell_j
 \end{cases}
\]
is nondecreasing and $C_*$-Lipschitz. The second case is absent when
$\ell_j=-\infty$. Fix $s<t$, put $a_j=K_j(c_j+s)$, and let
$M=\max_j(g_j(G_j)-c_j)$. Then
\[
 \Pp(s<M\le t)
 \le\Pp\{0<\max_j(g_N(G_j)-a_j)\le C_*(t-s)\}
 \le C_*(t-s).
\]
This includes a zero threshold $a_j$, since $g_N(G_j)>0$ almost surely.
Since $M$ is atomless, the same bound holds regardless of which
endpoints belong to the interval. The derivative form follows by
integration and difference quotients.
\end{proof}

For a common score $g$, condition~\eqref{eq:score-domination} is also
necessary up to an absolute constant for a bound uniform in dimension.
If $B(g)$ is the best constant over all dimensions, Gaussian covariance
matrices, and offsets, then
\begin{equation}\label{eq:optimal-score-constant}
 \frac1e\sup_{x<y}\frac{g_N(y)-g_N(x)}{g(y)-g(x)}
 \le B(g)\le
 \sup_{x<y}\frac{g_N(y)-g_N(x)}{g(y)-g(x)}.
\end{equation}
For the lower bound, partition $[x,y]$ into $N$ intervals on each of
which the increment of $g_N$ is the same.
On each interval $[t_i,t_{i+1}]$, put $p=\barPhi(t_i)$ and
$m=\lfloor1/p\rfloor$. Then $mp(1-p)^{m-1}\ge e^{-1}$, so the
maximum of $m$ independent standard normals falls in this interval with
probability at least $e^{-1}(1-e^{-[g_N(y)-g_N(x)]/N})$.
Apply the defining bound for $B(g)$, sum, and let $N\to\infty$.

For $g_N$ the strip constant is one. For $g_h$,
\[
 g_h'(t)=\frac{(t+\sqrt{t^2+h^2})^2}{2\sqrt{t^2+h^2}}.
\]
The ratio $\phi/(g_h'\barPhi)$ is bounded for each fixed $h>0$:
it tends to $1/2$ at $+\infty$ and to zero at $-\infty$.
Thus $h=1$ or $2$ gives an absolute constant in
\eqref{eq:score-strip}. The constant need not be uniform as
$h\downarrow0$.

\subsection{A smooth maximum with weighted derivative bounds}

The smooth maximum below satisfies two derivative bounds: one in a
weighted $\ell_p$ norm and one with a logarithmic dependence on the
number of coordinates.
For an integer $p\ge4$, let $h_p$ be the distribution function with
density $c_p t^{p-1}(1-t)^{p-1}$ on $(0,1)$. For independent variables
$U_i$ with this distribution, define
\begin{equation}\label{eq:smooth-max-definition}
 F_{\delta,p,m}(z)
 =\E\left[\lambda^{-1}\log\sum_{i=1}^m
       e^{\lambda(z_i+\delta U_i/2)}\right],
 \qquad\lambda=\frac{2\log(em)}{\delta}.
\end{equation}
For $w_i\ge0$, write
$B_k(F;w)=\sum_{i_1,\ldots,i_k}|\partial_{i_1}\cdots\partial_{i_k}F|
w_{i_1}\cdots w_{i_k}$.

\begin{lemma}\label{lem:smooth-max}
The function $F=F_{\delta,p,m}$ is convex, coordinatewise nondecreasing,
translation-equivariant and $C^{p-1}$. It satisfies
\begin{equation}\label{eq:smooth-max-envelope}
 \max_i z_i\le F(z)\le\max_i z_i+\delta,
 \qquad |F(z)-F(z')|\le\|z-z'\|_\infty.
\end{equation}
For $1\le k<p$,
\begin{equation}\label{eq:smooth-max-derivatives}
 B_k(F;w)\le C_{p,k}\delta^{1-k}\|w\|_p^k,
 \qquad
 B_k(F;\mathbf1)\le C_k\delta^{1-k}\{\log(em)\}^{k-1}.
\end{equation}
Consequently,
\begin{equation}\label{eq:smooth-max-mixed}
 |D^kF[u_1,\ldots,u_k]|
 \le C_{p,k}\delta^{1-k}\prod_{j=1}^k\|u_j\|_p.
\end{equation}
\end{lemma}
\begin{proof}
The usual bounds for the logarithm of a sum of exponentials prove
\eqref{eq:smooth-max-envelope}, convexity, monotonicity, and translation
equivariance. They also give $B_1(F;\mathbf1)=1$.
Put $s=\delta/2$ and $M_s(z)=\E\max_i(z_i+sU_i)$. For $1\le a<p$,
the beta distribution function satisfies
\begin{equation}\label{eq:beta-derivatives}
 |h_p^{(a)}(t)|\le C_{p,a}\{h_p(t)(1-h_p(t))\}^{1-a/p}.
\end{equation}
Indeed, on $(0,1/2]$, $h_p(t)\asymp_p t^p$ and
$|h_p^{(a)}(t)|\le C_{p,a}t^{p-a}$; reflection handles the other endpoint.

The distribution function of the maximum defining $M_s$ at $x$ is
$P(x,z)=\prod_i r_i$, where $r_i=h_p((x-z_i)/s)$. Put
$\rho=\min_i r_i$ and $S=\sum_i(1-r_i)$, initially with $0<r_i<1$.
A Leibniz term with $\ell$ distinct differentiated factors of orders
$a_1+\cdots+a_\ell=k$ is bounded, after summing its weighted absolute
values, by
\[
 C_{p,k}s^{-k}P\rho^{-k/p}
 \prod_{b=1}^{\ell}\sum_i(1-r_i)^{1-a_b/p}w_i^{a_b}
 \le C_{p,k}s^{-k}P\rho^{-k/p}S^{\ell-k/p}\|w\|_p^k.
\]
Allowing repeated indices in the product enlarges the sum and hence
preserves this upper bound. Since
$P\le e\rho e^{-S}$, $k<p$ and $\ell\ge1$, the factor multiplying
$s^{-k}\|w\|_p^k$ is bounded. There are only finitely many multiplicity
patterns for fixed $k$. Continuity extends the bound to $r_i=0$ or $1$. The derivatives of $P$ vanish
outside $\max_i z_i\le x\le\max_i z_i+s$. Integration of the expectation
formula therefore gives
\begin{equation}\label{eq:jitter-weighted}
 B_k(M_s;w)\le C_{p,k}s^{1-k}\|w\|_p^k.
\end{equation}

For independent standard Gumbel variables $V_i$, with mean $\gamma_E$,
the Gumbel-max identity gives
\begin{equation}\label{eq:gumbel-smoothing}
 F(z)=\E M_s(z+V/\lambda)-\gamma_E/\lambda.
\end{equation}
This identity follows by multiplying the Gumbel distribution functions;
see also \citet[Section~2]{MTM2014}. The uniform derivative bound
\eqref{eq:jitter-weighted} justifies differentiating the expectation.
This proves the first bound in \eqref{eq:smooth-max-derivatives}, with
absolute values taken before summing the derivative tensor entries.

For the second part, write $\pi_i=e^{\lambda z_i}/\sum_j e^{\lambda z_j}$.
The $k$th derivative of $\lambda^{-1}\log\sum_i e^{\lambda z_i}$ is
\[
 \lambda^{k-1}\sum_{\mathcal A}
 (-1)^{|\mathcal A|-1}(|\mathcal A|-1)!
 \prod_{A\in\mathcal A}
 \left(\pi_{i_{\min A}}\prod_{b\in A}\1_{\{i_b=i_{\min A}\}}\right),
\]
where $\mathcal A$ ranges over partitions of $\{1,\ldots,k\}$.
For each partition, the product over its blocks is nonnegative and
sums to one over all indices. Summing the absolute coefficients over
partitions bounds the absolute derivative sum by $C_k\lambda^{k-1}$.
Averaging over $U$ gives the required estimate. Such derivative estimates for the soft maximum are used in
\citet{CCK2013,CCKK2022}.
If every $u_j$ is nonzero, insert
$w=\sum_{j=1}^k|u_j|/\|u_j\|_p$ into the nonnegative tensor bound.
The $k!$ mixed terms give \eqref{eq:smooth-max-mixed}. If any $u_j=0$,
the bound is immediate.
\end{proof}
\subsection{A sufficient comparison for grouped coordinates}\label{sec:qualitative-core}

The theorem below supplies sufficiency for the grouped criterion in
Theorem~\ref{thm:legacy-grouped}. We prove it here together with the
scalar and Gaussian estimates used in the main arguments. The
correlation bound is needed only for the final Gaussian comparison;
Proposition~\ref{prop:core-comparison}, which replaces the observations,
does not require it.

Fix $2<q\le3$, $K\ge1$, and $0<\eta\le1$. Under the centering,
unit-variance, and marginal-moment assumptions of \eqref{eq:model},
suppose the coordinates have a partition $C_1,\ldots,C_s$ with
$R_{jk}\le1-\eta$ for distinct indices in the same group. Write
$\mathcal E_{n,q,K}(d,m,\eta)$ for the supremum of the rectangle error
over laws admitting at most $m$ such groups.

\begin{theorem}[Sufficient comparison for grouped coordinates]\label{thm:core-sufficiency}
For deterministic $1\le m_n\le d_n$,
\begin{equation}\label{eq:grouped-qualitative}
 \frac{m_n\{\log(ed_n/m_n)\}^{q/2}}{n^{(q-2)/2}}\to0
 \quad\Longrightarrow\quad
 \mathcal E_{n,q,K}(d_n,m_n,\eta)\to0.
\end{equation}
For prescribed partitions it is enough that
$n^{-(q-2)/2}\sum_c\{\log(e|C_c|)\}^{q/2}\to0$.
\end{theorem}

We first replace the independent observations to compare expectations
of smooth functions of the self-normalized ratios and their Gaussian
counterparts.
The weighted derivatives make marginal moment bounds sufficient for
this step. We then use the within-group correlation bound to control
the difference between a Gaussian maximum and its smooth approximation.

\subsubsection{Adding one observation and bounding scalar ratios}

Joint numerator--denominator expansions are basic tools in scalar
self-normalized approximation \citep{JSW2003,SZ2016,GSS2022}. Remove
one observation from a scalar sample, and let $S$ and $Q>0$ be the sum
and sum of squares of the remaining observations. Set $t=S/\sqrt Q$.
Adding an observation $y$ gives the identity
\begin{equation}
 \frac{S+y}{\sqrt{Q+y^2}}
   =v_t(y/\sqrt Q),\qquad
 v_t(u)=\frac{t+u}{\sqrt{1+u^2}}.
 \label{eq:core-insertion}
\end{equation}
Cauchy--Schwarz and direct differentiation give, respectively,
\begin{equation}
 (v_t(u)_+)^2\leq t_+^2+1,
 \qquad v_t'(0)=1,\qquad v_t''(0)=-t.
 \label{eq:core-insertion-cap}
\end{equation}
For $t\geq0$, the first inequality follows from
$(t+u)^2\leq(t^2+1)(1+u^2)$. For $t<0$, a positive numerator is at most $u_+$,
which gives $(v_t(u)_+)^2\leq1$.

Use the score $g=g_1$ from \eqref{eq:scores}.
Direct differentiation yields
\begin{align*}
 g_h'(t)&=t+\frac{2t^2+h^2}{2\sqrt{t^2+h^2}},\\
 g_h''(t)&=1+\frac{t(2t^2+3h^2)}{2(t^2+h^2)^{3/2}},\\
 g_h'''(t)&=\frac{3h^4}{2(t^2+h^2)^{5/2}}.
\end{align*}
Consequently
\begin{equation}
 \begin{gathered}
 0\leq g_h(t)-t_+^2\leq h^2/2,\quad
 0<g_h'(t)\leq2t_++h,\quad 0\leq g_h''(t)\leq2,\\
 g_h'(t)\geq t\quad(t\geq0),\qquad
 |g_h'''(t)|\leq3/(2h),\qquad
 g_h(v_t(u))-g_h(t)\leq1+h^2/2.
 \label{eq:core-score-bounds}
 \end{gathered}
\end{equation}
Strict monotonicity preserves the threshold event even when its
threshold is negative.

\begin{lemma}[Scalar denominator and tail bounds]
\label{lem:core-scalar}
Let $q>2$ and let $Y_1,\ldots,Y_N$ be independent centered scalar variables with
$\mathbb E Y_i^2=1$ and $\mathbb E|Y_i|^q\leq K$. Put
$Q_N=\sum_iY_i^2$ and $T_N=\sum_iY_i/\sqrt{Q_N}$, with $0/0=0$.
There are $c,C>0$, depending only on $q,K$, such that
\[
 \mathbb P(Q_N<N/2)\leq e^{-cN},\qquad
 \mathbb P(|T_N|>x)\leq C e^{-x^2/C}+e^{-cN},\qquad |T_N|\leq\sqrt N.
\]
For any $r>0$ and any $k$ coordinates of a vector of these ratios,
\begin{equation}
 \mathbb E\bigl(1+\max_{j\leq k}|T_{N,j}|\bigr)^r
 \leq C_r\{\log(ek)\}^{r/2}+C_r kN^{r/2}e^{-cN}.
 \label{eq:core-ratio-moments}
\end{equation}
\end{lemma}

\begin{proof}
Set $L=(4K)^{1/(q-2)}$ and $U_i=Y_i^2\mathbf1_{\{|Y_i|\leq L\}}$.
Then $\mathbb E U_i\geq3/4$ and $U_i^2\leq L^2U_i$. At $t=(4L^2)^{-1}$,
\[
 \mathbb E e^{-tU_i}
 \leq1-t\mathbb E U_i+\tfrac12t^2\mathbb E U_i^2
 \leq e^{-21t/32}.
\]
Exponential Markov therefore gives
$\mathbb P(Q_N<N/2)\leq e^{-5tN/32}$.
The scalar self-normalized exponential inequality of \citet[Lemma~3.2, (3.5)]{LSW2011} states, for these independent centered
variables,
\[
 \mathbb P\{|\textstyle\sum_iY_i|\geq u\sqrt{Q_N+5N}\}
 \leq\sqrt2 e^{-u^2/8}.
\]
On $Q_N\geq N/2$, $Q_N+5N\leq11Q_N$. This proves the stated tail bound;
Cauchy--Schwarz gives $|T_N|\le\sqrt N$. A union bound over $k$
coordinates, followed by integration of the tail up to $\sqrt N$, proves
\eqref{eq:core-ratio-moments} without coordinate independence.
\end{proof}

\subsubsection{Weighted expansion}
For the rest of Section~\ref{sec:qualitative-core}, set $h=h_4$ and
$F_\delta=F_{\delta,4,s}$ as in \eqref{eq:smooth-max-definition}.
Whenever the number of inputs changes, use that number as the third
subscript. Thus
\begin{equation}\label{eq:core-outer-bracket}
 \max_c z_c\le F_\delta(z)\le\max_c z_c+\delta,
 \qquad |F_\delta(z)-F_\delta(z')|\le\|z-z'\|_\infty.
\end{equation}
\begin{lemma}\label{lem:core-outer}
For $k=2,3$,
$|D^kF_\delta(z)[u_1,\ldots,u_k]|
\le C\delta^{1-k}\prod_{j=1}^k\|u_j\|_3$.
Also $|DF_\delta(z)[u]|\le\|u\|_\infty$.
\end{lemma}
\begin{proof}
Use Lemma~\ref{lem:smooth-max}, $\|u\|_4\le\|u\|_3$, and monotonicity
and translation equivariance for the first derivative.
\end{proof}

The next estimate supplies the Taylor remainder for Lindeberg
replacement \citep[Theorem~1.1 and Section~2]{CHAT2006}. Within each
group the weights sum to one and depend only on the retained sample.
Conditional expectation of the weighted remainder can therefore use the
marginal moment bound on the independent incoming observation.

\begin{lemma}[Global weighted second-order expansion]
\label{lem:core-weighted}
For each family $v$, let $w_{vj}\geq0$ and $\sum_jw_{vj}=1$.
Let positive $C^3$ functions $r_{vj}$ satisfy
\[
 r_{vj}(0)=1,\qquad r_{vj}(z)\leq U\ (z\in\mathbb R),\qquad
 |r_{vj}^{(k)}(z)|\leq C_0\quad(|z|\leq1,\ k=1,2,3),
\]
where $U\geq1$. Let $\Psi:(0,\infty)^b\to[-1,1]$ have mixed
$\ell_3$ derivative norms through order three at most $J\geq1$ on
$[1/2,\infty)^b$. Define
\[
 S_v=\sum_jw_{vj}r_{vj}(z_{vj}),\quad
 A_v=\sum_jw_{vj}r_{vj}'(0)z_{vj},\quad
 C_v=\frac12\sum_jw_{vj}r_{vj}''(0)z_{vj}^2.
\]
Then, for $2<q\leq3$,
\begin{equation}
 \left|\Psi(S)-\Psi(\mathbf1)-D\Psi(\mathbf1)[A+C]
       -\tfrac12D^2\Psi(\mathbf1)[A,A]\right|
 \leq C(C_0,q)JU\sum_{v,j}w_{vj}|z_{vj}|^q.
 \label{eq:core-weighted-remainder}
\end{equation}
Coordinates may occur in more than one family.
\end{lemma}

\begin{proof}
The scalar second-order remainder is bounded by $CU|z|^q$.
For $|z|\leq1$ this follows from the third derivative and $q\leq3$;
for $|z|>1$ it follows from boundedness of $r$ and the coefficient
bounds for its quadratic Taylor polynomial. Put $Q_v=\sum_jw_{vj}|z_{vj}|^q$, $Q=\sum_vQ_v$, and
$R=S-\mathbf1-A-C$. Weighted H\"older gives
\[
 \|A\|_3\leq CQ^{1/q},\qquad
 \|C\|_3\leq CQ^{2/q},\qquad
 \|R\|_3\leq CUQ.
\]
The first two inequalities use $3/q\geq1$ and $6/q\geq1$ when summing over
families. Choose $Q_0>0$ sufficiently small. If $Q\leq Q_0$, then
$V=\mathbf1+A+C\in[3/4,5/4]^b$.
If also $S\in[1/2,\infty)^b$, integrate the first derivative between $S$ and
$V$ and expand at $\mathbf1$. The total error is at most
\[
 CJUQ+CJ(Q^{3/q}+Q^{4/q})\leq CJUQ.
\]
If some $S_v<1/2$, then $|R_v|\geq1/4$, hence $UQ\geq c$.
Boundedness of $\Psi$ and the bounds on the retained polynomial give the same
conclusion. If $Q>Q_0$, boundedness and
$Q^{1/q}+Q^{2/q}\leq C(Q_0)Q$ prove the claim.
\end{proof}

\begin{remark}
The restriction $q\leq3$ in Lemma~\ref{lem:core-weighted} is substantive.
For $q>3$, a nonzero cubic remainder cannot be bounded by $C|z|^q$ near zero.
Improved rates under higher moments therefore require a separate treatment
of the cubic contribution.
\end{remark}

\subsubsection{Comparison at fixed smoothing parameters}

For the lower orthant $A_b=\{z_j\leq b_j\text{ for all }j\}$, with
finite thresholds $b_j$, and a partition $C_1,\ldots,C_s$, fix a
smoothing parameter $\lambda\geq1$ and define
\[
 A_c(z)=\lambda^{-1}\log\sum_{j\in C_c}
                 e^{\lambda\{g(z_j)-g(b_j)\}},\qquad
 B_c(z)=(2\lambda)^{-1}\log\sum_{j\in C_c}
                 e^{2\lambda\{g(z_j)-g(b_j)\}}.
\]
Set $A_* =\max_cA_c$, $B_* =\max_cB_c$ and
$Y(z)=\max_j\{g(z_j)-g(b_j)\}$.
Comparison of the maximum with the $\ell_1$ and $\ell_2$ norms gives
\begin{equation}
 Y\leq B_*\leq A_*,\qquad 2B_*-A_*\leq Y.
 \label{eq:core-double-temperature}
\end{equation}
For the last inequality, choose a group attaining $B_*$. With
$w_j=\lambda\{g(z_j)-g(b_j)\}$ in that group,
$\sum_je^{2w_j}\leq e^{\max_jw_j}\sum_je^{w_j}$ and $A_*\geq A_c$.

For $\tau>0$ define
\[
 \widetilde A=F_\tau((A_c)_c),\qquad
 \widetilde B=F_\tau((B_c)_c),\qquad
 H^-=h(-\widetilde A/\tau),\qquad
 H^+=h((3\tau-2\widetilde B+\widetilde A)/\tau).
\]
Since $g$ is strictly increasing, $A_b=\{Y\leq0\}$.
Equations~\eqref{eq:core-outer-bracket} and
\eqref{eq:core-double-temperature} give the pointwise bounds
\begin{equation}
 H^-\leq\mathbf1_{A_b}\leq H^+.
 \label{eq:core-tests-bracket}
\end{equation}
Indeed, $Y>0$ implies $\widetilde A>0$. If $Y\leq0$, then
$2\widetilde B-\widetilde A\leq Y+2\tau\leq2\tau$, so $H^+=1$.

\begin{proposition}[Row replacement at fixed smoothing parameters]
\label{prop:core-comparison}
Let $2<q\le3$, and let $X_1,\ldots,X_n$ be independent centered rows
with common covariance $R$, unit marginal variances, and
$\max_{i,j}\mathbb E|X_{ij}|^q\leq K$.
Fix $\delta>0$, $\lambda\geq1$, and a bounded $C^3$ function
$f:\mathbb R^2\to[-1,1]$ with bounded derivatives through order three. Set
\[
 H(z)=f\bigl(F_\delta((A_c(z))_c),F_\delta((B_c(z))_c)\bigr).
\]
For $n\geq C\log(2nd)$ sufficiently large,
\begin{align}
 |\mathbb EH(T_n)-\mathbb EH(G_R)|
 &\leq C_{q,K,\delta,\lambda,f}\biggl[
 n^{-(q-2)/2}\sum_{c=1}^s\{\log(e|C_c|)\}^{q/2}
 \notag\\*
 &\hspace{13mm}+nd e^{-cn}
 +\frac{\{\log(2nd)\}^{3/2}}{\sqrt n}+n^{-2}\biggr].
 \label{eq:core-comparison}
\end{align}
The constants do not depend on thresholds, number of groups, or covariance
rank. No correlation-gap assumption is needed in this proposition.
\end{proposition}

\begin{proof}
Replace the observation vectors successively by independent
$N_d(0,R)$ vectors. At each step, after removing the observation to be
replaced, the remaining vectors are independent, have the same
covariance, and have marginal $q$th moments bounded in terms of $q,K$.
Set $N=n-1$. For coordinate $j$, let $S_j$ and $Q_j$ be the sum and
sum of squares of the retained observations, and put
$t_j=S_j/\sqrt{Q_j}$. On
$\mathcal D=\{Q_j\geq N/2\text{ for every }j\}$, adding the incoming
observation is described by \eqref{eq:core-insertion}. The Taylor
polynomial of the resulting test function, through degree two in the
incoming row $y$, is
\begin{equation}
 H(t)+\sum_jH_j(t)\frac{y_j}{\sqrt{Q_j}}
 +\frac12\sum_{j,k}
 \{H_{jk}(t)-\mathbf1_{\{j=k\}}t_jH_j(t)\}
 \frac{y_jy_k}{\sqrt{Q_jQ_k}}.
 \label{eq:core-full-polynomial}
\end{equation}
Conditional on the retained observations, this polynomial has the same
expectation for the original and Gaussian incoming vectors. Their full
covariance matrices agree, including covariances between groups.

For $j\in C_c$ put $s_c=1+\max_{k\in C_c}|t_k|$,
$z_j=s_cy_j/\sqrt{Q_j}$ and, for $\nu=1,2$,
\[
 w_{c\nu j}=
 \frac{e^{\nu\lambda\{g(t_j)-g(b_j)\}}}
      {\sum_{k\in C_c}e^{\nu\lambda\{g(t_k)-g(b_k)\}}},\qquad
 r_{c\nu j}(z)=
 e^{\nu\lambda\{g(v_{t_j}(z/s_c))-g(t_j)\}}.
\]
The weights are measurable with respect to the retained observations
and sum to one within each family. By \eqref{eq:core-score-bounds}, $r_{c\nu j}\leq e^{3\lambda}$
globally. Differentiating $v_t$ gives
\[
 v_t'(u)=\frac{1-tu}{(1+u^2)^{3/2}},\quad
 v_t''(u)=\frac{-t-3u+2tu^2}{(1+u^2)^{5/2}},\quad
 v_t'''(u)=\frac{-3+12u^2+9tu-6tu^3}{(1+u^2)^{7/2}}.
\]
For $|z|\le1$, $|v_{t_j}(z/s_c)|\le Cs_c$. Evaluated at $u=z/s_c$,
the displayed derivatives satisfy $|v_{t_j}^{\prime}|\le C$ and
$|v_{t_j}^{(2)}|+|v_{t_j}^{(3)}|\le Cs_c$. Each differentiation with
respect to $z$ introduces a factor $s_c^{-1}$. Together with
\eqref{eq:core-score-bounds}, these bounds control the first three
$z$-derivatives of $g(v_{t_j}(z/s_c))$ by absolute constants. The
corresponding derivatives of $r_{c\nu j}$ are bounded in terms of
$\lambda$ only.

Let $S_{c\nu}=\sum_jw_{c\nu j}r_{c\nu j}(z_j)$.
After adding the observation, the smooth group maxima are
$A_c(t)+\lambda^{-1}\log S_{c1}$ and
$B_c(t)+(2\lambda)^{-1}\log S_{c2}$.
As a function of all positive $S_{c\nu}$, the test has mixed $\ell_3$
derivative bounds on $[1/2,\infty)^{2s}$ independent of $s$.
This follows from Lemma~\ref{lem:core-outer} and the bounded derivatives
of $f$ and of the logarithm on $[1/2,\infty)$. H\"older's inequality
bounds the products arising from derivatives of the coordinatewise
logarithm in $\ell_3$.

Lemma~\ref{lem:core-weighted} bounds the absolute remainder after
\eqref{eq:core-full-polynomial} by
\[
 C\sum_{c,\nu,j}w_{c\nu j}s_c^q\frac{|y_j|^q}{Q_j^{q/2}}.
\]
Conditional expectation uses only
$\mathbb E(|y_j|^q\mid\hbox{retained observations})\leq C_{q,K}$.
Applying Lemma~\ref{lem:core-scalar} to $s_c$ shows that the expected
remainder from one row is at most
\[
 CN^{-q/2}\sum_c\{\log(e|C_c|)\}^{q/2}+Cd e^{-cN}.
\]
The complement of $\mathcal D$ has probability at most $d e^{-cN}$ and is
handled directly using boundedness of the tests. Summing over the $n$ replacements gives the first two terms in
\eqref{eq:core-comparison}. The resulting statistic is the vector of
self-normalized ratios formed from Gaussian observations.

To compare this statistic with $G_R$, write
$W_j=n^{-1/2}\sum_iY_{ij}$ and $V_j=n^{-1}\sum_iY_{ij}^2$.
Then $W\sim N_d(0,R)$ and $nV_j\sim\chi_n^2$ marginally.
Gaussian tails and the scalar chi-square exponential bounds, followed by
union bounds, give with probability at least $1-Cn^{-2}$,
\[
 \max_j|W_j|\leq C\sqrt{\log(2nd)},\qquad
 \max_j|V_j-1|\leq C\sqrt{\log(2nd)/n}.
\]
Taking the constant in $n\geq C\log(2nd)$ sufficiently large makes the latter
bound at most $1/2$. Equation~\eqref{eq:core-score-bounds} then gives
\[
 \max_j|g(W_j/\sqrt{V_j})-g(W_j)|
 \leq C\{\log(2nd)\}^{3/2}/\sqrt n.
\]
The smooth maxima and $F_\delta$ are one-Lipschitz in their inputs for
the sup norm, uniformly in the smoothing parameter. Composing with $f$ and using boundedness on the
exceptional event proves the remaining terms.
\end{proof}

\subsubsection{Counting Gaussian constraints near their thresholds}

Here the Gaussian argument applies to any continuous strictly increasing
score satisfying
\begin{equation}\label{eq:score-growth}
 g(y)-g(x)\ge (y^2-x^2)/2,\qquad 0\le x<y.
\end{equation}
The constants in the following Gaussian estimates are uniform over
scores satisfying \eqref{eq:score-growth}. Both $g_h$ and $g_N$ satisfy
this condition. Write $G_j=\langle v_j,Z\rangle$, $|v_j|=1$, with $Z$
standard Gaussian in factor space. The estimates
\eqref{eq:core-sudakov}--\eqref{eq:core-low-count} apply to subsets in
which the inner products between distinct directions are at most $1-a$.

Fix the Gaussian polyhedron
\[
 K_0=\{x:\langle v_j,x\rangle\leq b_j\text{ for every }j\},\qquad
 \delta_0=\mathbb P(Z\in K_0)\geq\alpha>0.
\]
For $s>0$, put $H_s(y)=\mathbb P(y+\sqrt s Z\in K_0)$ and, for $0<t<1$,
$u_t(y)=\nabla\log H_{1-t}(y)$. Define
\[
 V_{K_0}(x)=\int_0^1
 \mathbb E\left|u_t(tx+\sqrt{t(1-t)}Z)\right|^2\,dt.
\]
The function $V_{K_0}$ is nonnegative and measurable, and may take the
value $+\infty$. Let $B$ be standard Brownian motion under Wiener
measure. Then $V_{K_0}(x)$ is a version of the conditional expectation
of the time integral $\int_0^1|u_t(B_t)|^2\,dt$ given $B_1=x$.
Indeed, $B_t$ conditional on $B_1=x$ has the law of
$tx+\sqrt{t(1-t)}Z$. The field $u_t$ is the F\"ollmer drift for Brownian
motion conditioned to end in $K_0$; the conditional expectation defining
$V_{K_0}$ is taken under the original Wiener measure. The next lemma
bounds the Gaussian integral of $V_{K_0}$ over $K_0$ and uses
$V_{K_0}(x)$ to count constraints near their thresholds at a fixed
point $x$.

\begin{lemma}[Entropy identity and a count of nearby constraints]
\label{lem:core-energy}
The function $V_{K_0}$ satisfies
\begin{equation}
 \int_{K_0}V_{K_0}\,d\gamma
   =2\delta_0\log(1/\delta_0).
 \label{eq:core-energy-identity}
\end{equation}
For every $x\in K_0$ with $V_{K_0}(x)<\infty$, every $r\geq0$, and
every subset $I$ such that $\langle v_j,v_k\rangle\le1-a$ for distinct
$j,k\in I$, where $a\in(0,1]$,
\begin{equation}
 \#\{j\in I:g(b_j)-g(\langle v_j,x\rangle)\leq r\}
 \leq\exp\left\{\frac{C_\alpha}{a}
                     (1+r+V_{K_0}(x))\right\}.
 \label{eq:core-energy-count}
\end{equation}
The function $V_{K_0}$ depends on the full conditioning polyhedron
$K_0$. It is unchanged when the subset $I$ or radius $r$ varies.
\end{lemma}

\begin{proof}
Apply the F\"ollmer entropy--energy identity to the terminal density
$\mathbf1_{K_0}/\delta_0$ \citep{FOLLMER1985}, in the form of
\citet[Theorem~2]{LEHEC2013}. Under Wiener measure, with
$M_t=H_{1-t}(B_t)$, it gives
\[
 \E\int_0^1 M_t|u_t(B_t)|^2\,dt
   =2\delta_0\log(1/\delta_0).
\]
Here $M_t$ is the conditional probability of $B_1\in K_0$ given the
Brownian path up to time $t$. Since $u_t(B_t)$ is measurable with
respect to that path, the tower property permits replacing $M_t$ by
$\mathbf1_{K_0}(B_1)$ inside the expectation. Conditioning on $B_1=x$
and using the Brownian bridge distribution above, Tonelli's theorem
proves
\eqref{eq:core-energy-identity}. The finite entropy justifies this
application even though the terminal density is an indicator.

To bound the number of nearby constraints, differentiate the Gaussian
convolution to obtain
\[
 \nabla\log H_s(y)
   =s^{-1}\{\mathbb E(W\mid W\in K_0)-y\},\qquad
 W\sim N(y,sI).
\]
The conditional mean lies in the closed convex set $K_0$. Thus, for
every nonempty set $J$ of constraint indices,
\begin{equation}
 |\nabla\log H_s(y)|^2
 \geq s^{-2}\bigl(\max_{j\in J}\{\langle v_j,y\rangle-b_j\}\bigr)_+^2.
 \label{eq:core-drift-bound}
\end{equation}
Fix $x\in K_0$ with $V_{K_0}(x)<\infty$, and write $w=1+r$ and
$B_r=4\sqrt w$.
The number of counted coordinates with $b_j\leq B_r$ is at most
$\exp(C_\alpha w/a)$ by \eqref{eq:core-low-count}.
For a counted coordinate with $b_j>B_r$, put
$x_j=\langle v_j,x\rangle$ and $d_j=b_j-x_j$.
Equation~\eqref{eq:score-growth} gives
\[
 x_j>b_j/2,\qquad 0\leq d_j\leq2r/b_j.
\]
Separate these coordinates into bands
$2^\ell B_r\leq b_j<2^{\ell+1}B_r$, with cardinalities $k_\ell$.
For each nonempty band, use the time interval
$J_\ell=[4^{-\ell}/16,4^{-\ell}/8]$.
These intervals are disjoint, and throughout $s\in J_\ell$,
\[
 d_j/\sqrt s+\sqrt s\,x_j\leq5\sqrt w.
\]
In the definition of $V_{K_0}(x)$, change variables from $t$ to
$s=1-t$, and let $V_\ell$ be the integral over $s\in J_\ell$.
Evaluate \eqref{eq:core-drift-bound} at
$y=(1-s)x+\sqrt{s(1-s)}Z$. With
$M_\ell=\max_{j\text{ in band }\ell}\langle v_j,Z\rangle$, this gives
\[
 V_\ell\geq(\log2)\mathbb E(M_\ell/\sqrt2-5\sqrt w)_+^2,
 \qquad \sum_\ell V_\ell\leq V_{K_0}(x).
\]
For this inequality, the positive part permits replacement of
$\sqrt{1-s}$ by $1/\sqrt2$ even if $M_\ell$ can be negative.
Jensen's inequality and \eqref{eq:core-sudakov} imply
$\log k_\ell\leq C(w+V_\ell)/a$.
For each nonempty band, choose one coordinate and restrict its
standard normal value to $(9\sqrt w,10\sqrt w)$. This gives
$V_\ell\geq ce^{-Cw}$.
There are therefore at most $Ce^{Cw}V_{K_0}(x)$ nonempty bands.
Summing their cardinalities and using $\log(1+v)\leq1+v$, then adding the
low-threshold coordinates, proves \eqref{eq:core-energy-count}.
All steps hold for the fixed $x$ simultaneously in $r$ and $I$.
\end{proof}

\begin{lemma}[Approximation of Gaussian group maxima]
\label{lem:core-localization}
Suppose that within each group all off-diagonal correlations are at
most $1-a$, where $a\in(0,1]$, and that
$\mathbb P(Y(G_R)\leq s)\geq\alpha>0$. If
$\lambda\geq C_\alpha/a$, then for every $u>0$,
\begin{equation}
 \mathbb P\{Y(G_R)\leq s,\ A_*(G_R)>s+u\}
 \leq\frac{C_\alpha}{a\lambda u}.
 \label{eq:core-localization}
\end{equation}
More generally, let $\mathcal I$ be a collection of groups and let
$K$ be a polyhedron in Gaussian factor space, with
$\Pp(Z\in K)\ge\alpha$, such that
\[
 K\subset\{x:g(\langle v_j,x\rangle)-g(b_j)\le s:
                     j\in\textstyle\bigcup_{c\in\mathcal I}C_c\}.
\]
Then the same bound holds for
$\Pp\{Z\in K,\ \max_{c\in\mathcal I}A_c(G_R)>s+u\}$.
The maximum is over the selected groups only. Additional constraints
remain part of the conditioning event $K$.
\end{lemma}

\begin{proof}
Positive sublevel probability ensures that every $g(b_j)+s$ is in the
range of $g$. Let $b_j(s)=g^{-1}(g(b_j)+s)$ and let $K_s$ be the corresponding full
polyhedron. On $K_s$, set $D_j=g(b_j)+s-g(G_j)\geq0$.
Apply Lemma~\ref{lem:core-energy} and integrate the number of scores
whose deficit is at most $r$. For every group, this gives
\begin{align*}
 \sum_{j\in C_c}e^{-\lambda D_j}
 &=\int_0^\infty\lambda e^{-\lambda r}
           \#\{j\in C_c:D_j\leq r\}\,dr\\
 &\leq2\exp\{C_\alpha(1+V_{K_s})/a\},
\end{align*}
provided $\lambda\geq2C_\alpha/a$.
Thus $A_*>s+u$ implies $1+V_{K_s}\geq c_\alpha a\lambda u$ after adjustment
of constants, with the trivial bound used when $a\lambda u$ is small.
Equation~\eqref{eq:core-energy-identity} bounds
$\mathbb E[\mathbf1_{K_s}(1+V_{K_s})]$ by an absolute constant.
Markov's inequality proves \eqref{eq:core-localization}. For the extension, use $V_K$ for the full conditioning polyhedron $K$.
Only the selected groups enter the count, whereas
\eqref{eq:core-energy-identity} still applies to $K$.
\end{proof}

\subsubsection{Proof of the sufficient comparison}

\begin{proof}[Proof of Theorem~\ref{thm:core-sufficiency}]
First consider a lower orthant with finite thresholds and Gaussian
probability in $[\varepsilon,1-\varepsilon]$. Choose $\tau>0$ sufficiently small, depending
only on $\varepsilon$, so that Lemma~\ref{lem:score-strip} gives
$\mathbb P(Y\leq-3\tau)\geq\varepsilon/2$.

The lower test function satisfies $H^-=1$ if $Y\leq-3\tau$ and
$A_*\leq-2\tau$.
Consequently
\[
 \mathbb P(Y\leq0)-\mathbb EH^-(G_R)
 \leq C\tau+
       \frac{C_\varepsilon}{\eta\lambda\tau}.
\]
For the upper test function, $H^+>0$ implies $2B_*-A_*<4\tau$.
On $Y>5\tau$ with $2B_*-A_*\leq4\tau$, put $B'=B_*-4\tau>\tau$.
Let $r_k=\tau(3/2)^k$ and partition this event according to
$r_k<B'\leq3r_k/2$.
On the event corresponding to the $k$th interval,
\[
 Y\leq s_k:=4\tau+3r_k/2,\qquad A_*>s_k+r_k/2.
\]
Each event $\{Y(G_R)\le s_k\}$ has probability at least $\varepsilon$.
Apply Lemma~\ref{lem:core-localization} and use
$\sum_{k\geq0}r_k^{-1}=3/\tau$ to obtain
\[
 \mathbb EH^+(G_R)-\mathbb P(Y\leq0)
 \leq C\tau+
       \frac{C_\varepsilon}{\eta\lambda\tau}.
\]
Use Lemma~\ref{lem:score-strip} on $0<Y\le5\tau$ and sum the bounds
from Lemma~\ref{lem:core-localization} over the remaining geometric
levels.

For fixed $\tau,\lambda$, Proposition~\ref{prop:core-comparison} applies to
both tests. If the partition has $s\leq m$ groups, concavity of
$x\mapsto\{\log(ex)\}^{q/2}$ on $[1,\infty)$ gives
\[
 \sum_c\{\log(e|C_c|)\}^{q/2}
 \leq s\{\log(ed/s)\}^{q/2}
 \leq C_qm\{\log(ed/m)\}^{q/2}.
\]
For the last inequality put $x=s/m\leq1$ and use boundedness of
$x\{1+\log(1/x)\}^{q/2}$.
The assumed criterion also implies $m=o(n^{(q-2)/2})$ and
$\log d=o(n^{(q-2)/q})$: write $\log d=\log m+\log(d/m)$.
For a prescribed partition, use instead
\[
 s\le\sum_c\{\log(e|C_c|)\}^{q/2},\qquad
 \log\max_c|C_c|
 \le\left(\sum_c\{\log(e|C_c|)\}^{q/2}\right)^{2/q}.
\]
The assumed sum is $o(n^{(q-2)/2})$, so
$\log d\le\log s+\log\max_c|C_c|=o(n^{(q-2)/q})$ in this case as well.
Thus all terms in \eqref{eq:core-comparison} vanish. In particular,
$(q-2)/q\le1/3$ controls the error from comparing the self-normalized
Gaussian statistic with $G_R$.
Using \eqref{eq:core-tests-bracket} in the two comparison directions gives
\[
 \limsup_{n\to\infty}
 |\mathbb P(T_n\leq b)-\mathbb P(G_R\leq b)|
 \leq C\tau+
       \frac{C_\varepsilon}{\eta\lambda\tau},
\]
uniformly over the admissible laws and thresholds whose Gaussian
probabilities belong to $[\varepsilon,1-\varepsilon]$.
First let $n\to\infty$ with all smoothing constants fixed, then let
$\lambda\to\infty$, and finally let $\tau\downarrow0$.
This order of limits is part of the argument; it does not yield a quantitative
rate without further estimates on the smoothing constants.

Write a two-sided rectangle event as $\{Y\le b\}$ using selected
coordinates of $(X,-X)$. Place the constraints involving $X_j$ and
$-X_j$ in separate copies of each original group.
There are at most $2d$ coordinates and $2m$ groups, with the same
within-group gap. No gap is imposed between the positive and negative
copies. Repeated directions with
correlation one are equal almost surely under the row law, because their
difference has variance zero, and are also equal in the Gaussian reference;
retain only their tighter constraint. Replacing $(d,m)$ by $(2d,2m)$ changes
the complexity by at most a fixed factor.

Finally a bound $e_\varepsilon$ on rectangles having Gaussian probability in
$[\varepsilon,1-\varepsilon]$ extends to a bound
$e_\varepsilon+\varepsilon$ on every rectangle.
For a rectangle of Gaussian probability below $\varepsilon$, enlarge it
continuously to probability $\varepsilon$ and use inclusion.
For probability above $1-\varepsilon$, intersect it with a moving coordinate
half-space until probability $1-\varepsilon$, and use the reverse inclusion.
Continuity of each Gaussian marginal suffices even for singular joint laws.
Infinite endpoints follow by monotone limits; zero-width rectangles are
included in the two-sign argument. Letting $\varepsilon\downarrow0$ proves
the result.
\end{proof}

\section{Grouped coordinates and bounded-variation factor models}\label{sec:extensions}

We first study a correlation bound imposed within groups, with
unrestricted dependence between groups. We state the complete criterion
for $2<q\le3$ before the bounds that use higher moments and individual
group sizes, and then give their proofs. The final subsection treats
representations by independent factors with a bound on the sum of the
factor functions' total variations. For that class, we give a coupling
bound and an exact dimension threshold.

\subsection{Within-group correlation bounds}
\label{sec:higher-groups}

For fixed $q>2$, $K>1$, and a positive within-group gap, the partition
connects two settings. A bounded number of groups permits exponentially
growing dimension; one group recovers the correlation-gap
condition on all pairs of distinct coordinates. If every group is a
singleton, the correlation restriction disappears, and uniform
approximation has the polynomial condition $d=o(n^{(q-2)/2})$.
The results below give the exponential scale and the conditions for
intermediate partitions.

Fix $q>2$, $K\ge1$, and $\eta\in(0,1]$. Under the centering,
unit-variance, and marginal-moment assumptions of \eqref{eq:model},
suppose the coordinates have a partition $\mathcal C$ into $s$ nonempty
groups such that $R_{jk}\le1-\eta$ for distinct indices in the same
group. Dependence between groups is unrestricted. As in
Section~\ref{sec:qualitative-core}, write
$\mathcal E_{n,q,K}(d,m,\eta)$ for the supremum of the rectangle error
over laws admitting such a partition into at most $m$ groups.
The partition is deterministic for each law and sample size. It may
depend on that law and on $n$, but not on the realized observations.

\subsubsection{A necessary and sufficient condition for \texorpdfstring{$2<q\le3$}{2 < q <= 3}}
\label{sec:legacy}

For $2<q\le3$, the dimension and the allowed number of groups determine
an exact criterion for uniform approximation over the class.
Theorem~\ref{thm:core-sufficiency} supplies sufficiency. We prove necessity
below using the same block construction as in the higher-moment case.

\begin{theorem}\label{thm:legacy-grouped}
Fix $2<q\le3$, $K>1$, and $\eta\in(0,1]$. For arbitrary deterministic
integers $1\le m_n\le d_n$,
\[
 \mathcal E_{n,q,K}(d_n,m_n,\eta)\to0
 \quad\Longleftrightarrow\quad
 m_n\{\log(ed_n/m_n)\}^{q/2}=o(n^{(q-2)/2}).
\]
\end{theorem}

For bounded $m_n$, this criterion is equivalent to
$\log d_n=o(n^{(q-2)/q})$. For $m_n=d_n$, it is equivalent to
$d_n=o(n^{(q-2)/2})$.

\subsubsection{Individual group sizes and higher moments}

For $q\ge3$ and fixed $K>1$, Theorem~\ref{thm:higher-groups} gives the
exact condition $\log d=o(n^{1/3})$ when the number of groups is bounded.
The singleton case again gives $d=o(n^{(q-2)/2})$. To cover intermediate
partitions, we first retain the individual group sizes and then give
bounds that depend only on the dimension and the number of groups.

For ordinary bounded sums, \citet[Corollary~5]{CCKK2026} use the sum of
$\sqrt{1-|R_{jk}|}$ over pairs for which $1-|R_{jk}|$ is positive and
below a chosen threshold. Here we impose a one-sided correlation bound
within each group and leave dependence between groups unrestricted.

For $q\ge3$, define
\begin{equation}\label{eq:group-profile}
 \Gamma_n(\mathcal C)
 =n^{1-q/2}\sum_{C\in\mathcal C}\{\log(e|C|)\}^{q/2},
 \qquad \nu_q=\max\{2/3,4/[3(q-2)]\}.
\end{equation}
\begin{theorem}\label{thm:group-profile}
Let the rows be iid, centered, with unit marginal variances, covariance
$R$, and $\max_j\E|X_j|^q\le K$.
For $q=3$, $\Gamma_n(\mathcal C_n)\to0$ implies
\[
 \sup_{A\in\calR_{d_n}}|\Pp(T_n\in A)-\Pp(G_R\in A)|\to0.
\]
For $q>3$, the same conclusion follows if
\begin{equation}\label{eq:group-profile-sufficient}
 \Gamma_n(\mathcal C_n)\to0,
 \qquad \log(2nd_n)\{\log(es_n)\}^{\nu_q}=o(n^{1/3}).
\end{equation}
These conclusions are uniform over families satisfying the displayed
conditions uniformly. They are sufficient conditions for a specified
partition, not necessary conditions for an individual sequence of laws.
\end{theorem}

For a condition expressed only through the dimension and the maximum
number of groups, put
\[
 V_n=\frac{m\{\log(ed/m)\}^{q/2}}{n^{(q-2)/2}}.
\]

\begin{theorem}\label{thm:higher-groups}
Fix $q\ge3$, $K>1$, and $\eta\in(0,1]$, and let $1\le m_n\le d_n$ be
arbitrary deterministic integer sequences. If
$\mathcal E_{n,q,K}(d_n,m_n,\eta)\to0$, then
\begin{equation}\label{eq:group-necessary}
 V_n\to0,\qquad \frac{\{\log(2nd_n)\}^3}{n}\to0.
\end{equation}
For $q=3$, $V_n\to0$ is also sufficient. For $q>3$, a sufficient condition is
\begin{equation}\label{eq:group-sufficient}
 V_n\to0,\qquad
 \log(2nd_n)\{\log(em_n)\}^{\nu_q}=o(n^{1/3}).
\end{equation}
Thus, if $\log(ed_n)=O(n^{1/3-\varepsilon})$ for fixed $\varepsilon>0$,
convergence is equivalent to $V_n\to0$. For bounded $m_n$, it is equivalent
to $\log d_n=o(n^{1/3})$; for $m_n=d_n$, it is equivalent to
$d_n=o(n^{(q-2)/2})$.
\end{theorem}

\begin{remark}\label{rem:group-profile}
Using the individual group sizes can give a sharper sufficient
condition than using only $(d,m)$.
For $q=4$, take one group of size $\lfloor e^{n^{1/4}}\rfloor$ and
$\lfloor n^{4/5}\rfloor$ singleton groups. Then
$\Gamma_n=O(n^{-1/5}+n^{-1/2})\to0$ and
$\log(2nd)\{\log(es)\}^{2/3}=O(n^{1/4}(\log n)^{2/3})=o(n^{1/3})$,
whereas $V_n\asymp n^{3/10}$. The bound using the specified partition
therefore gives convergence, while the bound using only $(d,m)$ does
not. This is consistent with \cref{thm:higher-groups}, which is uniform
over all laws admitting a partition into at most $m$ groups. A covering by sets
with the same within-set correlation bound also suffices: assigning
each coordinate to one covering set cannot increase the sum in
\eqref{eq:group-profile}.
\end{remark}

\subsubsection{A block construction with rare observations}

The following block construction gives lower bounds in both moment
ranges. We retain the empirical denominators throughout its analysis.

\begin{lemma}\label{lem:marked-block-max}
Let $J\to\infty$, let $t=\Phi^{-1}(1-J^{-1})=o(n^{1/4})$, and put
$A=\sqrt n/t$, $\lambda=n\pi$, and $v=A^2\pi(1-\pi)$.
Suppose that $\lambda\to0$, $(m-1)\lambda\to\tau\in(0,\infty)$, and
$t^2v\log(em)\to0$. In each of $J$ independent blocks, take the coordinates
\begin{equation}\label{eq:group-block-law}
 X_{ib0}=U_{ib},\qquad
 X_{ibj}=\frac{U_{ib}+A(I_{ibj}-\pi)}{\sqrt{1+v}},
 \quad 1\le j<m,
\end{equation}
where the $U_{ib}$ are independent symmetric signs and the marks
$I_{ibj}$ are mutually independent Bernoulli variables with parameter
$\pi$, independent of all signs. Let $G$ be centered Gaussian with
the population covariance of this vector. Then, for every fixed
$z\in\R$,
\begin{equation}\label{eq:group-block-gaussian}
 \Pp\{t(\max_{b,j}G_{bj}-t)\le z\}\to F(z):=\exp(-e^{-z}).
\end{equation}
Take the maximum of the self-normalized ratios over the coordinates
$X_{ib0}=U_{ib}$, which contain only the sign shared within their block,
and over columns with exactly one Bernoulli mark among their $n$
observations. The distribution function of this restricted maximum at
$t+z/t$ tends to
\begin{equation}\label{eq:group-block-sample}
 \exp\{-e^{-z}[1+(1-e^{-\tau})(e^{1/2}-1)]\}.
\end{equation}
\end{lemma}

\begin{proof}
The Gaussian coordinates in one block are
$G_b$ and $(G_b+\sqrt vZ_{bj})/\sqrt{1+v}$, with independent standard
normal factors. At $x=t+z/t>0$, the block exceedance probability lies
between $\barPhi(x)$ and
$\E\barPhi(x-\sqrt v\max_{j<m}(Z_{bj})_+)$. The normal hazard bound shows that the ratio of the block exceedance
probability to $\barPhi(x)$ is at most
$\E\exp\{Ct\sqrt v\max_{j<m}(Z_{bj})_+\}$.
The Gaussian maximum tail bound and $t^2v\log(em)\to0$ make this
expectation tend to one. Mills' inequalities give $h(x)/x\to1$ for the normal hazard, so
$J\barPhi(t+z/t)=\exp\{-\int_t^{t+z/t}h(u)\,du\}\to e^{-z}$.
Independence
of the blocks proves \eqref{eq:group-block-gaussian}.

Write $W_b=n^{-1/2}\sum_iU_{ib}$. For a column with exactly one mark, at row $i_*$, the numerator
divided by $\sqrt n$ and the squared denominator divided by $n$ are
\begin{equation}\label{eq:group-one-hit-denominator}
 W_b+\frac{1-\lambda}{t},\qquad
 1+t^{-2}+\frac{2U_{i_*b}}{t\sqrt n}
       -\frac{2\pi W_b}{t}+\frac{-2\pi+n\pi^2}{t^2}.
\end{equation}
The population normalization cancels from the ratio.
Hoeffding's inequality gives $\max_b|W_b|\le Ct$ with probability tending
to one. Uniformly over columns with exactly one mark,
\eqref{eq:group-one-hit-denominator} therefore gives
\begin{equation}\label{eq:group-one-hit-ratio}
 T_{n,bj}=\frac{W_b+t^{-1}}{\sqrt{1+t^{-2}}}+o_P(t^{-1}).
\end{equation}
Indeed, after multiplication by $t$, the error is bounded in probability
by $C\{\lambda+t/\sqrt n+t^2\pi+\pi+n\pi^2\}=o(1)$.

The probability that a block contains a column with exactly one mark is
\[
 \omega_n=1-\{1-\lambda(1-\pi)^{n-1}\}^{m-1}
 \longrightarrow\omega=1-e^{-\tau}.
\]
The indicators of these block events are independent across blocks
and independent of the sums $W_b$ of the signs shared within each block.
For the restricted maximum to be at most $x=t+z/t$, a block containing
such a column imposes the smaller threshold on $W_b$,
\[
 x\sqrt{1+t^{-2}}-t^{-1}
 =t+(z-1/2)/t+O_z(t^{-3}).
\]
The symmetric-sign relative tail estimate in \cref{cor:sign-tails}
applies uniformly at these thresholds. Thus the selected maximum has
limit \eqref{eq:group-block-sample}. The uniform error in
\eqref{eq:group-one-hit-ratio} does not affect this limit, by continuity
of the limiting cdf.
\end{proof}

\subsubsection{A scalar tail ratio for \texorpdfstring{$2<q\le3$}{2 < q <= 3}}

For the proof of Theorem~\ref{thm:legacy-grouped}, fix $2<q\le3$ and
$K>1$, and choose once and for all
\[
 \varepsilon_K=(K^{1/q}-1)/8,\qquad
 0<\rho\le\min(1/128,\varepsilon_K^q).
\]
If $U$ is a symmetric sign, $I$ is an independent Bernoulli variable of
parameter $\pi$, and $A^q\pi\le2\rho$, then
\begin{equation}\label{eq:legacy-rare-law}
 X=\frac{U+A(I-\pi)}{\sqrt{1+A^2\pi(1-\pi)}}
\end{equation}
is centered, has variance one, and
\[
 \|X\|_q\le1+(2\rho)^{1/q}\le1+2\varepsilon_K<K^{1/q}.
\]
Here $\E|I-\pi|^q\le\pi$. Independent marks add only their variance to
the diagonal of a base covariance matrix. In particular, if the base
covariance is $R_0$ and $v=A^2\pi(1-\pi)$, the exact normalized covariance
is $(R_0+vI)/(1+v)$. The normalization cancels from $T_n$, but is kept in
the Gaussian reference.

\begin{lemma}\label{lem:scalar-rare-tail}
Fix $\lambda_0>0$ and put
\[
 t=(\lambda_0 n^{(q-2)/2}/\rho)^{1/q},\qquad
 A=\sqrt n/t,\qquad \pi=\lambda_0/n=\rho A^{-q}.
\]
For iid observations from \eqref{eq:legacy-rare-law},
\begin{equation}\label{eq:legacy-scalar-tail}
 \frac{\Pp(T_n>t)}{\barPhi(t)}
 \longrightarrow c_{\lambda_0}
 :=\exp\{\lambda_0(e^{1/2}-2)\}\in(0,1).
\end{equation}
\end{lemma}

\begin{proof}
Let $N_n=\sum_iI_i$ and condition on a specified set of $k$ marked rows.
Put $W=n^{-1/2}\sum_iU_i$ and $H_k=\sum_{i:I_i=1}U_i$, so $|H_k|\le k$.
The common population scaling cancels from the ratio. Without this
scaling, the numerator divided by $\sqrt n$ and the squared denominator
divided by $n$ are
\begin{equation}\label{eq:scalar-rare-denominator}
 W+(k-\lambda_0)/t,\qquad
 Q=1+\frac{k}{t^2}+\frac{2H_k}{t\sqrt n}
       -\frac{2\pi W}{t}+\frac{-2k\pi+n\pi^2}{t^2}.
\end{equation}
For fixed $k$, on $|W|\le2t$,
\[
 t^2|Q-(1+k/t^2)|
 \le C_k\{t/\sqrt n+\pi t^2+\pi+n\pi^2\}\to0.
\]
Except on $W>2t$, the event $T_n>t$ lies between the two events
\[
 W>t+(\lambda_0-k/2)/t\pm o(t^{-1}).
\]
This follows by expanding
$t\sqrt{1+k/t^2}-(k-\lambda_0)/t$; the remainder is $O_k(t^{-3})$.
For fixed $k$, $W<-2t$ cannot make the numerator positive when $n$ is
large. Since $2<q\le3$ gives $t=o(n^{1/4})$, \cref{cor:sign-tails}
applies uniformly up to $2t$. Integrating the normal hazard yields
\begin{equation}\label{eq:scalar-rare-fixed-count}
 \frac{\Pp(T_n>t\mid N_n=k)}{\barPhi(t)}
 \longrightarrow e^{k/2-\lambda_0}.
\end{equation}
The same tail estimate makes $\Pp(W>2t)/\barPhi(t)$ negligible.

To average over the number of marks, we need a bound uniform in $k$. For large $n$,
$A(1-\pi)\ge2$, and for every sign array and every marked set,
\[
 Q=\frac1n\sum_i\{U_i+A(I_i-\pi)\}^2
 \ge1-2A\pi+\frac{k}{n}
       \{A^2(1-\pi)^2-2A(1-\pi)\}
 \ge1-2A\pi.
\]
Consequently $T_n>t$ implies
\[
 W>t-k/t-C/\sqrt n,
\]
with $C=2\lambda_0$ sufficient. If $k\le t^2/2$, this threshold is at
least $t/2-C/\sqrt n$. The uniform sign tail bound and the Gaussian hazard
bound give
\[
 \frac{\Pp(T_n>t\mid N_n=k)}{\barPhi(t)}\le Ce^{2k}.
\]
If $k>t^2/2$, the same bound follows from
$\barPhi(t)^{-1}\le C(1+t)e^{t^2/2}\le Ce^{2k}$.
Thus it holds for every $0\le k\le n$. Since
$\E e^{bN_n}\le\exp\{\lambda_0(e^b-1)\}$, taking $b=3$ gives
\[
 \sup_n\E[e^{2N_n}\1_{\{N_n>M\}}]\le Ce^{-M}\to0.
\]
Poisson convergence and \eqref{eq:scalar-rare-fixed-count} now yield
\[
 \E e^{N/2-\lambda_0}
 =\exp\{\lambda_0(e^{1/2}-2)\},\qquad N\sim\operatorname{Poisson}(\lambda_0),
\]
which proves \eqref{eq:legacy-scalar-tail}.
\end{proof}

\subsubsection{Proof of the criterion for \texorpdfstring{$2<q\le3$}{2 < q <= 3}}

\begin{proof}[Proof of \cref{thm:legacy-grouped}]
Sufficiency is \cref{thm:core-sufficiency}. Suppose that, on a subsequence,
\[
 m\{\log(ed/m)\}^{q/2}\ge\varepsilon n^{(q-2)/2}
\]
for a fixed $\varepsilon>0$. After taking a further subsequence, one of
the following three cases applies.

First suppose $m\to\infty$ and $m=o(n^{(q-2)/2})$. Choose
$0<c<\min(1,\varepsilon/4)$ and set
\[
 J=\left\lfloor\exp\left\{\tfrac12
       (cn^{(q-2)/2}/(m-1))^{2/q}\right\}\right\rfloor,
 \qquad t=\Phi^{-1}(1-J^{-1}),\qquad A=\sqrt n/t,
 \qquad \pi=\rho A^{-q}.
\]
Use $J$ independent blocks of the form \eqref{eq:group-block-law}.
Then $Jm\le d$, $(m-1)n\pi\to\rho c$, and
\[
 \frac{t^3}{\sqrt n}
 =O(m^{-3/q}n^{(q-3)/q})\to0,
 \qquad t^2v\log(em)\to0.
\]
The marginal moment bound follows from \eqref{eq:legacy-rare-law}. Grouping by
position across the independent blocks gives $m$ groups with zero
within-group correlations, so any $\eta\le1$ is allowed.
The denominator calculation for a single marked observation and the
Gaussian limit in \cref{lem:marked-block-max} apply. Their limiting cdfs differ at $z=0$,
which gives a positive discrepancy for the full sample maximum.

Next suppose $m\ge cn^{(q-2)/2}$ along a further subsequence. Choose
$0<c'<c$ and use one block with $r=\lfloor c'n^{(q-2)/2}\rfloor\le m$
coordinates, including the coordinate that contains only the sign
shared within the block. Take $A=\sqrt n$ and $\pi=\rho n^{-q/2}$,
and write $W=n^{-1/2}\sum_iU_i$ for the normalized sum of those signs.
The probability that any column contains at least two marks is
$O(rn^2\pi^2)=o(1)$, and the probability of no mark in the block tends
to $e^{-\rho c'}$. At zero threshold, only the numerator sign matters.
The coordinate containing only the shared sign imposes $W\le0$ when
there is no mark; a column with one mark imposes $W\le-1+n\pi$.
The scalar central limit theorem for $W$ therefore gives the following
limit for the cdf of the sample maximum at zero:
\[
 \tfrac12e^{-\rho c'}+(1-e^{-\rho c'})\Phi(-1).
\]
The maximum of the Gaussian reference differs from the Gaussian
factor shared by the block by $o_P(1)$, since $v\log(er)\to0$.
Its cdf at zero therefore tends to $1/2$.
This again gives a positive discrepancy. Assign one coordinate to each
of at most $m$ groups.

Finally suppose $m\le M$ is bounded. Choose $\lambda_0>0$ sufficiently
small relative to $\rho,\varepsilon,M$, and use independent coordinates
with the scalar law of \cref{lem:scalar-rare-tail}. Take
$d'=\lfloor\barPhi(t)^{-1}\rfloor$. Mills' inequalities give
\[
 \log d'=\tfrac12(\lambda_0/\rho)^{2/q}n^{(q-2)/q}
              +O(\log n)\le\log d.
\]
These independent coordinates form one group satisfying the
within-group correlation bound. At the common
threshold $t$, the cdfs of the sample and Gaussian maxima tend
respectively to $e^{-c_{\lambda_0}}$ and $e^{-1}$, and these limits
differ by
\cref{lem:scalar-rare-tail}.

In all three cases, append independent signs and use a cylinder rectangle
to reach the exact target dimension. The padding can be assigned to the
existing groups and preserves their correlation restriction. This proves
necessity for arbitrary deterministic sequences.
\end{proof}

\subsubsection{Necessary conditions for \texorpdfstring{$q\ge3$}{q >= 3}}

\begin{proof}[Necessity in \cref{thm:higher-groups}]
The independent bounded-skew construction in \cref{lem:fixed-skew-lower}
proves the second condition in \eqref{eq:group-necessary}. Suppose therefore
that $\log(ed)=o(n^{1/3})$ and, along a subsequence,
$m\{\log(ed/m)\}^{q/2}\ge\varepsilon n^{(q-2)/2}$.
If $m\ge cn^{(q-2)/2}$ on a further subsequence, use the polynomial
obstruction of \cref{thm:polynomial-envelope} in at most $m$ coordinates,
place each in a different group, and pad with independent signs.
Otherwise, $m=o(n^{(q-2)/2})$ and we may assume $m\to\infty$: for bounded
$m$, $\log(ed)=o(n^{1/3})$ and $q/6\le(q-2)/2$ would imply $V_n\to0$.

Choose $0<\zeta<\varepsilon/4$ and set
\[
 J=\left\lfloor\exp\left\{\tfrac12
       (\zeta n^{(q-2)/2}/m)^{2/q}\right\}\right\rfloor,
 \quad t=\Phi^{-1}(1-J^{-1}),\quad A=\sqrt n/t,
 \quad \pi=\rho A^{-q},
\]
where $\rho>0$ is fixed and sufficiently small that
$(1+2\rho^{1/q})^q\le K$. Then
\[
 \begin{gathered}
 J\to\infty,\qquad Jm\le d,\qquad
 t\sim(\zeta n^{(q-2)/2}/m)^{1/q},\\
 (m-1)n\pi\to\rho\zeta,\qquad
 t=o(n^{1/6}),\qquad t^2v\log(em)\to0.
 \end{gathered}
\]
Use \eqref{eq:group-block-law}. Minkowski's inequality verifies the
marginal moment bound, and grouping coordinates by their positions across the
independent blocks gives $m$ groups whose distinct coordinates have
correlation zero. The two limits in \cref{lem:marked-block-max}, at $z=0$,
then give a positive rectangle discrepancy. The full sample maximum is
at least the selected maximum, so the same lower bound applies to it.
Independent sign padding yields the prescribed dimension. This proves
$V_n\to0$ is necessary.
\end{proof}

\subsubsection{Smooth comparison using the group sizes}
\label{subsec:newgroup-jitter}

Fix $q>3$. Let $H=\lceil q\rceil-1$ and $p=\lceil q\rceil+1$, and use
the scalar cutoff $h_p$ from \cref{lem:smooth-max}. To smooth the maximum over at most $2s$ group summaries, take
$F_\delta=F_{\delta,p,v}$ from that lemma, where $v\le2s$ is the number
of summaries. For each $1\le r\le H+1$ and nonnegative vector $w$,
\begin{equation}\label{eq:newgroup-two-bounds}
 \sup_z\sum_{i_1,\ldots,i_r}
 |\partial_{i_1\cdots i_r}F_\delta(z)|\prod_{k=1}^rw_{i_k}
 \le C_r\|w\|_p^r,
 \qquad
 \sup_z\sum_{i_1,\ldots,i_r}|\partial_{i_1\cdots i_r}F_\delta(z)|
 \le C_r\{\log(es)\}^{r-1}.
\end{equation}
The constants may depend on $q,\delta$, but not on $s$.
Also $\max_i z_i\le F_\delta(z)\le\max_i z_i+\delta$.
These are direct consequences of \cref{lem:smooth-max}, since $H+1<p$.
The upper and lower test functions in
\eqref{eq:core-tests-bracket}, formed using smoothing parameters
$\lambda$ and $2\lambda$, therefore remain valid. The choice $p>q$
provides the derivatives of $h_p$ and the weighted estimates needed
below.

Fix the within-group smoothing parameter $\lambda\ge1$ and width
$\delta>0$. Let $f$ be either test function in
\eqref{eq:core-tests-bracket}, with $g=g_1$ as in
\eqref{eq:core-score-bounds}. Differentiating first within groups and
then through the smooth maximum over groups gives
\begin{equation}\label{eq:newgroup-test-derivative}
 \sum_{i_1,\ldots,i_r}|f_{i_1\cdots i_r}(t)|
 \le C_r\{\log(es)\}^{r-1}(1+\|t\|_\infty)^r,
 \qquad 1\le r\le H+1.
\end{equation}
Indeed, the absolute derivative sum of the inner log-sum-exp in score
coordinates is at most $C_r(1+\lambda)^{r-1}$. A product of $b\ge1$
nonconstant outer derivative tensors of total order at most $r$
contributes at most $\{\log(es)\}^{r-b}$, which is bounded by the power in
\eqref{eq:newgroup-test-derivative}.

\subsubsection{Row replacement and the weighted remainder}
\label{subsec:newgroup-hybrid}

Use the two replacements in \cref{lem:newcrit-reserve}: original rows
$(X,X^2)$ are first replaced by the augmented rows in
\eqref{eq:newcrit-augmented}, and then by $(G,1)$. The initial statistic
is $T_n$ and the final statistic has exactly the law $N(0,R)$.
After removing row $i_0$ in either replacement, write the retained
rows as $(y_i,h_i)$, where $h_{ij}$ is the contribution to the squared
denominator in coordinate $j$. Put
\[
 Q_j=\sum_{i\ne i_0}h_{ij},\qquad
 t_j=Q_j^{-1/2}\sum_{i\ne i_0}y_{ij},\qquad
 d_j=\sqrt{n/Q_j},\qquad B=\{\min_jQ_j\ge n/4\}.
\]
By \cref{lem:shr-retained}, uniformly over all three row types,
\begin{equation}\label{eq:newgroup-retained-tail}
 \Pp(B^c)\le Cde^{-cn},\qquad
 \Pp\{B,\max_{j\in I}|t_j|>x\}\le C|I|e^{-cx^2},\qquad x\ge1,
\end{equation}
for every nonempty deterministic coordinate set $I$. Thus, for
$s_C=1+\max_{j\in C}|t_j|$ and fixed $r>0$,
\begin{equation}\label{eq:newgroup-retained-moment}
 \E[\1_Bs_C^r]\le C_r\{\log(e|C|)\}^{r/2}.
\end{equation}
Work first on
$B\cap\{\max_j|t_j|\le C_0\sqrt{\log(2nd)}\}$.
When $\log(2nd)=o(n^{1/3})$, the complement contributes at most
$O(n^{-12})$ to each replacement of a bounded test, for sufficiently
large $C_0$.

For the incoming pair $(y,b')$, interpolate between its omission and
insertion along the path
\[
 x_j(z)=\left(t_j+\frac{zd_jy_j}{\sqrt n}\right)
       \left(1+\frac{z^2d_j^2b'_j}{n}\right)^{-1/2},\qquad 0\le z\le1.
\]
For $j\in C$, let
$v_j=Cs_C(1+|y_j|+\sqrt{b'_j})/\sqrt n$.
All three row types satisfy
$\E(1+|y_j|+\sqrt{b'_j})^q\le C_{q,K}$.
When $v_j\le1$, direct differentiation and the score bounds give
\[
 \left|\frac{d^r}{dz^r}\{g(x_j(z))-g(t_j)\}\right|
 \le C_rv_j^r,\qquad 1\le r\le H+1.
\]
To check this, write $A_j=d_jy_j/\sqrt n$ and $B'_j=d_j^2b'_j/n$.
Then $|A_j|+\sqrt{B'_j}\le Cv_j/s_C$; the first path derivative is
$O(v_j/s_C)$ and higher path derivatives are
$O_r(s_C(v_j/s_C)^r)$. The chain rule gives the assertion. The same bounds hold for positive-order $z$-derivatives of
$\exp\{\lambda_\nu[g(x_j(z))-g(t_j)]\}$, where
$\lambda_\nu\in\{\lambda,2\lambda\}$, provided the score increase
$g(x_j(z))-g(t_j)$ is bounded above. The Taylor coefficients of this
exponential at $z=0$ satisfy the corresponding bounds for every $v_j$.

For an original observation, \eqref{eq:core-insertion-cap} bounds the
increase of the score. For a row with a
Gaussian component, first restrict the incoming Gaussian to
$\max_j|G_j|\le C_1\sqrt{\log(2nd)}$.
For an incoming Gaussian row of the final type $(G,1)$,
$\max_jv_j=O(\log(2nd)/\sqrt n)=o(1)$.
For an augmented row, first insert $2IX_j$ and its square, then add
$1/2$ to the squared denominator and the Gaussian term to the numerator.
Increasing the squared denominator decreases the positive part of the
ratio. More precisely,
uniformly for $0\le z\le1$,
\[
 (x_j(z)_+)^2\le
 \left(\sqrt{t_{j,+}^2+1}+C|G_j|/\sqrt n\right)^2.
\]
The score increment is consequently at most
$C\{1+s_C|G_j|/\sqrt n+G_j^2/n\}=O(1)$ when the retained ratios and
the incoming Gaussian row satisfy the displayed bounds. This estimate
is uniform in $X_j$. The Taylor remainder of the scalar exponential after degree $H$ is
therefore $O(v_j^q)$. For $v_j\le1$, use the derivative bound. For
$v_j>1$, use the upper bound on the exponential and the bounds on its
Taylor coefficients at zero.

The excluded Gaussian tail event contributes $O(n^{-12})$ to the
bounded test and to the homogeneous terms. In a coefficient of $z^r$,
count $y_j$ with degree one and $b'_j$ with degree two. For $r\le H<q$,
the product of incoming variables then has bounded $L^{q/r}$ norm.
H\"older's inequality bounds its expectation on the omitted event by
$C\Pp\{\max_j|G_j|>C_1\sqrt{\log(2nd)}\}^{1-r/q}$.
Together with \eqref{eq:newgroup-test-derivative}, the bound imposed
on the retained ratios, and the Gaussian union bound, this is negligible for sufficiently large
$C_1$. No moment of order $2H$ is needed.

Let $w_{\nu j}$ be the exponential weights determined by the retained
observations, with $\nu$ indexing each group and its two smoothing
parameters. Put
$M_\nu=\sum_jw_{\nu j}v_j^q$ and $M=\sum_\nu M_\nu$.
Weighted H\"older gives an $\ell_p$ bound $CM^{r/q}$ for the vector of
degree-$r$ coefficients, because $pr/q\ge1$, and a bound $CM$ for the
vector of remainder sums. As a function of the weighted exponential
sums, the test is a composition of $F_\delta$, coordinate logarithms, and
the fixed scalar cutoff. Its mixed derivative tensors on
$[1/2,\infty)^{2s}$ have multilinear $\ell_p$ bounds independent of $s$.
This follows from the first bound in \eqref{eq:newgroup-two-bounds},
applied to the sum of normalized absolute input vectors, and the
coordinate-product inequality
$\|v_1\cdots v_r\|_p\le\prod_k\|v_k\|_p$.

If $M$ is small, the weighted sums and their Taylor approximations
are above $1/2$. The Taylor remainder for the outer function is
$O(M^{(H+1)/q}+M)=O(M)$; terms whose total path degree exceeds $H$ have
the same bound. If $M$ is bounded below, boundedness of the test and
$CM^{r/q}$ for $r\le H<q$ give the same conclusion. Writing $P_r$ for the
homogeneous path polynomial of degree $r$, we have
\begin{equation}\label{eq:newgroup-weighted-remainder}
 \left|f(x(1))-\sum_{r=0}^HP_r\right|
 \le C\sum_{\nu,j}w_{\nu j}v_j^q.
\end{equation}
The weights are independent of the incoming row and sum to one within
each family. Taking conditional expectations uses the incoming row's
marginal $q$th moments; \eqref{eq:newgroup-retained-moment} then shows
that the remainders over both replacements sum to at most
$C\Gamma_n(\mathcal C)+o(1)$.

For all three row types, the numerator has mean zero and covariance
$R$, and the squared-denominator contribution satisfies $\E b'_j=1$.
The first replacement also matches third numerator moments and
$\E(y_ib'_j)$ by \cref{lem:newcrit-reserve}. Conditional on the retained
observations, the expectations of its terms through degree three
therefore agree for the two incoming row types. In the second replacement, the expected cubic term equals that for an
original row. For $4\le r\le H$, marginal H\"older and
\eqref{eq:newgroup-test-derivative} bound the total contribution by
\begin{equation}\label{eq:newgroup-higher-orders}
 Cn^{1-r/2}\{\log(2nd)\}^{r/2}\{\log(es)\}^{r-1}.
\end{equation}
Every incoming product has total degree $r<q$ and bounded expectation;
products of the normalized retained weights have total mass one.
The same bounds and \eqref{eq:newgroup-retained-tail} justify replacing
the event $B\cap\{\max_j|t_j|\le C_0\sqrt{\log(2nd)}\}$ by $B$
when averaging these polynomials, with an error $O(n^{-10})$ per step.

\subsubsection{The signed cubic term}
\label{subsec:newgroup-cubic}

For coordinate indices $i,j,\ell$, let
$\kappa_{ij\ell}=\E X_iX_jX_\ell$. All derivatives of $f$ below are
evaluated at the retained ratio vector $t$. At $z=0$, the path
derivatives are
\[
 x'_j=d_jy_j/\sqrt n,\qquad x''_j=-t_jd_j^2b'_j/n,
 \qquad x'''_j=-3d_j^3y_jb'_j/n^{3/2}.
\]
Conditional on the retained observations, the expected cubic term
vanishes for the Gaussian row and agrees for the augmented and
original rows. Thus the expected cubic coefficient
in the second replacement, multiplied by $6n^{3/2}$, is
\begin{equation}\label{eq:newgroup-cubic-polynomial}
 \sum_{ij\ell}d_id_jd_\ell f_{ij\ell}\kappa_{ij\ell}
 -3\sum_{ij}d_id_j^2t_jf_{ij}\kappa_{ijj}
 -3\sum_jd_j^3f_j\kappa_{jjj}.
\end{equation}
After multiplication by $\1_B$, the last two sums have $L^1$ norms
at most $C\log(es)\{\log(2nd)\}^{3/2}$ and
$C\sqrt{\log(2nd)}$, respectively.
For the first sum, set
\[
 U_{ij}=\begin{cases}
 X_iX_j,&q\ge4,\\
 X_iX_j\1_{\{|X_i|\le M,\ |X_j|\le M\}},&3<q<4,
 \end{cases}
\]
and take $(v_M,e_M)=(1,0)$ for $q\ge4$ and
$(v_M,e_M)=(M^{(4-q)/2},M^{3-q})$ for $3<q<4$. Then
\begin{equation}\label{eq:newgroup-tensor-truncation}
 \|U_{ij}\|_2\le Cv_M,
 \qquad |\kappa_{ij\ell}-\E(U_{ij}X_\ell)|\le Ce_M.
\end{equation}
Cauchy--Schwarz reduces the first bound to fourth moments, truncated
at $M$ when $3<q<4$. For the second bound when $3<q<4$, apply
H\"older to the three
coordinate factors and the indicator of
$\{|X_i|>M\}\cup\{|X_j|>M\}$, whose probability is at most $2KM^{-q}$.

In an isonormal realization, write $Z_i=W(X_i)$ and $W_{ij}=W(U_{ij})$.
The Gaussian union bound gives
$\|\max_{ij}|W_{ij}|\|_2\le Cv_M\sqrt{\log(2nd)}$.
Condition on the retained original components and indicators in the
second replacement. By \cref{lem:newcrit-reserve},
$t_j=t_{0j}+\theta d_jZ_j$ with $\theta\ge1/2$, while $B$ and the $d_j$
are fixed. The Gaussian integration-by-parts identity
\eqref{eq:newcrit-one-ibp} gives
\begin{equation}\label{eq:newgroup-ibp}
 \E\left[\1_B\sum_{ij\ell}
     d_id_jd_\ell f_{ij\ell}\E(U_{ij}X_\ell)\right]
 =\theta^{-1}\E\left[\1_B\sum_{ij}d_id_jf_{ij}W_{ij}\right].
\end{equation}
Polynomial growth suffices for this integration-by-parts identity.
The conditioning fixes $B$ and the denominators, so they are not
differentiated. Cauchy--Schwarz,
\eqref{eq:newgroup-test-derivative}, and
\eqref{eq:newgroup-retained-moment} bound the right side by
$C\log(es)v_M\{\log(2nd)\}^{3/2}$.
The error from truncating the third-moment tensor is bounded by
$C\{\log(es)\}^2e_M\{\log(2nd)\}^{3/2}$.

For $3<q<4$, choose $M=\{\log(es)\}^{2/(q-2)}$. Then both
$\log(es)v_M$ and $\{\log(es)\}^2e_M$ equal
$\{\log(es)\}^{2/(q-2)}$. For $q\ge4$, no truncation is needed.
Thus, with $b_q=\max\{1,2/(q-2)\}$, the cubic contributions summed
over the second replacement are bounded by
\begin{equation}\label{eq:newgroup-cubic-error}
 Cn^{-1/2}\{\log(2nd)\}^{3/2}\{\log(es)\}^{b_q}.
\end{equation}
Third-moment matching cancels the expected cubic contribution in
the first replacement. In the second replacement, the coefficient of the
independent Gaussian component stays bounded below.

\subsubsection{Completion of the comparison}

\begin{proof}[Proof of \cref{thm:group-profile} and sufficiency in
\cref{thm:higher-groups}]
For $q>3$, the preceding estimates bound the difference between the
expectations of either fixed test function by
\[
 C\left[\Gamma_n(\mathcal C)
 +n^{-1/2}\{\log(2nd)\}^{3/2}\{\log(es)\}^{b_q}
 +\sum_{r=4}^H n^{1-r/2}\{\log(2nd)\}^{r/2}\{\log(es)\}^{r-1}\right]
 +o(1).
\]
Each term $\{\log(e|C|)\}^{q/2}$ is at least one. Hence
$\Gamma_n\to0$ implies $s=o(n^{(q-2)/2})$, so $\log(es)=O(\log n)$.
The second condition in \eqref{eq:group-profile-sufficient} is exactly
what makes the cubic bound vanish, because $\nu_q=2b_q/3$.
It also gives $\log(2nd)=o(n^{1/3})$. For $r\ge4$, the remaining terms
are $O(n^{1-r/3}(\log n)^{r-1})=o(1)$.

For $q=3$, retain the sum over individual group sizes in
\cref{prop:core-comparison}. If $M_0=\max_{C\in\mathcal C}|C|$, then
\[
 \log d\le\log s+\log M_0
 \le O(\log n)+(\sqrt n\,\Gamma_n)^{2/3}=o(n^{1/3}),
\]
so the remaining terms in \eqref{eq:core-comparison} vanish.

For a lower orthant with Gaussian probability in
$[\varepsilon,1-\varepsilon]$, apply \cref{lem:core-localization} and
the pointwise bounds \eqref{eq:core-tests-bracket}. The resulting error
is at most $C_\varepsilon\delta+C_\varepsilon/(\eta\lambda\delta)+o(1)$,
independently of $s$. Let $n\to\infty$, then $\lambda\to\infty$, and
finally $\delta\downarrow0$. A rectangle uses at most two signed copies
of each group, increasing \eqref{eq:group-profile} by at most a factor two.
The logarithms change by bounded additive constants. Gaussian marginal
continuity extends the bound from rectangles with Gaussian probability
in $[\varepsilon,1-\varepsilon]$ to all rectangles, including infinite
endpoints and zero widths, at an additional cost at most $\varepsilon$.
Then let $\varepsilon\downarrow0$. This proves \cref{thm:group-profile}.

To obtain upper bounds uniform over the class, use concavity of
$x\mapsto(q/2+\log x)^{q/2}$ on $[1,\infty)$ and monotonicity of
$s\mapsto s\{q/2+\log(d/s)\}^{q/2}$ for $1\le s\le d$. They give
\[
 \sum_{C\in\mathcal C}\{\log(e|C|)\}^{q/2}
 \le C_qm\{\log(ed/m)\}^{q/2}.
\]
If $\log(ed)=O(n^{1/3-\varepsilon})$ and $V_n\to0$, then
$m=o(n^{(q-2)/2})$ and hence $\log(em)=O(\log n)$. Thus the fixed
power of $\log(em)$ in \eqref{eq:group-sufficient} is dominated by
$n^\varepsilon$. If $m$ is bounded, those factors are constant.
This proves the corresponding conclusions. For $m=d$,
apply \cref{thm:polynomial-envelope}, since the restriction is vacuous
when every coordinate can be a separate group. When $m$ diverges, these arguments do not establish necessity of the
additional factor $\{\log(em)\}^{\nu_q}$.
\end{proof}

\subsection{An exact threshold for bounded-variation factor models}
\label{sec:factor}

In the representation below, the independent factors are uniform
variables on $[0,1]$. For each coordinate, the total variations of its
factor functions have a finite sum. A fixed bound on these sums implies
bounded coordinates. The upper bound combines the empirical-process
KMT coupling \citep{KMT1975,MAJOR1990} with Stieltjes integration,
concentration of the denominators, and Gaussian anti-concentration.
The resulting comparison uses the original random denominators and
the population covariance; no correlation gap is required.

For a variation bound that grows as $a^{-1/(q-2)}$, the dimension scale
in Corollary~\ref{cor:factor-threshold} involves
$a^{2/[3(q-2)]}$. For $q>3$, this improves on the power
$a^{4/[3(q-2)]}$ for the general correlation-gap class in
Theorem~\ref{thm:main-critical}. The improvement concerns the smaller
factor class and requires the additional assumptions stated below.

\begin{theorem}\label{thm:factor}
Let $X_1,\ldots,X_n$ be independent copies of a centered vector
$X\in\R^d$ with $\E X_j^2=1$. Suppose that
\begin{equation}\label{eq:factor-representation}
 X_j=\sum_{h\ge1}\left\{f_{jh}(U_h)-\int_0^1 f_{jh}(u)\,du\right\},
 \qquad
 \max_{j\le d}\sum_{h\ge1}\operatorname{TV}(f_{jh})\le B<\infty,
\end{equation}
where the $U_h$ are independent uniform random variables on $[0,1]$ and
all $f_{jh}$ have bounded variation. Finitely or countably many factors
are allowed. If $G\sim N_d(0,\E XX^\top)$, then, for $n,d\ge2$,
\begin{equation}\label{eq:factor-coupling}
 \distrect(T_n,G)
 \le 1\wedge
 \frac{CB\log(2nd)\sqrt{\log(2d)}}{\sqrt n},
\end{equation}
where $C$ is an absolute constant.
\end{theorem}

\begin{proof}
Put $v_{jh}=\operatorname{TV}(f_{jh})$. Each centered summand in
\eqref{eq:factor-representation} has absolute value at most $v_{jh}$, so the series in
\eqref{eq:factor-representation} converges absolutely and $|X_j|\le B$.
In particular, $B\ge1$. We use right-continuous representatives of the
functions; changing their values at the countably many discontinuities
does not change their distributions.

For each fixed $n$, apply the uniform empirical-process KMT coupling
independently to every factor. The required bound is given in
\citet[Theorem~3]{KMT1975}; see also the full-interval case of
\citet[Theorem~A]{MAJOR1990}. Let $\mathbb B_h$ be the coupled Brownian
bridge, and let $E_h$ be its sup-norm distance from the normalized
empirical process of the $h$th factor. For absolute $c,C>0$,
\[
 \Pp\{\sqrt n E_h>C\log n+x\}\le Ce^{-cx},\qquad x\ge0.
\]
Consequently, for a sufficiently small absolute $c_0>0$,
\begin{equation}\label{eq:factor-kmt-mgf}
 \E e^{c_0\sqrt n E_h}\le Cn^C.
\end{equation}
For countably many factors, take the product of the coupling spaces.

The covariance between $-\int_0^1\mathbb B_h(t)\,df_{jh}(t)$ and
the corresponding functional with $j$ replaced by $k$ is
$\Cov(f_{jh}(U_h),f_{kh}(U_h))$. Indeed, integrate the Brownian-bridge kernel
$\min(s,t)-st$ against the two finite signed Stieltjes measures. The same
kernel is the covariance of $\1_{\{U_h\le s\}}-s$ and
$\1_{\{U_h\le t\}}-t$. Finite total variation and boundedness of the kernel
justify these integrals.

Define $G_j=-\sum_{h\ge1}\int_0^1\mathbb B_h(t)\,df_{jh}(t)$.
The covariance calculation gives $\Cov G=\E XX^\top$. Stieltjes
integration by parts then gives, for $S=n^{-1/2}\sum_{i=1}^nX_i$,
\begin{equation}\label{eq:factor-path-coupling}
 |S_j-G_j|\le\sum_{h\ge1}v_{jh}E_h.
\end{equation}
In the countable case, each Gaussian summand has both its $L^1$ and $L^2$
norms bounded by $Cv_{jh}$. The series thus converges absolutely almost
surely and in $L^1$, and also converges in $L^2$. The covariance identity
and \eqref{eq:factor-path-coupling} follow by finite truncation.

Apply Jensen's inequality with weights $v_{jh}/B$, placing their unused
mass at zero. By \eqref{eq:factor-kmt-mgf},
\[
 \E\exp\{c_0\sqrt n|S_j-G_j|/B\}\le Cn^C.
\]
A union bound and Markov's inequality give
\begin{equation}\label{eq:factor-sum-coupling}
 \Pp\left\{\|S-G\|_\infty>
          \frac{CB\log(2nd)}{\sqrt n}\right\}\le\frac Cn.
\end{equation}
The factor coupling preserves dependence between coordinates, and
the union bound requires no coordinate independence.

It remains to control the denominators. We may assume that
$B\log(2nd)\sqrt{\log(2d)/n}$ is smaller than a sufficiently small absolute
constant, since otherwise \eqref{eq:factor-coupling} is trivial.
The inequalities $\E X_j^4\le B^2$ and $|X_j^2-1|\le CB^2$, followed by
Bernstein's inequality and a union bound, imply that, outside an event of
probability $C/n$,
\begin{equation}\label{eq:factor-squares}
 \max_{j\le d}\left|\frac1n\sum_{i=1}^nX_{ij}^2-1\right|
 \le C\left\{B\sqrt{\frac{\log(2nd)}n}
                +\frac{B^2\log(2nd)}n\right\}\le\frac12.
\end{equation}
The Gaussian union bound and \eqref{eq:factor-sum-coupling} give
$\|G\|_\infty+\|S\|_\infty\le C\sqrt{\log(2nd)}$ outside an event of
probability at most $C/n$. Since $x\mapsto x^{-1/2}$ is Lipschitz on $[1/2,3/2]$,
\[
 \|T_n-S\|_\infty
 \le C\sqrt{\log(2nd)}
       \left\{B\sqrt{\frac{\log(2nd)}n}
                   +\frac{B^2\log(2nd)}n\right\}
 \le\frac{CB\log(2nd)}{\sqrt n}.
\]
The last step uses $B\sqrt{\log(2nd)/n}\le c$. Combining this with
\eqref{eq:factor-sum-coupling} gives the same bound for
$\|T_n-G\|_\infty$, except on an event of probability $C/n$.
Apply Nazarov's inequality through
\eqref{eq:gaussian-coupling-transfer} to obtain
\eqref{eq:factor-coupling}. Unit marginal variances suffice for this
step, including when the covariance is singular
\citep{NAZ2003,NAZPROOF2017}.
\end{proof}

\subsubsection{A matching lower bound}
\label{sec:skew-boundary}

Fix $q>2$ and $K>1$, and set
\[
 D_+=n^{1/3}a^{2/[3(q-2)]}.
\]
In this subsection we assume
\begin{equation}\label{eq:skew-gap-assumptions}
 0<a_n<1\text{ eventually},\qquad a_n\to0,
 \qquad a_nD_+\to\infty.
\end{equation}
The condition $a_nD_+\to\infty$ ensures that the binomial
moderate-deviation expansion applies at this scale.

\begin{theorem}\label{thm:skew-block}
Under \eqref{eq:skew-gap-assumptions},
\begin{align*}
 \Del(d_n,a_n)\to0
   &\quad\Longrightarrow\quad \log(2d_n)=o(D_+),\\
 \liminf_n\frac{\log(2d_n)}{D_+}>0
   &\quad\Longrightarrow\quad \liminf_n\Del(d_n,a_n)>0,\\
 \frac{\log(2d_n)}{D_+}\to\infty
   &\quad\Longrightarrow\quad \Del(d_n,a_n)\to1.
\end{align*}
For every fixed $C_0>0$, the lower-bound examples use two independent
indicator factors per coordinate, with the sum of the absolute
coefficients at most $C_0a^{-1/(q-2)}$. The constants may depend on $C_0$.
\end{theorem}

\begin{proof}
Fix $C_0>0$. First take $M^2$ coordinates, where $M=\lfloor e^r\rfloor$ and
\begin{equation}\label{eq:skew-entropy-range}
 cD_+\le r\le D_+(aD_+)^{1/4}
\end{equation}
for some fixed $c>0$. Choose
$0<\omega<\min(1,C_0/2)$ with $1+\omega^{q-2}<K$, and put
\begin{equation}\label{eq:skew-parameters}
 \begin{gathered}
 s_0=\sqrt{1-a},\qquad B=\omega a^{-1/(q-2)},\qquad
 \pi(1-\pi)=a/B^2,\quad 0<\pi<1/2,\\
 u=\Phi^{-1}(1-M^{-1}),\qquad e_n=Bu^2/\sqrt n.
 \end{gathered}
\end{equation}
The smaller root for $\pi$ exists eventually. Since $\pi\sim a/B^2$ and
$D_+^3=na^{2/(q-2)}$,
\[
 \frac{D_+^2}{n\pi}\asymp\frac1{aD_+},\qquad
 \frac{B^2}{n}=\frac{\omega^2}{D_+^3}.
\]
Thus
\[
 \frac{r^2}{n\pi}\le\frac{C}{(aD_+)^{1/2}}\to0,
 \qquad
 \frac{B\sqrt r}{\sqrt n}
 \le\frac{C(aD_+)^{1/8}}{D_+}\to0.
\]
In particular,
\begin{equation}\label{eq:skew-scales}
 u\sim\sqrt{2r},\qquad
 e_nu\asymp(r/D_+)^{3/2},\qquad
 \frac{u^4}{n\pi}\to0,\qquad \frac{Bu}{\sqrt n}\to0.
\end{equation}
Also $u^4/n\to0$.

Let $U_{ib}$ be independent symmetric signs and let $I_{ibj}$ be mutually
independent Bernoulli variables with parameter $\pi$, independent of the
signs. Set
\begin{equation}\label{eq:skew-law}
 Y_{ibj}=B(I_{ibj}-\pi),\qquad X_{ibj}=s_0U_{ib}+Y_{ibj},
 \qquad 1\le b,j\le M.
\end{equation}
The coordinates are centered with variance one. Their correlations are
$1-a$ within a block and zero between blocks. Moreover, $B\pi\to0$ and
$\pi B^q\to\omega^{q-2}$, so conditioning on the Bernoulli variable gives
$\E|X_{ibj}|^q\to1+\omega^{q-2}<K$. The exact Gaussian reference is
\[
 G_{bj}=s_0Z_b+\sqrt a Z_{bj},
\]
with independent standard Gaussian factors.

We first locate the Gaussian maximum. Let $V=\max_{j\le M}Z_{bj}$ in one
block, and write $h(x)=\phi(x)/\barPhi(x)$. Monotonicity of the normal
hazard implies, for every real $y$,
\[
 M\barPhi(u+y)\le e^{-h(u)y}.
\]
Hence $\Pp(V-u>x)\le e^{-h(u)x}$ for $x\ge0$. With
$r_0=h(u)\sqrt a/s_0<h(u)$ eventually, tail integration gives
\begin{equation}\label{eq:skew-mgf}
 \E e^{r_0(V-u)}\le1+\frac{r_0}{h(u)-r_0}=1+o(1).
\end{equation}
Put $t_0=(s_0+\sqrt a)u$. For fixed $z\in\R$, the probability that a
single block exceeds $t_0+s_0z/u$, multiplied by $M$, tends to $e^{-z}$.
For the upper bound, condition on $V$ and use the hazard inequality and
\eqref{eq:skew-mgf}, noting $h(u)/u\to1$. For the lower bound, restrict to
$|V-u|\le A/u$, use the uniform normal tail ratio there, and then let
$A\to\infty$. This last passage follows from
$\Pp\{u(V-u)\le y\}=\Phi(u+y/u)^M\to\exp(-e^{-y})$.
Independence of the blocks gives
\begin{equation}\label{eq:skew-gaussian-limit}
 \Pp\{\max_{b,j}G_{bj}\le t_0+s_0z/u\}
 \longrightarrow F(z):=\exp(-e^{-z}).
\end{equation}
If $e_nu\to\infty$, the same conditioning and a union bound give
\begin{equation}\label{eq:skew-gaussian-far}
 \Pp\{\max_{b,j}G_{bj}>t_0+e_n/24\}
 \le e^{-h(u)e_n/(24s_0)}(1+o(1))\to0.
\end{equation}

Choose a block $b_*$ maximizing
$S_{U,b}=n^{-1/2}\sum_iU_{ib}$, with deterministic tie-breaking. The
symmetric-sign tail estimate in \cref{cor:sign-tails} applies because
$u^4/n\to0$, and yields
\begin{equation}\label{eq:skew-sign-max}
 \Pp\{u(\max_bS_{U,b}-u)\le z\}\to F(z).
\end{equation}
The selected block depends only on the signs and is therefore
independent of all Bernoulli variables $I_{ibj}$. Within that block,
choose the coordinate maximizing
$S_{Y,b_*j}=n^{-1/2}\sum_iY_{ib_*j}$. Set
$\sigma=\sqrt{n\pi(1-\pi)}=\sqrt{na}/B$ and apply
\cref{lem:binomial-tail} at $x=u+u^2/(12\sigma)$. Shifting the normal
threshold from $u$ to $x$ changes its log tail probability by
$-u^3/(12\sigma)+o(1)$, while the binomial cubic correction is
$(1-2\pi)u^3/(6\sigma)+o(1)$. Since
$u^3/\sigma=e_nu/\sqrt a\to\infty$, $M$ times this binomial tail tends to
infinity. At $x=2u$, it tends to zero, since the cubic correction is
$o(u^2)$ and the Gaussian exponent is $-2u^2$. Thus
\begin{equation}\label{eq:skew-private-max}
 \Pp\{\sqrt a u+e_n/12\le\max_jS_{Y,b_*j}\le2\sqrt a u\}\to1.
\end{equation}

To control the empirical denominators, use the identity
\[
 Y^2=a+B(1-2\pi)Y,\qquad
 X^2-1=B(1-2\pi)Y+2s_0UY.
\]
The expression for $X^2-1$ is centered, has variance at most
$C(B^2a+a)$, and has absolute value at most $CB^2$ eventually. Bernstein's inequality and a
union bound over $M^2$ coordinates give
\begin{equation}\label{eq:skew-denominator}
 \max_{b,j}\left|\frac1n\sum_iX_{ibj}^2-1\right|
 =O_P\left(\frac{B\sqrt a\,u}{\sqrt n}+\frac{B^2u^2}n\right),
 \qquad
 u\max_{b,j}\left|\frac1n\sum_iX_{ibj}^2-1\right|=o_P(e_n).
\end{equation}
Indeed, dividing the first bound times $u$ by $e_n$ gives
$O_P(\sqrt a+Bu/\sqrt n)=o_P(1)$. The selected numerator is $O_P(u)$ by
\eqref{eq:skew-sign-max} and \eqref{eq:skew-private-max}. Its original
self-normalized ratio therefore differs from that numerator by $o_P(e_n)$,
so, with probability tending to one,
\begin{equation}\label{eq:skew-sample-shift}
 \max_{b,j}T_{n,bj}\ge s_0\max_bS_{U,b}+\sqrt a u+e_n/13.
\end{equation}
If $\liminf e_nu>0$, choose $0<\eta<\liminf e_nu/13$. At $t_0$, the
cdf of the sample maximum has upper limit at most $F(-\eta)<F(0)$,
whereas the cdf of the Gaussian maximum tends to $F(0)$. If
$e_nu\to\infty$, the former cdf at $t_0+e_n/24$ tends to zero, and
\eqref{eq:skew-gaussian-far} shows that the latter tends to one.

For general $d$, take
\[
 r=\min\{\tfrac12\log d,\ D_+(aD_+)^{1/4}\}
\]
on a subsequence where $\log(2d)/D_+$ is bounded below. This satisfies
\eqref{eq:skew-entropy-range}; if the dimension ratio diverges, then
$r/D_+\to\infty$. There are $\lfloor e^r\rfloor^2\le d$ constructed
coordinates. Append independent symmetric signs to reach dimension $d$. This gives
both lower bounds, and a subsequence argument proves necessity.

Finally, represent the sign shared within each block by
$2(\1_{\{U\le1/2\}}-1/2)$ and the coordinate-specific Bernoulli term by
$\1_{\{V\le\pi\}}-\pi$, using independent uniforms. The total coefficient
magnitude is $2\sqrt{1-a}+\omega a^{-1/(q-2)}\le C_0a^{-1/(q-2)}$
eventually. The appended symmetric signs satisfy the same variation bound.
\end{proof}

\begin{corollary}\label{cor:factor-threshold}
Fix $q>2$, $K>1$, $C_0>0$, and assume \eqref{eq:skew-gap-assumptions}.
Let $\Delta^{\mathrm{BV}}_{n,q,K,C_0}(d,a)$ be the supremum of the rectangle
error over laws in the marginal-moment and correlation-gap class that
also admit \eqref{eq:factor-representation} with $B=C_0a^{-1/(q-2)}$.
Then
\[
 \Delta^{\mathrm{BV}}_{n,q,K,C_0}(d_n,a_n)\to0
 \quad\Longleftrightarrow\quad
 \log(2d_n)=o\{n^{1/3}a_n^{2/[3(q-2)]}\}.
\]
A positive lower limit of the dimension ratio gives a positive lower
limit of the error, and divergence of that ratio gives error tending to
one. The constants may depend on $C_0$.
\end{corollary}

\begin{proof}
Since $aD_+\to\infty$, eventually $a\ge D_+^{-1}$. Together with
$D_+^3=na^{2/(q-2)}$, this gives $D_+^{3+2/(q-2)}\ge n$ and hence
$D_+\gg\log n$. If $\log(2d)=o(D_+)$, then $\log(2nd)=o(D_+)$ and
$B\{\log(2nd)\}^{3/2}/\sqrt n\to0$. Apply \cref{thm:factor}.
The lower-bound examples in \cref{thm:skew-block} belong to this factor
class. They prove the two nonconvergence assertions and, by a
subsequence argument, necessity.
\end{proof}

\appendix
\section{Details for the inference remarks}\label{sec:inference-details}

\subsection{Transfer from Studentized to self-normalized statistics}

We use the notation of Section~\ref{sec:inference-scope}.
The identity $s_{n,j}^2=n^{-1}\sum_iX_{ij}^2-\bar X_j^2$
gives \eqref{eq:studentization-transform} on $E_n$.
Applied coordinatewise, $k_n$ is continuous and strictly increasing.
The inverse image of a coordinate rectangle under $k_n$ is
a coordinate rectangle or the empty set. Hence
\[
 \distrect(T_n,G)
 \le \distrect(U_n,G)+\Pp(E_n^c)
       +\distrect(k_n(G),G).
\]
To bound the last term, use $|k_n(x)-x|\le |x|^3/(2n)$.
The event $\max_j|G_j|\le\sqrt{2\log(2dn)}$ has probability
at least $1-1/n$ by a Gaussian union bound. On this event,
the coordinate displacement is at most
$C\log^{3/2}(edn)/n$. Gaussian rectangle anti-concentration
\citep{NAZPROOF2017,CCK2015} therefore gives
\[
 \distrect(k_n(G),G)\le C\frac{\log^2(edn)}n,
\]
which proves \eqref{eq:studentization-transfer}. The argument
uses unit marginal variances and permits singular $R_P$.

\subsection{Uniform marginal calibration}

Let $\delta=\min(q-2,1)$. The moment assumption gives
$\E|X_j|^{2+\delta}\le K^{(2+\delta)/q}$ uniformly in the
coordinate and the underlying law. Apply
\citet[Theorem~2.3]{JSW2003} to both $X_{ij}$ and $-X_{ij}$.
Uniformly over $P\in\mathcal P_{d_n,q}(K,0)$ and $j$,
\[
 \Pp(|T_{n,j}|>x)
 \le 2\barPhi(x)
       \{1+C_{q,K}(1+x)^{2+\delta}n^{-\delta/2}\}
\]
for $0\le x\le c_{q,K}n^{\delta/[2(2+\delta)]}$.
The growth condition in \eqref{eq:marginal-calibration} places
$c_n=\Phi^{-1}(1-\alpha_n/(2d_n))$ in this range and gives
$(1+c_n)^{2+\delta}n^{-\delta/2}=o(1)$.
Summing the scalar bounds over $j$ proves
\eqref{eq:marginal-calibration}, with the same $o(1)$ for every
law in the class. This deduction uses no independence between
coordinates.

\par\medskip
\noindent\textbf{AI Statement.} During the preparation of this manuscript,
the authors used ChatGPT and Codex (OpenAI) for language editing and
LaTeX formatting, and to assist with the exposition and organization
of proofs. These tools were also used to check references, notation,
cross-references, and mathematical derivations. The authors reviewed and
revised the manuscript, verified the mathematical arguments, and take
full responsibility for its content.

\bibliographystyle{apalike}
\bibliography{snclt}
\end{document}